\documentclass[10pt]{amsproc}

\usepackage{amsmath}
\usepackage{amssymb}
\usepackage[all]{xy}
\usepackage{mathrsfs}
\usepackage{graphicx}
\usepackage{epsf}

\DeclareMathAlphabet{\scr}{U}{eus}{m}{n}

\newcommand\N{{\mathbb N}}
\newcommand\Z{{\mathbb Z}}
\newcommand\Q{{\mathbb Q}}
\newcommand\F{{\mathbb F}}

\newcommand\A{{\mathbb A}}
\newcommand\G{{\mathbb G}}
\newcommand\C{{\mathbb C}}
\newcommand\R{{\mathbb R}}

\newcommand\cA{{\mathcal A}}

\newcommand\cF{{\mathcal F}}

\newcommand\cH{{\mathcal H}}

\newcommand\cO{{\mathcal O}}

\newcommand\cV{{\mathcal V}}
\newcommand\cL{{\mathcal L}}

\newcommand\Mloc{{\rm\bf M}^{\rm loc}}

\newcommand\xt{[\hspace{-.12em}[ t ]\hspace{-.12em} ] } 
\newcommand\xT{(\!( t )\!)}

\newcommand\Ql{\overline{\mathbb Q}_\ell}

\newcommand\Fl{\mathcal Fl}

\newcommand\tens\otimes

\newcommand\lto{\longrightarrow}

\newcommand{\Grass}{\mathop{\rm Grass}\nolimits}

\newcommand\Gal{\mathop{\rm Gal}}

\newtheorem{theorem}{Theorem}[section]
\newtheorem{lemma}[theorem]{Lemma}

\theoremstyle{definition}
\newtheorem{definition}[theorem]{Definition}

\newtheorem{xca}[theorem]{Exercise}
\newtheorem{prop}[theorem]{Proposition}
\theoremstyle{remark}
\newtheorem{rem}[theorem]{Remark}
\newtheorem{cor}[theorem]{Corollary}
\newtheorem{conj}[theorem]{Conjecture}

\numberwithin{equation}{subsection}

\begin{document}
\ifx\href\undefined\else\hypersetup{linktocpage=true}\fi 

\title{Introduction to Shimura varieties with \\ bad reduction of parahoric type}
\author{Thomas J. Haines}
\thanks{Research partially supported by NSF grant DMS 0303605.}
\address{Mathematics Department, University of Maryland, College Park, MD 20742-4015}  
\email{tjh@math.umd.edu}


\maketitle

\pagestyle{myheadings}
\markboth{THOMAS J. HAINES}{SHIMURA VARIETIES WITH PARAHORIC LEVEL STRUCTURE}

\section{Introduction}

This survey article is intended to introduce the reader to several important concepts relating to Shimura varieties with parahoric level structure at $p$.  The main tool is the Rapoport-Zink local model \cite{RZ}, which plays an important role in several aspects of the theory.  We discuss local models attached to general linear and symplectic groups, and we illustrate their relation to Shimura varieties in two examples: the  {\em simple} or {\em ``fake'' unitary} Shimura varieties with parahoric level structure, and the Siegel modular varieties with $\Gamma_0(p)$-level structure.  In addition, we present some applications of local models to questions of flatness, stratifications of the special fiber, and the determination of the semi-simple local zeta functions for simple Shimura varieties. 

There are several good references for material of this sort that already exist in the literature.  This survey  has a great deal of overlap with two articles of Rapoport: \cite{R1} and \cite{R2}.   A main goal of this paper is simply to make more explicit some of the ideas expressed very abstractly in those papers.   Hopefully it will shed some new light on the earlier seminal works of Rapoport-Zink \cite{RZ82}, and Zink \cite{Z}.  This article is also closely related to some recent work of H. Reimann \cite{Re1}, \cite{Re2}, \cite{Re3}.  

Good general introductions to aspects of the Langlands program which might be consulted while reading parts of this report are those of Blasius-Rogawski \cite{BR}, and T. Wedhorn \cite{W2}.

Several very important developments have taken place in the theory of Shimura varieties with bad reduction, which are completely ignored in this report.  In particular, we mention the book of Harris-Taylor \cite{HT} which uses bad reduction of Shimura varieties to prove the local Langlands conjecture for ${\rm GL}_n(\Q_p)$, and the recent work of L. Fargues and E. Mantovan \cite{FM}.

Most of the results stated here are well-known by now (although scattered around the literature, with differing systems of conventions).  However, the author took this opportunity to present a few new results (and some new proofs of old results).  For example, there is the proof of the non-emptiness of the Kottwitz-Rapoport strata in any connected component of the Siegel modular and ``fake'' unitary Shimura varieties with Iwahori-level structure (Lemmas \ref{Genestier_surjectivity}, \ref{connectedness_of_fibers}), some foundational relations between Newton strata, Kottwitz-Rapoport strata, and affine Deligne-Lusztig varieties (Prop. \ref{KR_vs_Newton_strata_via_DL}), and the verification of the conjectural non-emptiness of the basic locus in the ``fake'' unitary case (Cor. \ref{basic_locus_nonempty}).  The main new proofs relate to topological flatness of local models attached to Iwahori subgroups of unramified groups (see $\S \ref{flatness}$) and to the description of the nonsingular locus of Shimura varieties with Iwahori-level structure (see $\S \ref{smooth_locus}$).  Finally, some of the results explained here (especially in $\S \ref{ss_for_simple}$) are background material necessary for the author's as yet unpublished joint work with B.C. Ng\^{o} \cite{HN3}.

I am grateful to U. G\"{o}rtz, R. Kottwitz, B.C. Ng\^{o}, G. Pappas, and M. Rapoport for all they have taught me about this subject over the years.  I thank them for their various comments and suggestions on an early version of this article.  
Also, I heartily thank U. G\"{o}rtz for providing the figures.  Finally, I thank the Clay Mathematics Institute for sponsoring the June 2003 {\em Summer School on Harmonic Analysis, the Trace Formula, and Shimura Varieties}, which provided the opportunity for me to write this survey article.

\small
\tableofcontents
\normalsize

\section{Notation} \label{notation}
 
\subsection{Some field-theoretic notation}  Fix a rational prime $p$.  We let $F$ denote a non-Archimedean local field of residual characteristic $p$, with ring of integers ${\mathcal O}$.  Let ${\mathfrak p} \subset \mathcal O$ denote the maximal ideal, and fix a uniformizer $\pi \in {\mathfrak p}$.  The residue field $\mathcal O/ \mathfrak p$ has cardinality $q$, a power of $p$. 

We will fix an algebraic closure $k$ of the finite field ${\mathbb F}_p$.  The Galois group ${\rm Gal}(k/{\mathbb F}_p)$ has a canonical generator (the Frobenius automorphism), given by $\sigma(x) = x^p$.  For each positive integer $r$, we denote by $k_r$ the fixed field of $\sigma^r$.  Let $W(k_r)$ (resp. $W(k)$) denote the ring of Witt vectors of $k_r$ (resp. $k$), with fraction field $L_r$ (resp. $L$).  We also use the symbol $\sigma$ to denote the Frobenius automorphism of $L$ induced by that on $k$. 

We fix throughout a rational prime $\ell \neq p$, and a choice of algebraic closure $\mathbb Q_\ell \subset \Ql$.

\subsection{Some group-theoretic notation}

The symbol ${\bf G}$ will always denote a connected reductive group over $\Q$ (sometimes defined over $\Z$).  Unless otherwise indicated, $G$ will denote the base-change ${\bf G}_F$, where $F$ is a suitable local field (usually, $G = {\bf G}_{\Q_p}$).

Now let $G$ denote a connected reductive $F$-group.  Fix once and for all a Borel subgroup $B$ and a maximal torus $T$ contained in $B$. We will usually assume $G$ is split over $F$, in which case we can even assume $G,B$ and $T$ are defined and split over the ring $\mathcal O$.  For ${\rm GL}_n$ or ${\rm GSp}_{2n}$, $T$ will denote the usual ``diagonal'' torus, and $B$ will denote the ``upper triangular'' Borel subgroup.

We will often refer to ``standard'' parahoric and ``standard'' Iwahori subgroups.  For the group $G = {\rm GL}_n$ (resp. $G = {\rm GSp}_{2n}$), the ``standard'' hyperspecial maximal compact subgroup of $G(F)$ will be the subgroup $G(\mathcal O)$.  The ``standard'' Iwahori subgroup will be inverse image of $B(\cO/\mathfrak p)$ under the reduction modulo $\mathfrak p$ homomorphism
$$
G(\mathcal O) \rightarrow G(\mathcal O/ \mathfrak p).
$$
A ``standard'' parahoric will be defined similarly as the inverse image of a standard (= upper triangular) parabolic subgroup.

For ${\rm GL}_n(F)$, the standard Iwahori is the subgroup stabilizing the standard lattice chain (defined in $\S \ref{Iwahori_subgroups}$).  The standard parahorics are stabilizers of standard partial lattice chains.  Similar remarks apply to the group ${\rm GSp}_{2n}(F)$.  The symbols $I$ or $I_r$ or $K^{\bf a}_p$ will always denote a standard Iwahori subgroup of $G(F)$ defined in terms of our fixed choices of $B \supset T$, and $G(\cO)$ as above (for some local field $F$).  Often (but not always) $K$ or $K_r$ or $K_p^0$ will denote our fixed hyperspecial maximal compact subgroup $G(\cO)$.  

We have the associated spherical Hecke algebra $\cH_K := C_c(K\backslash G(F)/K)$, a convolution algebra of $\C$-valued (or $\Ql$-valued) functions on $G(F)$ where convolution is defined using the measure giving $K$ volume 1.  Similarly, $\cH_I := C_c(I\backslash G(F)/I)$ is a convolution algebra using the measure giving $I$ volume 1.  For a compact open subset $U \subset G(F)$, ${\mathbb I}_U$ denotes the characteristic function of the set $U$.

The {\em extended affine Weyl group} of $G(F)$ is defined as the group
$\widetilde{W} = N_{G(F)}T/T(\cO)$.  The map $X_*(T) \rightarrow T(F)/T(\cO)$ given by $\lambda \mapsto \lambda(\pi)$ is an isomorphism of abelian groups.  The finite Weyl group $W_0 := N_{G(F)}T/T(F)$ can be identified with $N_{G(\cO)}T/T(\cO)$, by choosing representatives of $N_{G(F)}T/T(F)$ in $G(\cO)$.   Thus we have a canonical isomorphism
$$
\widetilde{W} = X_*(T) \rtimes W_0.
$$
We will denote elements in this group typically by the notation $t_\nu w$ (for $\nu \in X_*(T)$ and $w \in W_0$).

Our choice of $B \supset T$ determines a unique {\em opposite} Borel subgroup $\bar{B}$ such that $B \cap \bar{B} = T$.  We have a notion of $B$-positive (resp. $\bar{B}$-positive) root $\alpha$ and coroot $\alpha^\vee$.  Also, the group $W_0$ is a Coxeter group generated by the simple reflections $s_\alpha$ in the vector space $X_*(T) \otimes \R$ through the walls fixed by the $B$-positive (or $\bar{B}$-positive) simple roots $\alpha$.  Let $w_0$ denote the unique element of $W_0$ having greatest length with respect to this Coxeter system.

We will often need to consider $\widetilde{W}$ as a subset of $G(F)$.  We choose the following conventions.  For each $w \in W_0$, we fix once and for all a lift in the group $N_{G(\cO)}T$.  We identify each 
$\nu \in X_*(T)$ with the element $\nu(\pi) \in T(F) \subset G(F)$.

Let ${\bf a}$ denote the alcove in the building of $G(F)$ which is fixed by the Iwahori $I$, or equivalently, the unique alcove in the apartment corresponding to $T$ which is contained in the $\bar{B}$-positive (i.e., the $B$-{\em negative}) Weyl chamber, whose closure contains the origin 
(the vertex fixed by the maximal compact subgroup 
$G(\cO)$) \footnote{This choice of base alcove results from our convention of embedding $X_*(T) \hookrightarrow G(F)$ by the rule $\lambda \mapsto \lambda(\pi)$; to see this, consider how vectors spanning the standard periodic lattice chain in $\S \ref{Iwahori_subgroups}$ are identified with vectors in $X_*(T) \otimes \R$.}.

The group $\widetilde{W}$ permutes the set of affine roots $\alpha + k$ ($\alpha$ a root, $k \in \Z$) (viewed as affine linear functions on $X_*(T) \otimes \R$), and hence permutes (transitively) the set of alcoves.  Let $\Omega$ denote the subgroup which stabilizes the base alcove ${\bf a}$.  Then we have a semi-direct product
$$
\widetilde{W} = W_{\rm aff} \rtimes \Omega,
$$
where $W_{\rm aff}$ (the {\em affine Weyl group}) is the Coxeter group generated by the reflections $S_{\rm aff}$ through the walls of ${\bf a}$.  In the case where $G$ is an almost simple group of rank $l$, with simple $B$-positive roots $\alpha_1, \dots, \alpha_l$, then $S_{\rm aff}$ consists of the $l$ simple reflections $s_i = s_{\alpha_i}$ generating $W_0$, along with one more simple affine reflection $s_0 = t_{-\widetilde{\alpha}^\vee} s_{\widetilde{\alpha}}$, where $\widetilde{\alpha}$ is the highest $B$-positive root.  

The Coxeter system $(W_{\rm aff}, S_{\rm aff})$ determines a length function $\ell$ and a Bruhat order $\leq$ on $W_{\rm aff}$, which extend naturally to $\widetilde{W}$: for $x_i \in W_{\rm aff}$ and $\sigma_i \in \Omega$ ($i = 1,2$), we define $x_1 \sigma_1 \leq x_2 \sigma_2 $ in $\widetilde{W}$ if and only if $\sigma_1 = \sigma_2$ and $x_1 \leq x_2$ in $W_{\rm aff}$.  Similarly, we set $\ell(x_1 \sigma_1) = \ell(x_1)$.

By the Bruhat-Tits decomposition, the inclusion $\widetilde{W} \hookrightarrow G(F)$ induces a bijection
$$
\widetilde{W} = I \backslash G(F) / I.
$$

In the function-field case (e.g., $F = \F_p((t))$), the affine flag variety ${\mathcal Fl} = G(F)/I$ is naturally an ind-scheme, and the closures of the $I$-orbits ${\mathcal Fl}_w : = IwI/I$ are determined by the Bruhat order on $\widetilde{W}$:
$$
{\mathcal Fl}_x \subset \overline{{\mathcal Fl}_y} \Longleftrightarrow  x \leq y.
$$

Similar statements hold for the affine Grassmannian, ${\rm Grass} = G(F)/G(\cO)$.  Now the $G(\cO)$-orbits ${\mathcal Q}_\lambda := G(\cO)\lambda G(\cO)/G(\cO)$ are given (using the Cartan decomposition) by the $B$-dominant 
coweights $X_+(T)$:
$$
X_+(T) \leftrightarrow G(\cO) \backslash G(F) / G(\cO).
$$
By definition, $\lambda$ is $B$-dominant if $\langle \alpha, \lambda \rangle \geq 0$ for all $B$-positive roots $\alpha$.  Here $\langle \cdot, \, \cdot \rangle: X^*(T) \times X_*(T) \rightarrow \Z$ is the canonical duality pairing.

The closure relations in ${\rm Grass}$ are given by the partial order on $B$-dominant coweights $\lambda$ and $\mu$:
$$
{\mathcal Q}_\lambda \subset \overline{{\mathcal Q}_\mu} \Longleftrightarrow \lambda \preceq \mu,
$$
where $\lambda \prec \mu$ means that $\mu - \lambda$ is a sum of $B$-positive coroots.

Unless otherwise stated, a {\em dominant} coweight $\lambda \in X_*(T)$ will always mean one that is $B$-dominant.

For the group $G = {\rm GL}_n$ or ${\rm GSp}_{2n}$, there is a $\Z_p$-ind-scheme $M$ which is a deformation of the affine Grassmannian ${\rm Grass}_{{\mathbb Q}_p}$ to the affine flag variety $\Fl_{\F_p}$ for the underlying $p$-adic group $G$ (see \cite{HN1}, and Remark $\ref{deformation_remark}$).  A very similar deformation ${\rm Fl}_X$ over a smooth curve $X$ (due to Beilinson) exists for any group $G$ in the function field setting, and has been extensively studied by Gaitsgory \cite{Ga}.  For any dominant coweight $\lambda \in X_+(T)$, the symbol $M_\lambda$ will always denote the $\Z_p$-scheme which is the scheme-theoretic closure in $M$ of ${\mathcal Q}_\lambda \subset {\rm Grass}_{\Q_p}$. 

\subsection{Duality notation}

If $A$ is any abelian scheme over a scheme $S$, we denote the dual abelian scheme by $\widehat{A}$.  The existence of $\widehat{A}$ over an arbitrary base is a delicate matter; see \cite{CF}, $\S$I.1.

If $M$ is a module over a ring $R$, we denote the dual module by $M^\vee = {\rm Hom}_R(M,R)$.  Similar notation applies to quasi-coherent $\cO_S$-modules over a scheme $S$.

If $G$ is a connected reductive group over a local field $F$, then the Langlands dual (over $\C$ or $\Ql$) will be denoted $\widehat{G}$.  The Langlands $L$-group is the semi-direct product $^LG = \widehat{G} \rtimes W_F$, where $W_F$ is the Weil-group of $F$.

\subsection{Miscellaneous notation}

We will use the following abbreviation for elements of ${R}^n$ (here $R$ can be any set): let $a_1, \dots, a_r$ be a sequence of positive integers whose sum is $n$.  A vector of the form $(x_1, \dots, x_1, x_2, \dots, x_2, \dots, x_r,\dots, x_r)$, where for $i = 1, \dots, r$, the element $x_i$ is repeated $a_i$ times, will be denoted by 
$$
(x_1^{a_1},x_2^{a_2}, \dots, x_r^{a_r}).
$$

We will denote by $\A$ the adeles of $\Q$, by $\A_f$ the finite adeles, and by 
$\A^p_f$ the finite adeles away from $p$ (with the exception of two instances, 
where $\A$ denotes affine space!).

\section{Iwahori and parahoric subgroups} \label{Iwahori_subgroups}

\subsection{Stabilizers of periodic lattice chains}

We discuss the definitions for the groups ${\rm GL}_n$ and ${\rm GSp}_{2n}$.

\subsubsection{The linear case}
For each $i \in \{1, \dots, n \}$, let $e_i$ denote the $i$-th standard vector $(0^{i-1},1,0^{n-i})$ in $F^n$, and let $\Lambda_i \subset F^n$ denote the free $\mathcal O$-module with basis $\pi^{-1}e_1, \dots, \pi^{-1}e_i, e_{i+1}, \dots, e_n$.  We consider the diagram
$$
\Lambda_\bullet: ~ \Lambda_0 \lto \Lambda_1 \lto \cdots \lto \Lambda_{n-1} \lto \pi^{-1}\Lambda_0,
$$
where the morphisms are inclusions.  The lattice chains $\pi^n\Lambda_\bullet$ ($n \in \Z$) fit together to form an infinite complete lattice chain $\Lambda_i$, ($i \in \Z$).  If we identify each $\Lambda_i$ with $\mathcal O^n$, then the diagram above becomes
$$\xymatrix{
\mathcal O^n \ar[r]^{m_1} & \mathcal O^n \ar[r]^{m_2} & \cdots  \ar[r]^{m_{n-1}} & \mathcal O^n \ar[r]^{m_n} & \mathcal O^n,}
$$
where $m_i$ is the morphism given by the diagonal matrix 
$$
m_i = {\rm diag}(1, \dots, \pi, \dots, 1),
$$
the element $\pi$ appearing in the $i$th place.  We define the {\em (standard) Iwahori subgroup} $I$ of ${\rm GL}_n(F)$  by
$$
I = \bigcap_{i} \,\, {\rm Stab}_{{\rm GL}_n(F)}\, ( \Lambda_i).
$$
Similarly, for any non-empty subset $J \subset \{0,1,\dots, n-1 \}$, we define the {\em parahoric subgroup} of ${\rm GL}_n(F)$ corresponding to the subset $J$ by
$$
{\mathcal P}_J = \bigcap_{i \in J} \,\, {\rm Stab}_{{\rm GL}_n(F)} \, (\Lambda_i).
$$
Note that $\mathcal P_J$ is a compact open subgroup of ${\rm GL}_n(F)$, and that $\mathcal P_{\{0\}} = {\rm GL}_n(\mathcal O)$ is a hyperspecial maximal compact subgroup, in the terminology of Bruhat-Tits, cf. \cite{Tits}.

\subsubsection{The symplectic case}
The definitions for the group of symplectic similitudes ${\rm GSp}_{2n}$ are similar.  We define this group using the alternating matrix 
$$
\label{matrix}
\tilde{I} = \begin{bmatrix} 0 & \tilde{I}_n \\ -\tilde{I}_n & 0 \end{bmatrix}, 
$$
where $\tilde{I}_n$ denotes the $n \times n$ matrix with 1 along the anti-diagonal and 0 elsewhere.  Let $(x, y) := x^t \, \tilde{I} y$ denote the corresponding alternating pairing on $F^{2n}$.  For an $\mathcal O$-lattice $\Lambda \subset F^{2n}$, we define $\Lambda^\perp := \{ x \in F^{2n} ~ | ~ (x,y) \in \mathcal O \,\,\, \mbox{for all $y \in \Lambda$} \}$.  The lattice $\Lambda_0$ is self-dual (i.e., $\Lambda_0^\perp = \Lambda_0$).  Consider the infinite lattice chain in $F^{2n}$
$$
\cdots \lto \Lambda_{-2n} = \pi \Lambda_0 \lto \cdots \lto \Lambda_{-1} \lto \Lambda_0 \lto \cdots \lto \Lambda_{2n} = \pi^{-1}\Lambda_0 \lto \cdots
$$
We have $\Lambda^\perp_i = \Lambda_{-i}$ for all $i \in \Z$.  Now we define the {\em (standard) Iwahori subgroup} $I$ of ${\rm GSp}_{2n}(F)$ by
$$
I = \bigcap_{i} \,\, {\rm Stab}_{{\rm GSp}_{2n}(F)} \, (\Lambda_i).
$$
For any non-empty subset $J \subset \{ -n, \dots, -1, 0, 1, \dots, n \}$ such that $i \in J \Leftrightarrow -i \in J$, we define the {\em parahoric subgroup} corresponding to $J$ by 
$$
{\mathcal P}_J = \bigcap_{i \in J} \,\, {\rm Stab}_{{\rm GSp}_{2n}(F)} \, (\Lambda_i).
$$

\subsection{Bruhat-Tits group schemes} \label{BT_group_schemes}

In Bruhat-Tits theory, parahoric groups are defined as the groups ${\mathfrak G}_{\Delta_J}^0(\mathcal O)$, where ${\mathfrak G}^0_{\Delta_J}$ is the neutral component of a group scheme ${\mathfrak G}_{\Delta_J}$, defined and smooth over ${\mathcal O}$, which has generic fiber the $F$-group $G$, and whose ${\mathcal O}$-points are the subgroup of $G(F)$ fixing the facet ${\Delta}_J$ of the Bruhat-Tits building corresponding to the set $J$; see \cite{BT2}, p. 356.  By \cite{Tits}, 3.4.1 (see also \cite{BT2}, 1.7) we can characterize the group scheme ${\mathfrak  G}_{\Delta_J}$ as follows: it is the unique ${\mathcal O}$-group scheme ${\mathcal P}$ satisfying the following three properties:
\begin{enumerate}
\item ${\mathcal P}$ is smooth and affine over ${\mathcal O}$;
\item The generic fiber ${\mathcal P}_F$ is $G_F$;
\item For any unramified extension $F'$ of $F$, letting ${\mathcal O}_{F'} \subset F'$ denote ring of integers, the group ${\mathcal P}({\mathcal O}_{F'}) \subset G(F')$ is the subgroup of elements which fix the facet $\Delta_J$ in the Bruhat-Tits building of $G_{F'}$.
\end{enumerate}

Let us show how automorphism groups of periodic lattice chains $\Lambda_\bullet$ give a concrete realization of the Bruhat-Tits parahoric group schemes, in the ${\rm GL}_n$ and ${\rm GSp}_{2n}$ cases.  We will discuss the Iwahori subgroups of ${\rm GL}_n$ in some detail, leaving for the reader the obvious generalizations to parahoric subgroups of ${\rm GL}_n$ (and ${\rm GSp}_{2n}$).

For any $\mathcal O$-algebra $R$, we may consider the diagram $\Lambda_{\bullet,R} = \Lambda_\bullet \otimes_{\mathcal O} R$, and we may define the ${\mathcal O}$-group scheme ${\rm Aut}$ whose $R$-points are the isomorphisms of the diagram $\Lambda_{\bullet,R}$.  More precisely, an element of ${\rm Aut}(R)$ is an $n$-tuple of $R$-linear automorphisms
$$
(g_0, \dots, g_{n-1}) \in {\rm Aut}(\Lambda_{0,R}) \times \cdots \times {\rm Aut}(\Lambda_{n-1,R})
$$
such that the following diagram commutes
$$
\xymatrix{
\Lambda_{0,R} \ar[r] \ar[d]^{g_0} & \cdots \ar[r] & \Lambda_{n-1,R} \ar[r] \ar[d]^{g_{n-1}} & \Lambda_{n,R} \ar[d]^{g_0} \\ \Lambda_{0,R} \ar[r] & \cdots \ar[r] & \Lambda_{n-1,R} \ar[r] & \Lambda_{n,R}.} 
$$

The group functor ${\rm Aut}$ is obviously represented by an affine group scheme, also denoted ${\rm Aut}$.  Further, it is not hard to see that ${\rm Aut}$ is formally smooth, hence smooth, over ${\mathcal O}$.  To show this, one has to check the lifting criterion for formal smoothness: if $R$ is an ${\mathcal O}$-algebra containing a nilpotent ideal ${\mathcal J} \subset R$, then any automorphism of $\Lambda_\bullet \otimes_{\mathcal O} R/{\mathcal J}$ can be lifted to an automorphism of $\Lambda_{\bullet} \otimes_{\mathcal O} R$.  This is proved on page 135 of \cite{RZ}.  Thus ${\rm Aut}$ satisfies condition (1) above.

Alcoves in the Bruhat-Tits building for ${\rm GL}_n$ over $F$ can be described as complete periodic $\mathcal O$-lattice chains in $F^{n}$ 
$$
\cdots \lto {\mathcal L}_0 \lto {\mathcal L}_1 \lto \cdots \lto {\mathcal L}_{n-1} \lto {\mathcal L}_n = \pi^{-1}{\mathcal L}_0 \lto \cdots
$$
where the arrows are inclusions.  We can regard $\Lambda_\bullet$ as the base alcove in this building.
It is clear that since $\pi$ is invertible over the generic fiber $F$, the automorphism $g_0$ determines the other $g_i$'s over $F$ and so ${\rm Aut}_F = {\rm GL}_n$.   By construction we have ${\rm Aut}({\mathcal O}) = {\rm Stab}_{{\rm GL}_n(F)} \, (\Lambda_\bullet)$.  This is unchanged if we replace $F$ by an unramified extension $F'$.  It follows that ${\rm Aut}$ satisfies conditions (2) and (3) above.  Thus, by uniqueness, ${\rm Aut} = {\mathfrak G}_{\Delta_J}$ for $J = \{0,\dots, n-1 \}$.

Further, one can check that the special fiber ${\rm Aut}_{k}$ is an extension of the Borel subgroup $B_k$ by a connected unipotent group over $k$; hence the special fiber is connected.  It follows that ${\rm Aut}$ is a connected group scheme (cf. \cite{BT2}, 1.2.12).  So in this case ${\rm Aut} = {\mathfrak G}_{\Delta_J} = {\mathfrak G}^0_{\Delta_J}$.  We conclude that ${\rm Aut}({\mathcal O})$ is the Bruhat-Tits Iwahori subgroup fixing the base alcove $\Lambda_\bullet$.

\begin{xca} 1) Check the lifting criterion which shows ${\rm Aut}$ is formally smooth directly for the case $n=2$, by explicit calculations with $2 \times 2$ matrices.

\noindent 2)  By identifying ${\rm Aut}({\mathcal O})$ with its image in ${\rm GL}_n(F)$ under the inclusion $g_\bullet \mapsto g_0$, show that the Iwahori subgroup is the preimage of $B(k)$ under the canonical surjection ${\rm GL}_n(\mathcal O) \rightarrow {\rm GL}_n(k)$.

\noindent 3)  Prove that ${\rm Aut}_k \rightarrow {\rm Aut}(\Lambda_{0,k})$, $g_\bullet \mapsto g_0$, has image $B_k$, and kernel a connected unipotent group.

\end{xca}

\section{Local models} \label{local_models}

Given a certain triple $(G,\mu, K_p)$ consisting of a $\Z_p$-group $G$, a minuscule coweight $\mu$ for $G$, and a parahoric subgroup $K_p \subset G(\Q_p)$, one may construct a projective $\Z_p$-scheme ${\bf M}^{\rm loc}$ which (\'{e}tale locally) models the singularities found in the special fiber of a certain Shimura variety with parahoric-level structure at $p$.  The advantage of ${\bf M}^{\rm loc}$ is that it is defined in terms of linear algebra and is therefore easier to study than the 
Shimura variety itself.  These schemes are called ``local models'', or sometimes ``Rapoport-Zink local models''; the most general treatment is given in \cite{RZ}, but in special cases they were also investigated in \cite{DP} and \cite{deJ}.

In this section we recall the definitions of local models associated to ${\rm GL}_n$ and ${\rm GSp}_{2n}$.  For simplicity, we limit the discussion to models for Iwahori-level structure.  In each case, the local model is naturally associated to a dominant minuscule coweight $\mu$, which we shall always mention.  In fact, it turns out that if the Shimura data give rise to $(G,\mu)$, then the Rapoport-Zink local model ${\bf M}^{\rm loc}$ (for Iwahori-level structure) can be identified with the space $M_{-w_0\mu}$ mentioned in $\S \ref{notation}$.  See Remark \ref{deformation_remark} below.

\subsection{Linear case}  

We use the notation $\Lambda_\bullet$ from section \ref{Iwahori_subgroups} to denote the ``standard'' lattice chain over ${\mathcal O} = \Z_p$ here.  

Fix an integer $d$ with $1 \leq d \leq n-1$.  We define the scheme ${\bf M}^{\rm loc}$ by defining its $R$-points for any $\Z_p$-algebra $R$ as the set of isomorphism classes of commutative diagrams 
$$
\xymatrix{
\Lambda_{0,R} \ar[r] &\cdots \ar[r] &  \Lambda_{n-1,R} \ar[r] &  \Lambda_{n,R} \ar[r]^{p} & \Lambda_{0,R} \\   
{\mathcal F}_0 \ar[u] \ar[r] & \cdots \ar[r] & {\mathcal F}_{n-1} \ar[u] \ar[r] & {\mathcal F}_{n} \ar[u] 
\ar[r]^{p} & {\mathcal F}_0 \ar[u],}
$$
where the vertical arrows are inclusions, and each ${\mathcal F}_i$ is an $R$-submodule of $\Lambda_{i,R}$ which is Zariski-locally on $R$ a direct factor of corank $d$.  It turns out that this is identical to the space $M_{-w_0\mu}$ of $\S \ref{notation}$, where $\mu = (0^{n-d},(-1)^{d})$.  
It is clear that ${\bf M}^{\rm loc}$ is represented by a closed subscheme of a product of Grassmannians over $\Z_p$, hence it is a projective $\Z_p$-scheme.  One can also formulate the moduli problem using quotients of rank $d$ instead of submodules of corank $d$.

Another way to formulate the same moduli problem which is sometimes useful (see \cite{HN1}) is given by  
adding an indeterminate $t$ to the story (following a suggestion of G. Laumon).  We replace the ``standard'' 
lattice chain term $\Lambda_{i, R}$ with ${\mathcal V}_{i,R} := \alpha^{-i}R[t]^n$, where $\alpha$ is the $n \times n$ matrix
$$
\begin{pmatrix}
0 & 1 & &  \\ & \ddots & \ddots & \\ & & 0 & 1 \\ t + p & & & 0 \end{pmatrix} \in {\rm GL}_n(R[t,t^{-1},(t + p)^{-1}]).$$

One can identify ${\bf M}^{\rm loc}(R)$ with the set of chains 
$$
{\mathcal L}_\bullet = ({\mathcal L}_0 \subset {\mathcal L}_1 \subset \cdots \subset {\mathcal L}_n = (t + p)^{-1} {\mathcal L}_0)$$ 
of $R[t]$-submodules of $R[t,t^{-1}, (t+ p)^{-1}]^n$ satisfying the following properties
\begin{enumerate}
\item for all $i = 0, \dots, n-1$, we have $t{\mathcal V}_{i,R} \subset {\mathcal L}_i \subset {\mathcal V}_{i,R}$;
\item as an $R$-module, ${\mathcal L}_i/t{\mathcal V}_{i,R}$ is locally a direct factor of 
${\mathcal V}_{i,R}/t{\mathcal V}_{i,R}$ of corank $d$.
\end{enumerate}

\begin{rem} \label{deformation_remark}
With the second definition, it is easy to see that the geometric generic fiber ${\bf M}^{\rm loc}_{\overline{\Q}_p}$ is contained in the affine Grassmannian ${\rm GL}_n(\overline{\Q}_p\xT)/{\rm GL}_n(\overline{\Q}_p\xt) $, and the geometric special fiber $\Mloc_{\overline{\F}_p}$ is contained in the affine flag variety 
${\rm GL}_n(\overline{\F}_p\xT)/I_{\overline{\F}_p}$, where $I_{\overline{\F}_p} := {\rm Aut}(\overline{\F}_p\xt)$ is the Iwahori 
subgroup of ${\rm GL}_n(\overline{\F}_p\xt)$ corresponding to the upper triangular Borel subgroup $B \subset {\rm GL}_n$.  
Moreover, it is possible to view $\Mloc$ as a piece of a $\Z_p$-ind-scheme $M$ which forms a deformation of the affine Grassmannian to the affine flag variety over the base ${\rm Spec}(\Z_p)$, in analogy with Beilinson's deformation ${\rm Fl}_X$ over a smooth $\F_p$-curve $X$ in the function field case (cf. \cite{Ga}, \cite{HN1}).  
In fact, letting $e_0$ denote the base point in the affine Grassmannian ${\rm GL}_n(\Q_p((t)))/{\rm GL}_n(\Q_p[[t]])$, and for $\lambda \in X_+(T)$, letting ${\mathcal Q}_\lambda$ denote the ${\rm GL}_n(\Q_p[[t]])$-orbit of the point $\lambda(t)e_0$, it turns out that $\Mloc$ coincides with the scheme-theoretic closure $M_{-w_0\mu}$ of ${\mathcal Q}_{-w_0\mu}$ taken in the model $M$.  A similar statement holds in the symplectic case.  

The identification ${\bf M}^{\rm loc} = M_{-w_0\mu}$ is explained in $\S \ref{KR_strata}$.  
It is a consequence of the determinant condition and the ``homology'' definition of our local models; also the flatness of ${\bf M}^{\rm loc}$ (see $\S \ref{flatness}$) plays a role.
\end{rem}

\subsection{The symplectic case}

Recall that for our group ${\rm GSp}_{2n} = {\rm GSp}(V, (\cdot, \cdot))$, we have an identification of $X_*(T)$ with the lattice $\{ (a_1, \dots, a_n,b_n,\dots, b_1) \in \Z^{2n} ~ | ~ \exists \,\, c \in \Z, \,\, a_i + b_i = c, \forall i \}$.  The Shimura coweight that arises here has the form $\mu = (0^n, (-1)^n)$.

For the group ${\rm GSp}_{2n}$, the symbol $\Lambda_\bullet$ now denotes the self-dual $\Z_p$-lattice chain in $\Q_p^{2n}$, discussed in section \ref{Iwahori_subgroups} in the context of ${\rm GSp}_{2n}$.  Let $(\cdot, \cdot)$ denote the alternating pairing on $\Z_p^{2n}$ discussed in that section, and let the dual $\Lambda^\perp$ of a lattice $\Lambda$ be defined using $(\cdot,\cdot)$.  As above, there are (at least) two equivalent ways to define $\Mloc(R)$ for a $\Z_p$-algebra $R$.  
We define $\Mloc(R)$ to be the set of isomorphism classes of diagrams
$$
\xymatrix{
\Lambda_{0,R} \ar[r]  & \cdots \ar[r] & \Lambda_{n-1,R} \ar[r] & \Lambda_{n,R} \\ 
\cF_0 \ar[u] \ar[r] & \cdots  \ar[r] & \cF_{n-1} \ar[u] \ar[r] & \cF_n \ar[u],}
$$
where the vertical arrows are inclusions of $R$-submodules with the following properties:
\begin{enumerate}
\item for $i = 0, \dots, n$, Zariski-locally on $R$ the submodule $\cF_i$ is a direct factor of $\Lambda_{i,R}$ of corank $n$;
\item $\cF_0$ is isotropic with respect to $(\cdot,\cdot)$ and $\cF_n$ is isotropic with respect to $p(\cdot, \cdot)$.  
\end{enumerate}

As in the linear case, this can be described in a way more transparently connected to affine flag varieties.  In this case, ${\mathcal V}_{i,R}$ has the same meaning as in the linear case, except that the ambient space is now $R[t,t^{-1},(t + p)^{-1}]^{2n}$.  We may describe $\Mloc(R)$ as the set of chains 
$$
{\mathcal L}_\bullet = ({\mathcal L}_0 \subset {\mathcal L}_1 \subset \cdots \subset {\mathcal L}_n)
$$
of $R[t]$-submodules of $R[t,t^{-1},(t+ p)^{-1}]^{2n}$ satisfying the following properties
\begin{enumerate}
\item for $i = 0, 1, \dots, n$, $t{\mathcal V}_{i,R} \subset {\mathcal L}_i \subset {\mathcal V}_{i,R}$;
\item as $R$-modules, ${\mathcal L}_i/t{\mathcal V}_{i,R}$ is locally a direct factor of 
${\mathcal V}_{i,R}/t{\mathcal V}_{i,R}$ of corank $n$;
\item ${\mathcal L}_0$ is self-dual with respect to $t^{-1}(\cdot,\cdot)$, and ${\mathcal L}_n$ is self-dual with respect to $t^{-1}(t + p)(\cdot,\cdot)$.
\end{enumerate}

\subsection{Generic and special fibers} \label{generic_and_special_fibers}

In the linear case with $\mu = (0^{n-d},(-1)^{d})$, the generic fiber of $M_{-w_0\mu}$ is the Grassmannian ${\rm Gr}(d,n)$ of $d$-planes in $\Q_p^n$.  In the symplectic case with $\mu = (0^n,(-1)^n)$, the generic fiber of $M_{-w_0\mu}$ is the Grassmannian of isotropic $n$-planes in $\Q_p^{2n}$ with respect to the alternating pairing $(\cdot,\cdot)$.

In each case, the special fiber is a union of finitely many Iwahori-orbits $I  w I /I $ in the affine flag variety, indexed by elements $w$ in the extended affine Weyl group $\widetilde{W}$ for ${\rm GL}_n$ 
(resp. ${\rm GSp}_{2n}$) ranging over the so-called $-w_0\mu$-{\em admissible} subset ${\rm Adm}(-w_0\mu)$, \cite{KR}.  Let $\lambda \in X_+(T)$.  Then by definition
$$
{\rm Adm}(\lambda) = \{ w \in \widetilde{W} ~ | ~ \exists \,\, \nu \in W\lambda, \,\,\,\mbox{such that 
$w \leq t_\nu$}  \}.
$$
Here, $W\lambda$ is the set of conjugates of $\lambda$ under the action of the finite Weyl group $W_0$, and $t_\nu \in \widetilde{W}$ is the translation element corresponding to $\nu$, and $\leq$ denotes the Bruhat order on $\widetilde{W}$.  Actually (see $\S \ref{construction_of_KR_stratification}$), the set that arises naturally from the moduli problem is the $-w_0\mu$-{\em permissible} subset ${\rm Perm}(-w_0\mu) \subset \widetilde{W}$ from \cite{KR}.  Let us recall the definition of this set, following loc. cit.  Let $\lambda \in X_+(T)$ and suppose $t_\lambda \in W_{\rm aff} \tau$, for $\tau \in \Omega$.  Then ${\rm Perm}(\lambda)$ consists of the elements $x \in W_{\rm aff}\tau$ such that $x(a) - a \in {\rm Conv}(\lambda)$ for every vertex $a \in {\bf a}$.  Here ${\rm Conv}(\lambda)$ denotes the convex hull of $W\lambda$ in $X_*(T) \otimes \R$.

The strata in the special fiber of $\Mloc = M_{-w_0\mu}$ are naturally indexed by the set ${\rm Perm}(-w_0\mu)$, which agrees with ${\rm Adm}(-w_0\mu)$ by the following non-trivial combinatorial theorem due to Kottwitz and Rapoport.

\begin{theorem} [\cite{KR}; see also \cite{HN2}] \label{Perm=Adm}
For every minuscule coweight $\lambda$ of either ${\rm GL}_n$ or ${\rm GSp}_{2n}$, we have the equality$$
{\rm Perm}(\lambda) = {\rm Adm}(\lambda).
$$
\end{theorem}

Using the well-known correspondence between elements in the affine Weyl group and the set of alcoves in the standard apartment of the Bruhat-Tits building, one can ``draw'' pictures of ${\rm Adm}(\mu)$ for low-rank groups.  Figures 1 and 2 depict this set for $G = {\rm GL}_3$, $\mu= (-1,0,0)$, and $G = {\rm GSp}_4$, $\mu = (-1,-1,0,0)$ \footnote{Note that in Figure 2 there is an alcove of length one contained in the Bruhat-closure of three distinct translations.  This already tells us something about the singularities: the special fiber of the Siegel variety for ${\rm GSp}_4$ is not a union of divisors with normal crossings; see 
$\S \ref{KR_strata}$.}.  Actually, we draw the image of ${\rm Adm}(\mu)$ in the apartment for ${\rm PGL}_3$ (resp. ${\rm PGSp}_4$); the base alcove is labeled by $\tau$. 
\begin{figure}[hbt]
\epsfxsize=4 in
\centerline{\epsfbox{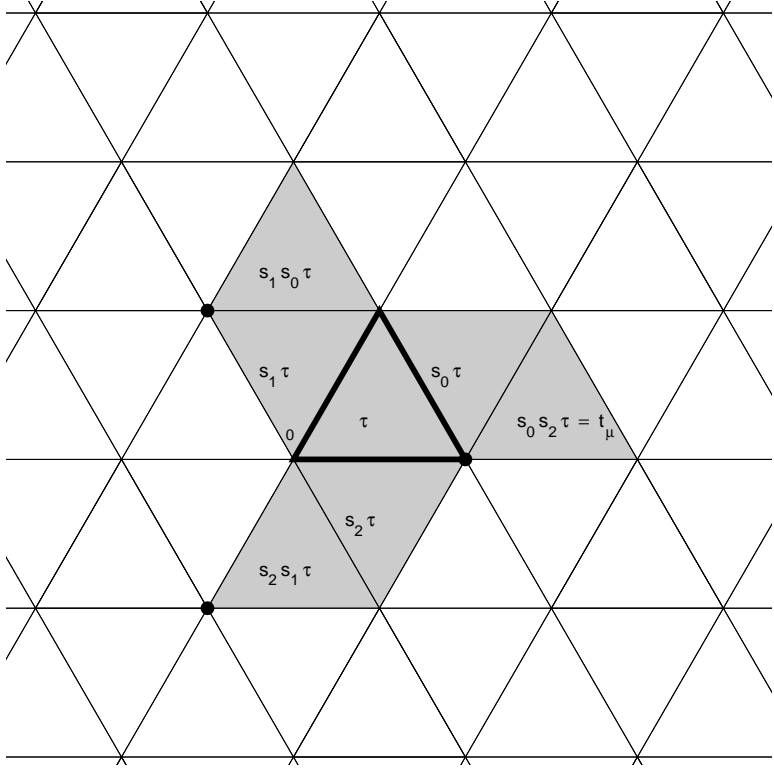}}
\caption{The admissible alcoves ${\rm Adm}(\mu)$ for ${\rm GL}_3, \,\,\mu = (-1,0,0)$. 
The base alcove is labeled by $\tau$.}
\end{figure}

\medskip
\begin{figure}[hbt]
\epsfxsize=4 in
\centerline{\epsfbox{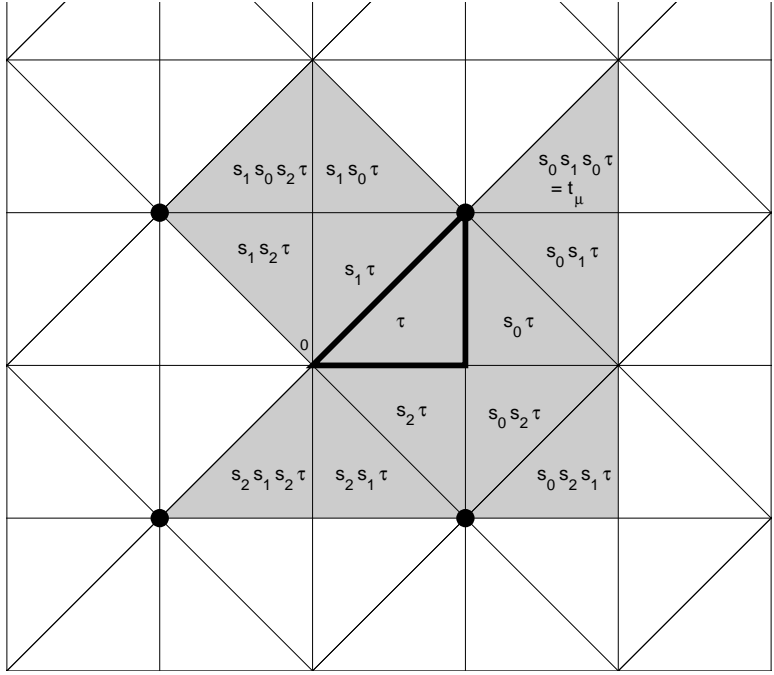}}
\caption{The admissible alcoves ${\rm Adm}(\mu)$ for ${\rm GSp}_4, \,\,\mu = (-1,-1,0,0)$.}
\end{figure}

\subsection{Computing the singularities in the special fiber of $\Mloc$}

In certain cases, the singularities in $\Mloc_{\bar{\F}_p}$ can be analyzed directly by writing down equations.  As the simplest example of how this is done, we analyze the local model for ${\rm GL}_2$, $\mu = (0,-1)$.  For a $\Z_p$-algebra $R$, we are looking at the set of pairs $(\cF_0,\cF_1)$ of locally free rank 1 $R$-submodules of $R^2$ such that the following diagram commutes
$$
\xymatrix{
R \oplus R \ar[r]^{\begin{tiny} \begin{bmatrix} p & 0 \\ 0 & 1 \end{bmatrix} \end{tiny}} & R \oplus R \ar[r]^{\begin{tiny} \begin{bmatrix} 1 & 0 \\ 0 & p \end{bmatrix} \end{tiny}} & R \oplus R \\
\cF_0 \ar[u] \ar[r] & \cF_1 \ar[u] \ar[r] & \cF_0 \ar[u] .
}
$$
Obviously this functor is represented by a certain closed subscheme of ${\mathbb P}^1_{\Z_p} \times {\mathbb P}^1_{\Z_p}$.  Locally around a fixed point $(\cF_0,\cF_1) \in {\mathbb P}^1(R) \times {\mathbb P}^1(R)$ we choose coordinates such that $\cF_0$ is represented by the homogeneous column vector $[ 1 : x]^t$ and $\cF_1$ by the vector $[y:1]^t$, for $x,y \in R$.  We see that $(\cF_0,\cF_1) \in \Mloc(R)$ if and only if 
$$
xy = p,
$$
so $\Mloc$ is locally the same as ${\rm Spec}(\Z_p[X,Y]/(XY-p)$, the usual deformation of ${\mathbb A}^1_{\Q_p}$ to a union of two ${\mathbb A}_{\F_p}^1$'s which intersect transversally at a point.  Indeed, $\Mloc$ is globally this kind of deformation:
\begin{itemize}
\item In the generic fiber, the matrices are invertible and so $\cF_0$ uniquely determines $\cF_1$; thus $\Mloc_{\Q_p} \cong {\mathbb P}^1_{\Q_p}$;
\item In the special fiber, $p=0$ and one can check that $\Mloc_{\F_p}$ is the union of the closures of two Iwahori-orbits in the affine flag variety ${\rm GL}_2(\F_p \xT)/I_{\F_p}$, each of dimension 1, which meet in a point.  Thus $\Mloc_{\F_p}$ is the union of two ${\mathbb P}^1_{\F_p}$'s meeting in a point. 
\end{itemize}

We refer to the work of U. G\"{o}rtz for many more complicated calculations of this kind: \cite{Go1}, 
\cite{Go2}, \cite{Go3}, \cite{Go4}.

\section{Some PEL Shimura varieties with parahoric level structure at $p$} \label{PEL_parahoric_section}

\subsection{PEL-type data} \label{PEL_type_data} 

Given a Shimura datum $({\bf G},\{h\},{\bf K})$ one can construct a Shimura variety $Sh({\bf G},h)_{\bf K}$ which has a canonical model over the reflex field ${\bf E}$, a number field determined by the datum.  We write $G$ for the $p$-adic group ${\bf G}_{\Q_p}$.  Let us assume that the compact open subgroup ${\bf K} \subset {\bf G}(\A_f)$ is of the form ${\bf K} = K^pK_p$, where $K^p \subset {\bf G}(\A^p_f)$ is a sufficiently small compact open subgroup, and $K_p \subset G(\Q_p)$ is a parahoric subgroup.  

Let us fix once and for all embeddings $\overline{\Q} \hookrightarrow \C$, and $\overline{\Q} \hookrightarrow \overline{\Q}_p$.  We denote by ${\mathfrak p}$ the corresponding place of ${\bf E}$ over $p$ and by $E = {\bf E}_{\mathfrak p}$ the completion of ${\bf E}$ at ${\mathfrak p}$.  

If the Shimura datum comes from PEL-type data, then it is possible to define a moduli problem (in terms of chains of abelian varieties with additional structure) over the ring $\cO_E$.  This moduli problem is representable by a 
quasi-projective $\cO_E$-scheme whose generic fiber is the base-change to $E$ of the initial Shimura variety $Sh({\bf G},h)_{\bf K}$ (or at least a finite union of Shimura varieties, one of which is the canonical model $Sh({\bf G},h)_{\bf K}$).  This is done in great generality in Chapter 6 of \cite{RZ}.  Our aim in this section is only to make somewhat more explicit the definitions in loc. cit., in two very special cases attached to the linear and symplectic groups.

First, let us recall briefly PEL-type data.  Let $B$ denote a finite-dimensional semi-simple $\Q$-algebra with 
positive involution $\iota$.  Let $V \neq 0$ be a finitely-generated left $B$-module, and let $( \cdot, \cdot )$ be a non-degenerate alternating form $V \times V \rightarrow \Q$ on the underlying $\Q$-vector space, such that $(bv, w ) = ( v, b^\iota w)$, for $b \in B$, $v,w \in V$.  The form $( \cdot, \cdot )$ determines a ``transpose'' involution on ${\rm End}(V)$, denoted by $*$ (so viewing the left-action of $b$ as an element of ${\rm End}(V)$, we have $b^\iota = b^*$).  We denote by ${\bf G}$ the $\Q$-group whose points in a $\Q$-algebra $R$ are 
$$
\{ g \in {\rm GL}_{B \otimes R}(V \otimes R) ~ | ~ g^* g = c(g) \in R^\times \}.
$$
We assume ${\bf G}$ is a {\em connected} reductive group; this means we are excluding the orthogonal case.  Consider the $\R$-algebra $C := {\rm End}_B(V) \otimes \R$.  We let $h_0 : \C \rightarrow C$ denote an $\R$-algebra homomorphism satisfying $h_0(\overline{z}) = h_0(z)^*$, for $z \in \C$.  We fix a choice of $i = \sqrt{-1}$ in $\C$ once and for all, and we assume the symmetric bilinear form $( \cdot\, , \, h_0(i)\,\cdot ) : V_\R \times V_\R \rightarrow \R$ is positive definite.  Let $h$ denote the inverse of the restriction of $h_0$ to $\C^\times$.  Then $h$ induces an algebraic homomorphism 
$$
h : \C^\times \rightarrow {\bf G}(\R)
$$
of real groups which defines on $V_{\R}$ a Hodge structure of type $(1,0) + (0,1)$ (in the terminology of \cite{Del2}, section 1) and which satisfies the usual Riemann conditions with respect to $( \cdot, \cdot )$ (see \cite{Ko92}, Lemma 4.1).  
For any choice of (sufficiently small) compact open subgroup ${\bf K}$, the data $({\bf G}, h, {\bf K})$ determine a (smooth) Shimura variety over a reflex field ${\bf E}$; cf. \cite{Del}.  

We recall that $h$ gives rise to a minuscule coweight
$$
\mu := \mu_h : \G_{m,\C} \rightarrow {\bf G}_\C
$$
as follows: the complexification of the real group $\C^\times$ is the torus $\C^\times \times \C^\times$, the factors being indexed by the two $\R$-algebra automorphisms of $\C$; we assume the first factor corresponds to the identity and the second to complex conjugation.  Then we define 
$$\mu(z) := h_\C(z,1).$$  
By definition of Shimura data, the homomorphism $h: \C^\times \rightarrow {\bf G}_\R$ is only specified up to ${\bf G}(\R)$-conjugation, and therefore $\mu$ is only well-defined up to ${\bf G}(\C)$-conjugation.  However, this conjugacy class is at least defined over the reflex field ${\bf E}$ (in fact we {\em define} ${\bf E}$ as the field of definition of the conjugacy class of $\mu$).  Via our choice of field embeddings $\overline{\C} \hookleftarrow \overline{\Q} \hookrightarrow \overline{\Q}_p$, we get a well-defined conjugacy class of minuscule coweights
$$
\mu : \G_{m,\overline{\Q}_p} \rightarrow G_{\overline{\Q}_p},
$$
which is defined over $E$.  

The argument of \cite{Ko84}, Lemma (1.1.3) shows that $E$ is contained in any subfield of $\overline{\Q}_p$ which splits $G$.  Therefore, when $G$ is split over $\Q_p$ (the case of interest in this report), it follows that $E = \Q_p$ and the conjugacy class of $\mu$ contains a $\Q_p$-rational and $B$-dominant element, usually denoted also by the symbol $\mu$.  It is this same $\mu$ which was mentioned in the definitions of local models in section \ref{local_models}. 

For use in the definition to follow, 
we decompose the $B_\C$-module $V_\C$ as $V_\C = V_1 \oplus V_2$, where $h_0(z)$ acts by $z$ on $V_1$ and by $\overline{z}$ on $V_2$, for $z \in \C$.   Our conventions imply that $\mu(z)$ acts by $z^{-1}$ on $V_1$ and by $1$ on $V_2$ ($z \in \C^\times$).  We choose $E' \subset \overline{\Q}_p$ a finite extension field $E' \supset E$ over which this decomposition is defined:
$$
V_{E'} = V_1 \oplus V_2.
$$ 
(We are implicitly using the diagram $\C \hookleftarrow \overline{\Q} \hookrightarrow \overline{\Q}_p$ to make sense of this.)

Recall that we are interested in defining an $\cO_E$-integral model for $Sh({\bf G},h)_{\bf K}$ in the case where $K_p \subset {\bf G}(\Q_p)$ is a parahoric (more specifically, an Iwahori) subgroup.

To define an integral model over $\cO_E$, we need to specify certain additional data.  We suppose $\cO_B$ is a $\Z_{(p)}$-order in $B$ whose $p$-adic completion $\cO_B \otimes \Z_p$ is a maximal order in $B_{\Q_p}$, stable under the involution $\iota$.  Using the terminology of \cite{RZ}, 6.2, we assume we are given a self-dual multichain ${\mathcal L}$ of $\cO_B \otimes \Z_p$-lattices in $V_{\Q_p}$ (the notion of multichain ${\mathcal L}$ is a generalization of the lattice chain $\Lambda_\bullet$ appearing in section \ref{local_models}; specifying ${\mathcal L}$ is equivalent to specifying a parahoric subgroup, namely $K_p := {\rm Aut}({\mathcal L})$, of ${\bf G}(\Q_p)$).  We can then give the definition of a model $Sh_{K^p}$ that depends on the above data and the choice of a small compact open subgroup $K^p$ \footnote{In the sequel, we sometimes drop the subscript $K^p$ on $Sh_{K^p}$, or replace it with the subscript $K_p$, depending on whether $K^p$, or $K_p$ (or both) is understood.}.

\begin{definition} \label{general_definition}
 A point of the functor $Sh_{K^p}$ with values in the $\cO_E$-scheme $S$ is given by the following set of data up to isomorphism \footnote{We say $\{ A_{\Lambda} \}$ is isomorphic to $\{ A'_\Lambda \}$ if there is a compatible family of prime-to-$p$ isogenies $A_\Lambda \rightarrow A'_\Lambda$ 
which preserve all the structures.}.
\begin{enumerate}
\item An ${\mathcal L}$-set of abelian $S$-schemes $A = \{ A_{\Lambda} \}$,  $\,\, \Lambda \in {\mathcal L}$, compatibly endowed with an action of $\cO_B$:
$$
i : \cO_B \otimes \Z_{(p)} \rightarrow {\rm End}(A) \otimes \Z_{(p)};
$$
\item A $\Q$-homogeneous principal polarization $\overline{\lambda}$ of the ${\mathcal L}$-set $A$;
\item A $K^p$-level structure
$$
\bar{\eta}: V \otimes \A^p_f \cong H_1(A, \A^p_f) \,\,\, \mbox{mod $\, K^p$}
$$
that respects the bilinear forms on both sides up to a scalar in $(\A^p_f)^\times$, and commutes with the 
$B   = \cO_B \otimes \Q$-actions.
\end{enumerate}
We impose the condition that under 
$$
i : \cO_B \otimes \Z_{(p)} \rightarrow {\rm End}(A) \otimes \Z_{(p)},
$$
we have $i(b^\iota)  = \lambda^{-1} \circ (i(b))^\vee \circ \lambda$; in other words, $i$ intertwines $\iota$ and the Rosati involution on ${\rm End}(A) \otimes \Z_{(p)}$ determined by $\overline{\lambda}$.    
In addition, we impose the following determinant condition: for each $b \in \cO_B$ and $\Lambda \in {\mathcal L}$:  
$$
{\rm det}_{\cO_S}(b, {\rm Lie}(A_\Lambda)) = {\rm det}_{E'}(b, V_1).
$$
\end{definition}

We will not explain all the notions entering this definition; we refer to loc. cit., Chapter 6 as well as \cite{Ko92}, section 5, for complete details.  However, in the simple examples we make explicit below, these notions will be made concrete and their importance will be highlighted.  For example, an ${\mathcal L}$-set of abelian varieties $\{A_\Lambda \}$ comes with a family of ``periodicity isomorphisms''
$$
\theta_a : A^a_{\Lambda} \rightarrow A_{a\Lambda},
$$
see \cite{RZ}, Def. 6.5, and we will describe these explicitly in the examples to follow.

Note that one can see from this definition why some of the conditions on PEL data are imposed.  For example, since the Rosati involution is always positive (see \cite{Mu}, section 21), we see that the involution $\iota$ on $B$ must be positive for the moduli problem to be non-empty.

\subsection{Some ``fake'' unitary Shimura varieties} \label{simple_setup}

This section concerns the so-called ``simple'' or ``fake unitary'' Shimura varieties investigated by Kottwitz in \cite{Ko92b}.  They are indeed ``simple'' in the sense that they are compact Shimura varieties for which there are no problems due to endoscopy (see loc. cit.).

Kottwitz made assumptions ensuring that the local group ${\bf G}_{\Q_p}$ be unramified, and that the level structure at $p$ be given by a hyperspecial maximal compact subgroup.  Here we will work in a situation where ${\bf G}_{\Q_p}$ is {\em split}, but we only impose parahoric-level structure at $p$.  For simplicity, we explain only the case $F_0 = \Q$ (notation of loc. cit.).

\subsubsection{The group-theoretic set-up}

Let $F$ be an imaginary quadratic extension of $\Q$, and let $(D,*)$ be a division algebra with center $F$, of dimension $n^2$ over $F$, together with an involution $*$ which induces on $F$ the non-trivial element of ${\rm Gal}(F/\Q)$.  Let ${\bf G}$ be the $\Q$-group whose points in a commutative $\Q$-algebra $R$ are
$$
\{ x \in D \otimes_\Q R ~ | ~ x^*x \in R^\times \}.
$$
The map $x \mapsto x^* x$ is a homomorphism of $\Q$-groups ${\bf G} \rightarrow \G_m$ whose kernel ${\bf G}_0$ is an inner form of a unitary group over $\Q$ associated to $F/\Q$.  Let us suppose we are given an $\R$-algebra homomorphism 
$$
h_0 : \C \rightarrow D \otimes_\Q \R
$$
such that $h_0(z)^* = h_0(\overline{z})$ and the involution $x \mapsto h_0(i)^{-1} x^* h_0(i)$ is positive. 

Given the data $(D,*,h_0)$ above, we want to explain how to find the PEL-data $(B,\iota, V, (\cdot, \cdot),h_0)$ used in the definition of the scheme $Sh_{K^p}$.

Let $B = D^{opp}$ and let $V = D$ be viewed as a left $B$-module, free of rank 1, using right multiplications.  
Thus we can identify $C := {\rm End}_B(V)$ with $D$ (left multiplications).  For $h_0 : \C \rightarrow C_\R$ we use the homomorphism $h_0 : \C \rightarrow D \otimes_\Q \R$ we are given.

Next, one can show that there exist elements $\xi \in D^\times$ such that $\xi^* = -\xi$ and the involution $x \mapsto \xi x^* \xi^{-1}$ is positive.  To see this, note that the Skolem-Noether theorem implies that the involutions of the second type on $D$ are precisely the maps of the form $x \mapsto bx^*b^{-1}$, for $b \in D^\times$ such that $b(b^*)^{-1}$ lies in the center $F$.  Since positive involutions of the second kind exist (see \cite{Mu}, p. 201-2), for some such $b$ the involution $x \mapsto bx^*b^{-1}$ is positive.  We have 
$N_{F/\Q}(b(b^*)^{-1}) = 1$, so by Hilbert's Theorem 90, we may alter any such $b$ by an element in $F^\times$ so that $b^* = b$.  There exists $\epsilon \in F^\times$ such that $\epsilon^* = -\epsilon$.  We then put $\xi = \epsilon b$.

We define the positive involution $\iota$ by $x^\iota := \xi x^* \xi^{-1}$, for $x \in B = D^{opp}$.

Now we define the non-degenerate alternating pairing $(\cdot, \cdot) : D \times D \rightarrow \Q$ by
$$
(x,y) = tr_{D/\Q}(x\xi y^*).
$$
It is clear that $(bx,y) = (x,b^\iota y)$ for any $b \in B = D^{opp}$, remembering that the left action of $b$ is right multiplication by $b$.  We also have $(h_0(z)\, x,y) = (x, h_0(\overline{z}) \, y)$, since $h_0(z) \in D$ acts by left 
multiplication on $D$.

Finally, we claim that $(\cdot \, , \, h_0(i) \, \cdot)$ is always positive or negative definite; thus we can always arrange for it to be positive definite by replacing $\xi$ with $-\xi$ if necessary.   To prove the definiteness, choose an isomorphism
$$
D \otimes_\Q \R ~ \widetilde{\rightarrow} ~ M_n(\C)
$$
such that the positive involution $x \mapsto x^\iota$ goes over to the standard positive involution $X \mapsto \overline{X}^t$ on $M_n(\C)$.  Let $H \in M_n(\C)$ be the image of $\xi h_0(i)^{-1}$ under this isomorphism, so that the symmetric pairing $\langle x,y \rangle = (x \, , \,h_0(i) y)$ goes over to the pairing$$
\langle X,Y \rangle  = tr_{M_n(\C)/\R}(X \,  \overline{Y}^t \, H).
$$  
We conclude by invoking the following exercise for the reader. 

\begin{xca}  The matrix $H$ is Hermitian and either positive or negative definite.  If positive definite, we then 
have $tr(X \, \overline{X}^t \, H) > 0$ whenever $X \neq 0$.
\end{xca}
(Hint: use the argument of \cite{Mu}, p. 200.)

\subsubsection{The minuscule coweight $\mu$} \label{mu}
How is the minuscule coweight $\mu$ described in terms of the above data?  Recall our decomposition
$$
D_\C = V_1 \oplus V_2.
$$
The homomorphism $h_0$ makes $D_\C$ into a $\C \otimes_\R \C$-module.  Of course
\begin{align*}
\C \otimes_\R \C ~ &\widetilde{\rightarrow} ~ \C \times \C \\
z_1 \otimes z_2 &\mapsto (z_1z_2, \overline{z}_1 z_2), 
\end{align*}
which induces the above decomposition of $D_\C\,\,$ ($h_0(z_1 \otimes 1)$ acts by $z_1$ on $V_1$ and by $\overline{z}_1$ on $V_2$).  

The factors $V_1$, $V_2$ are stable under right multiplications of 
$$D_\C = D \otimes_{F,\nu} \C \times D \otimes _{F,\nu^*} \C,$$ 
where $\nu,\nu^* : F \hookrightarrow \C$ are the two embeddings. (We may assume our fixed choice $\overline{\Q} \hookrightarrow \C$ extends $\nu$.)  Also $h_0(z \otimes 1)$ is the endomorphism given by left multiplication by a certain element of $D_\C$.  We can choose an isomorphism 
$$
D \otimes_{F,\nu} \C \times D \otimes_{F,\nu^*} \C \cong M_n(\C) \times M_n(\C)$$
such that $h_0(z \otimes 1)$ can be written explicitly as 
$$
h_0(z \otimes 1) = {\rm diag}(\overline{z}^{n-d}, z^{d}) \times {\rm diag}(z^{n-d}, \overline{z}^{d}),
$$
for some integer $d$, $0 \leq d \leq n$.  (One can then identify $V_1$ resp. $V_2$ as the span of certain columns of the two matrices.)  We know that $\mu(z) = h_{\C}(z,1)$ acts by $z^{-1}$ on $V_1$ and by $1$ on $V_2$.  Hence we can identify $\mu(z)$ as
$$
\mu(z) = {\rm diag}(1^{n-d}, (z^{-1})^{d}) \times {\rm diag}((z^{-1})^{n-d}, 1^{d}).
$$ 
We may label this by $(0^{n-d}, (-1)^{d}) \in \Z^n$, via the usual identification applied to the first factor.  

Here is another way to interpret the number $d$.  Let $W$ (resp. $W^*$) be the (unique up to isomorphism, $n$-dimensional) simple right-module for $D \otimes_{F,\nu} \C$ (resp. $D \otimes_{F,\nu^*} \C$).  Then as right $D_\C$-modules we have 
$$V_1 = W^{d} \oplus (W^*)^{n-d}, \,\,\,\, \mbox{resp. $V_2 = W^{n-d} \oplus (W^*)^d$}.$$

Finally, let us remark that if we choose the identification of $D \otimes_\Q \R = D \otimes_{F, \nu} \C$ with $M_n(\C)$ in such a way that the positive involution $x \mapsto h_0(i)^{-1}x^* h_0(i)$ goes over to $X \mapsto \overline{X}^t$, then we get an isomorphism
$$
G(\R) ~ \cong ~ {\rm GU}(d,n-d).
$$
(See also \cite{Ko92b}, section 1.)

In applications, it is sometimes necessary to prescribe the value of $d$ ahead of time (with the additional constraint that $1 \leq d \leq n-1$).  However, it can be a delicate matter to arrange things so that a prescribed value of $d$ is achieved.  To see how this is done for the case of $d = 1$, see \cite{HT}, Lemma I.7.1.

\subsubsection{Assumptions on $p$ and integral data} \label{integral_data}

We first make some assumptions on the prime $p$ \footnote{It would make sense to include in our discussion another case, where $p$ remains inert in $F$, and where the group ${\bf G}_{\Q_p}$ is a quasi-split unitary group associated to the extension $F_p/\Q_p$.  However, we shall postpone discussion of this case to a future occasion.}, and then we specify the integral data at $p$.  

\noindent {\em First assumption on $p$}: The prime $p$ splits in $F$ as a product of distinct prime ideals 
$$(p) = {\mathfrak p} {\mathfrak p}^*,$$ 
where ${\mathfrak p}$ is the prime determined by our fixed choice of embedding $\overline{\Q} \hookrightarrow \overline{\Q}_p$, and 
${\mathfrak p}^*$ is its image under the non-trivial element of ${\rm Gal}(F/\Q)$.

Under this assumption $F_{\mathfrak p} = F_{{\mathfrak p}^*} = \Q_p$.  Further the algebra $D_{\Q_p}$ is a product
$$
D \otimes \Q_p = D_{\mathfrak p} \times D_{{\mathfrak p}^*}
$$
where each factor is a central simple $\Q_p$-algebra.  We have $D_{\mathfrak p} ~ \widetilde{\rightarrow} ~ 
D^{opp}_{{\mathfrak p}^*}$ via $*$.  Therefore for any $\Q_p$-algebra $R$ we can identify the group $G(R)$ with the group 
$$
\{ (x_1,x_2) \in (D_{\mathfrak p} \otimes R)^\times \times (D_{{\mathfrak p}^*} \otimes R)^\times ~ | ~ 
x_1 = c (x^*_2)^{-1}, \,\, \mbox{for some $c \in R^\times$} \}.
$$
Therefore there is an isomorphism of $\Q_p$-groups $G \cong D_{\mathfrak p}^\times \times \G_m$ given by $(x_1,x_2) \mapsto (x_1,c)$.

\noindent{\em Second assumption on $p$}: The algebra $D_{\Q_p}$ {\em splits}:  $D_{\mathfrak p} \cong M_n(\Q_p)$. 

In this case the involution $*$ becomes isomorphic to the involution on $M_n(\Q_p) \times M_n(\Q_p)$ given by 
$$
* : (X,Y) \mapsto (Y^t,X^t).
$$

Our assumptions imply that $G = {\rm GL}_n \times \G_m$, a split $p$-adic group (and thus $E := {\bf E}_{\mathfrak p} = \Q_p$).  Why is this helpful?  As we shall see, this allows us to use the local models for ${\rm GL}_n$ described in section \ref{local_models} to describe the reduction modulo $p$ of the Shimura variety $Sh_{K^p}$, see $\S \ref{fake_unitary_local_models}$.  Also, we can  use the description in \cite{HN1} of nearby cycles on such models to compute the semi-simple local zeta function at $p$ of 
$Sh_{K^p}$, see \cite{HN3} and Theorem \ref{HN3_main_thm}.  One expects that this is still possible in the general case (where $D_{\mathfrak p}^\times$ is not a split group), but there the crucial facts about nearby cycles on the corresponding local models are not yet established.

\noindent {\em Integral data.}
We need to specify a $\Z_{(p)}$-order $\cO_B \subset B$ and a self-dual multichain ${\mathcal L} = \{ \Lambda \}$ of $\cO_B \otimes \Z_p$-lattices.  To give a multichain we need to specify first a (partial) $\Z_p$-lattice chain in 
$V_{\Q_p} = D_{\mathfrak p} \times D_{{\mathfrak p}^*}$.  We do this one factor at a time.  First, we may fix an isomorphism
\begin{equation} \label{eq:decomposition}
D_{\Q_p} = D_{\mathfrak p} \times D_{{\mathfrak p}^*}  \cong M_n(\Q_p) \times M_n(\Q_p)
\end{equation}
such that the involution $x \mapsto x^\iota = \xi x^* \xi^{-1}$ goes over to $(X,Y) \mapsto (Y^t,X^t)$.  
So $\xi$ gets identified with an element of the form $(\chi^t, -\chi)$, for $\chi \in {\rm GL}_n(\Q_p)$
\footnote{For use in $\S \ref{ss_for_simple}$, note that if we multiply $\xi$ by any integral power of $p$, we change neither its properties nor the isomorphism class of the symplectic space $V, (\cdot, \cdot)$.  Hence we may assume $\chi^{-1} \in {\rm GL}_n(\Q_p) \cap M_n(\Z_p)$.}, and our pairing $(x,y) = tr_{D/\Q}(x \xi y^*) = tr_{D/\Q}(x y^\iota \xi)$ goes over to 
$$
\langle (X_1,X_2),(Y_1,Y_2) \rangle = tr_{D_{\Q_p}/\Q_p}( X_1 Y_2^t\chi^t, -X_2Y_1^t\chi).
$$

Next we define a (partial) $\Z_p$-lattice chain $\Lambda^*_{-n} \subset \cdots \subset \Lambda^*_0 = 
p^{-1}\Lambda^*_{-n}$ in $D_{{\mathfrak p}^*}$ by setting
$$
\Lambda^*_{-i} = \chi^{-1} \, {\rm diag}(p^i, 1^{n-i}) \, M_n(\Z_p),
$$
for $i = 0, 1, \dots, n$.  We can then extend this by periodicity to define $\Lambda^*_{i}$ for all $i \in \Z$.  
Similarly, we define the $\Z_p$-lattice chain $\Lambda_0 \subset \cdots \subset \Lambda_n = p^{-1}\Lambda_0$ in $D_{\mathfrak p}$ by setting
$$
\Lambda_i = {\rm diag}((p^{-1})^i,1^{n-i}) \, M_n(\Z_p),
$$
for $i = 0, 1, \cdots, n$ (and then extending by periodicity to define for all $i$).  We note that 
$$
(\Lambda_i \oplus \Lambda^*_i)^\perp = \Lambda_{-i} \oplus \Lambda^*_{-i},
$$
where $\perp$ is defined in the usual way using the pairing $(\cdot, \cdot)$.  
Setting $\cO_B \subset B$ to be the unique maximal $\Z_{(p)}$-order such that under our fixed identification $D_{\Q_p} \cong M_n(\Q_p) \times M_n(\Q_p)$, we have
$$
\cO_B \otimes \Z_p ~ \widetilde{\rightarrow} ~ M^{opp}_n(\Z_p) \times M^{opp}_n(\Z_p),
$$
one can now check that ${\mathcal L} := \{ \Lambda \oplus \Lambda^* \}$ is a self-dual multichain of $\cO_B \otimes \Z_p$-lattices.  It is clear that $(\cO_B \otimes \Z_p)^\iota = \cO_B \otimes \Z_p$.

\subsubsection{The moduli problem} \label{simple_moduli_problem}

We have now constructed all the data that enters into the definition of $Sh_{K^p}$.  By the determinant condition, the abelian varieties have (relative) dimension ${\rm dim}(V_1) = n^2$.  An $S$-point in our moduli space is a chain of abelian schemes over $S$ of relative dimension $n^2$, equipped with $\cO_B \otimes \Z_{(p)}$-actions, indexed by ${\mathcal L}$ (we set $A_i = A_{\Lambda_i \oplus \Lambda^*_i}$ for all $i \in \Z$)
$$
\xymatrix{
\cdots \ar[r]^{\alpha} & A_{0} \ar[r]^{\alpha} &  A_{1} \ar[r]^{\alpha} & \cdots \ar[r]^{\alpha}  & A_{n} 
\ar[r]^{\alpha} & \cdots 
}
$$
such that 
\begin{itemize}
\item  each $\alpha$ is an isogeny of height $2n$ (i.e., of degree $p^{2n}$); 
\item there is a ``periodicity isomorphism'' $\theta_p : A_{i+n} \rightarrow A_{i}$ such that for each $i$ the composition 
$$
\xymatrix{
A_i \ar[r]^{\alpha} & A_{i+1} \ar[r]^{\alpha} & \cdots \ar[r]^{\alpha} & A_{i + n} \ar[r]^{\theta_p} & A_i}
$$
is multiplication by $p: A_i \rightarrow A_i$;
\item the morphisms $\alpha$ commute with the $\cO_B \otimes \Z_{(p)}$-actions;
\item the determinant condition holds: for every $i$ and $b \in \cO_B$, 
$$
{\rm det}_{\cO_S}(b, {\rm Lie}(A_i)) = {\rm det}_{E'}(b, V_1).
$$
\end{itemize}
(See \cite{RZ}, Def. 6.5.)

In addition, we have a principal polarization and a $K^p$-level structure (see \cite{RZ}, Def. 6.9).  Giving a polarization is equivalent to giving a commutative diagram whose vertical arrows are isogenies
$$
\xymatrix{
\cdots \ar[r]^{\alpha} & A_{-1} \ar[r]^{\alpha} \ar[d] & A_{0} \ar[r]^{\alpha} \ar[d] &  A_{1} \ar[r]^{\alpha} 
\ar[d] & 
\cdots \ar[r]^{\alpha}  & A_{n} \ar[r]^{\alpha} \ar[d] & \cdots \\ 
\cdots \ar[r]^{\alpha^\vee} & \widehat{A}_{1} \ar[r]^{\alpha^\vee}   & \widehat{A}_{0} \ar[r]^{\alpha^\vee}  &  
\widehat{A}_{-1} \ar[r]^{\alpha^\vee}  & \cdots \ar[r]^{\alpha^\vee}  & 
\widehat{A}_{-n}  \ar[r]^{\alpha^\vee}   & \cdots \\
}
$$
such that for each $i$ the quasi-isogeny
$$
A_i \rightarrow \widehat{A}_{-i} \rightarrow \widehat{A}_i
$$
is a rational multiple of a polarization of $A_i$.  If up to a $\Q$-multiple the vertical arrows are all isomorphisms, we say the polarization is principal.

The fact that ${\rm End}_B(V)$ is a division algebra implies that the moduli space $Sh_{K^p}$ is proper over $\cO_E$ (Kottwitz verified the valuative criterion of properness in the case of maximal hyperspecial level structure in \cite{Ko92} p. 392, using the theory of N\'{e}ron models; the same proof applies here.)

\subsection{Siegel modular varieties with $\Gamma_0(p)$-level structure} \label{Siegel}
 
The set-up is much simpler here.  The group ${\bf G}$ is ${\rm GSp}(V)$ where $V$ is the standard symplectic space $\Q^{2n}$ with the alternating pairing $(\cdot, \cdot)$ given by the matrix $\tilde{I}$ in \ref{matrix}.  We have $B = \Q$ with involution $\iota = {\rm id}$, and $h_0 : \C \rightarrow {\rm End}(V_\R)$ is defined as the unique $\R$-algebra homomorphism such that 
$$
h_0(i) = \tilde{I}.
$$
For the multichain ${\mathcal L}$ we use the standard complete self-dual lattice chain $\Lambda_\bullet$ in 
$\Q_p^{2n}$ that appeared in section \ref{local_models}.  We take $\cO_B = \Z_{(p)}$.  

The group $G = {\rm GSp}_{2n,\Q_p}$ is split, so again we have $E = \Q_p$, so $\cO_E = \Z_p$.  It turns out that the minuscule coweight $\mu$ is 
$$
\mu = (0^n,(-1)^n),
$$
in other words, the same that appeared in the definition of local models in the symplectic case in \ref{local_models}. 

The moduli problem over $\Z_p$ can be expressed as follows.  For a $\Z_p$-scheme $S$, an $S$-point is an element of the set of 4-tuples (taken up to isomorphism) 
$$\cA_{K^p}(S) = \{ (A_\bullet, \lambda_0, \lambda_n, \bar{\eta}) \}$$
consisting of 
\begin{itemize}
\item a chain $A_\bullet$ of (relative) $n$-dimensional abelian varieties $A_0 \overset{\alpha}{\rightarrow} A_1 \overset{\alpha}{\rightarrow} \cdots \overset{\alpha}{\rightarrow} A_n$ such that each morphism $\alpha: A_i \rightarrow A_{i+1}$ is an isogeny of degree $p$ over $S$;
\item principal polarizations $\lambda_0: A_0 ~ \widetilde{\rightarrow} ~ \widehat{A}_0$ and $\lambda_n : A_n ~ \widetilde{\rightarrow} ~ \widehat{A}_n$ such that the composition of 
$$
\xymatrix{
A_0 \ar[r]^{\alpha} & \cdots \ar[r]^{\alpha} & A_n \ar[d]^{\lambda_n} \\
\widehat{A}_0 \ar[u]^{\lambda_0^{-1}} & \ar[l]_{\alpha^\vee} \cdots & \ar[l]_{\alpha^\vee} \widehat{A}_n}
$$
starting and ending at any $A_i$ or $\widehat{A}_i$ is multiplication by $p$;
\item a level $K^p$-structure $\bar{\eta}$ on $A_0$.
\end{itemize}

\xca Show that the above description of $\cA_{K^p}$ is equivalent to the definition given in Definition 6.9 of \cite{RZ} (cf. our Def. \ref{general_definition}) for the group-theoretic data $(B,\iota, V,...)$ we described above.

Note that the only information imparted by the determinant condition in this case is that ${\rm dim}(A_i) = n$ for every $i$.

There is another convenient description of the same moduli problem, used by de Jong \cite{deJ}.  We define another moduli problem $\cA'_{K^p}$ whose $S$-points is the set of 4-tuples
$$
\cA'_{K^p}(S) = \{ (A_0, \lambda_0, \bar{\eta}, H_\bullet)  \}
$$
consisting of 
\begin{itemize}
\item an $n$-dimensional abelian variety $A_0$ with principal polarization $\lambda_0$ and $K^p$-level structure $\bar{\eta}$;
\item a chain $H_\bullet$ of finite flat group subschemes of $A_0[p] := {\rm ker}(p: A_0 \rightarrow A_0)$ over $S$
$$
(0) = H_0 \subset H_1 \subset \cdots \subset H_n \subset A_0[p]
$$
such that $H_i$ has rank $p^i$ over $S$ and $H_n$ is totally isotropic with respect to the Riemann form $e_{\lambda_0}$, defined by the diagram
\end{itemize}
$$
\xymatrix{
A_0[p] \times A_0[p] \ar[r]^{e_{\lambda_0}} \ar[d]^{{\rm id} \times \lambda_0} & {\mathbb \mu}_p  \\ 
A_0[p] \times \widehat{A}_0[p] \ar[r]^{\cong}  & A_0[p] \times \widehat{A_0[p]} \ar[u]^{can}.
}
$$
Here $\widehat{A_0[p]} = {\rm Hom}(A_0[p],\G_m)$ denotes the Cartier dual of the finite group scheme $A_0[p]$ and {\em can} denotes the canonical pairing (which takes values in the $p$-th roots of unity group subscheme ${\mathbb \mu}_p \subset \G_m$).  (See \cite{Mu}, section 20, or \cite{Mi}, section 16.)

The isomorphism $\cA ~ \widetilde{\rightarrow} ~ \cA'$ is given by 
$$
(A_\bullet, \lambda_0, \lambda_n, \bar{\eta}) \mapsto (A_0, \lambda_0, \bar{\eta}, H_\bullet); \,\,\,\,\,\,\, H_i := {\rm ker}[\alpha^i : A_0 \rightarrow A_i].
$$
The inverse map is given by setting $A_i = A_0/H_i$ (the condition on $H_n$ allows us to define 
a principal polarization $\lambda_n : A_0/H_n ~ \widetilde{\rightarrow} ~ \widehat{(A_0/H_n)}$ using $\lambda_0$).

In \cite{deJ}, de Jong analyzed the singularities of $\cA$ in the case $n = 2$, and deduced that the model $\cA$ is flat in that case (by passing from $\cA$ to a local model $\Mloc$ according to the procedure of section \ref{relating_to_local_models} and then by writing down equations for $\Mloc$).
 
In the sequel, we will denote the model $\cA$ (and $\cA'$) by the symbol $Sh$, sometimes adding the subscript $K_p$ when the level-structure at $p$ is not already understood.

\section{Relating Shimura varieties and their local models} \label{relating_to_local_models}

\subsection{Local model diagrams} \label{local_model_diagrams}

Here we describe the desiderata for local models of Shimura varieties.  Quite generally, consider a diagram of finite-type $\cO_E$-schemes
$$
\xymatrix{
{\mathcal M} & \ar[l]_{\varphi} \widetilde{\mathcal M} \ar[r]^{\psi} & {\mathcal M}^{\rm loc} }.
$$

\begin{definition}  \label{def_local_model_diagrams} We call such a diagram a {\em local model diagram} provided the following conditions are satisfied:
\begin{enumerate}
\item the morphisms $\varphi$ and $\psi$ are smooth and $\varphi$ is surjective;
\item \'{e}tale locally ${\mathcal M} \cong {\mathcal M}^{\rm loc}$: there exists an \'{e}tale covering $V \rightarrow {\mathcal M}$ and a section $s : V \rightarrow \widetilde{\mathcal M}$ of $\varphi$ over $V$ such that $\psi \circ s : V \rightarrow {\mathcal M}^{\rm loc}$ is \'{e}tale.
\end{enumerate}
\end{definition}

In practice ${\mathcal M}$ is the scheme we are interested in, and ${\mathcal M}^{\rm loc}$ is somehow simpler to study; $\widetilde{\mathcal M}$ is just some intermediate scheme used to link the other two.  Every property that is local for the \'{e}tale topology is shared by ${\mathcal M}$ and ${\mathcal M}^{\rm loc}$.  For example, if ${\mathcal M}^{\rm loc}$ is {\em flat} over ${\rm Spec}(\cO_E)$, then so is ${\mathcal M}$.  The singularities in ${\mathcal M}$ and ${\mathcal M}^{\rm loc}$ are the same.  

\subsection{The general definition of local models}

We briefly recall the general definition of local models, following \cite{RZ}, Def. 3.27.  We suppose we have data $G,\mu, V, V_1, \dots$ coming from a PEL-type data as in section \ref{PEL_type_data}.  We assume $\mu$ and $V_1$ are defined over a finite extension $E' \supset E$.  We suppose we are given a self-dual multichain of $\cO_B \otimes \Z_p$-lattices ${\mathcal L} = \{ \Lambda \}$.  

\begin{definition}[\cite{RZ}, 3.27] \label{general_local_model_def}
A point of $\Mloc$ with values in an $\cO_E$-scheme $S$ is given by the following data.
\begin{enumerate}
\item A functor from the category ${\mathcal L}$ to the category of $\cO_B \otimes_{\Z_p} \cO_S$-modules on $S$
$$
\Lambda \mapsto t_\Lambda, \,\,\,\,\,\, \Lambda \in {\mathcal L};
$$
\item  A morphism of functors $\psi_\Lambda: \Lambda \otimes_{\Z_p} \cO_S \rightarrow t_\Lambda$.
\end{enumerate}
We require the following conditions are satisfied:
\begin{enumerate}
\item [(i)] $t_\Lambda$ is a locally free $\cO_S$-module of finite rank.  For the action of $\cO_B$ on $t_\Lambda$ we have the determinant condition
$$
{\rm det}_{\cO_S}(a; t_\Lambda) = {\rm det}_{E'}( a; V_1), \,\,\,\,\,\,\,\, a \in \cO_B;
$$
\item [(ii)] the morphisms $\psi_\Lambda$ are surjective;
\item [(iii)] for each $\Lambda$ the composition of the following map is zero:
$$
\xymatrix{
t^\vee_\Lambda \ar[r]^{\psi^\vee_\Lambda \,\,\,\,\,\,\,\,\,\,\,\,}  & (\Lambda \otimes \cO_S)^\vee \ar[r]_{\cong}^{(\cdot,\cdot)} & \Lambda^\perp \otimes \cO_S \ar[r]^{\,\,\,\,\,\,\,\,\,\,\, \psi_{\Lambda^\perp}} & t_{\Lambda^\perp}.
}
$$
\end{enumerate}
\end{definition} 

It is clear that one can associate to any PEL-type Shimura variety $Sh = Sh_{K_p}$ a scheme 
$\Mloc$ (just use the same PEL-type data and multichain ${\mathcal L}$ used to define $Sh_{K_p}$).  It is less clear why the resulting scheme $\Mloc$ really is a local model for $Sh_{K_p}$, in the sense described above.  
We shall see this below, thus justifying the terminology ``local model''.  Then we will show that in our two examples  -- the ``fake'' unitary and the Siegel cases -- this definition agrees with the concrete ones defined for ${\rm GL}_n$ and ${\rm GSp}_{2n}$ in section \ref{local_models}.

\subsection{Constructing local model diagrams for Shimura varieties} \label{constructing_local_model_diagrams}

\subsubsection{The abstract construction}

For one of our models $Sh = Sh_{K_p}$ from $\S \ref{PEL_parahoric_section}$, we want to construct a local model diagram
$$
\xymatrix{
Sh & \ar[l]_{\varphi} \widetilde{Sh} \ar[r]^{\psi} & \Mloc.
}$$
In the following we use freely the notation of the appendix, $\S \ref{appendix}$. 
For an abelian scheme $a: A \rightarrow S$, let $M(A)$ be the locally free $\cO_S$-module dual to the de Rham cohomology
$$
M^\vee(A) = H^1_{DR}(A/S) := R^1 a_*(\Omega^\bullet_{A/S}).
$$
This is a locally free $\cO_S$-module of rank $2 \, {\rm dim}(A/S)$. We have the Hodge filtration
$$
0 \rightarrow {\rm Lie}(\widehat{A})^\vee \rightarrow M(A) \rightarrow {\rm Lie}(A) \rightarrow 0.
$$
This is dual to the usual Hodge filtration on de Rham cohomology
$$
0 \subset \omega_{A/S} := a_* \Omega^1_{A/S}  \subset H^1_{DR}(A/S).
$$
We shall call $M(A)$ the {\em crystal associated to $A/S$} (this is perhaps non-standard terminology).   
If $A$ carries an action of $\cO_B$, then by functoriality so does $M(A)$.  Note also that $M(A)$ is covariant as a functor of $A$.  So if ${\mathcal L}$ denotes a self-dual multichain of $\cO_B \otimes \Z_p $-lattices, and $\{ A_\Lambda \}$ denotes an ${\mathcal L}$-set of abelian schemes over $S$ with $\cO_B$-action and polarization (as in Definition \ref{general_definition}), then applying the functor $M(\cdot)$ gives us a polarized multichain $\{M(A_\Lambda) \}$ of $\cO_B \otimes \cO_S $-modules of type $({\mathcal L})$, in the sense of \cite{RZ}, Def. 3.6, 3.10, 3.14.  

One key consequence of the conditions imposed in loc. cit., Def. 3.6, is that locally on $S$ there is an isomorphism of polarized multichains of $\cO_B \otimes \cO_S$-modules
$$
\gamma_\Lambda : M(A_\Lambda) ~ \widetilde{\rightarrow} ~  \Lambda \otimes_{\Z_p} \cO_S.
$$
In fact we have the following result which guarantees this.

\begin{theorem} [\cite{RZ}, Theorems 3.11, 3.16] \label{RZ_thm_3.16}
Let ${\mathcal L} = \{ \Lambda \}$ be a (self-dual) multichain of $\cO_B \otimes \Z_p$-lattices in $V$.
Let $S$ be any $\Z_p$-scheme where $p$ is locally nilpotent.  Then any (polarized) multichain $\{ M_\Lambda \}$ of $\cO_B \otimes_{\Z_p} \cO_S$-modules of type $({\mathcal L})$ is locally (for the \'{e}tale topology 
on $S$) isomorphic to the (polarized) multichain $\{ \Lambda \otimes_{\Z_p} \cO_S \}$.

Moreover, the functor ${\rm Isom}$
$$
T \mapsto {\rm Isom}(\{ M_\Lambda \otimes \cO_T \} , \{ \Lambda \otimes \cO_T \}),
$$
is represented by a smooth affine scheme over $S$.
\end{theorem}

The analogous statements hold for any $\Z_p$-scheme $S$, see \cite{P}.  In particular for such $S$ we have a smooth affine group scheme ${\mathcal G}$ over $S$ given by 
$$
{\mathcal G}(T) = {\rm Aut}(\{ \Lambda \otimes \cO_T \}),
$$
and the functor ${\rm Isom}$ is obviously a left-torsor under ${\mathcal G}$.
This generalizes the smoothness of the groups ${\rm Aut}$ in section \ref{BT_group_schemes}.  
Moreover, by the same arguments as in section \ref{BT_group_schemes}, for $S = {\rm Spec}(\Z_p)$ the 
group ${\mathcal G}_{\Z_p}$ is a Bruhat-Tits parahoric group scheme corresponding to the parahoric 
subgroup of 
$G(\Q_p) = {\mathcal G}(\Q_p)$ which stabilizes the multichain ${\mathcal L}$ \footnote{More precisely, the {\em connected component} of ${\mathcal G}$ is the Bruhat-Tits group scheme.  As G. Pappas points out, in some cases (e.g. the unitary group for ramified quadratic extensions), the stabilizer group ${\mathcal G}$ is not connected.} .

In the special case of lattice chains for ${\rm GSp}_{2n}$, the theorem was proved by de Jong \cite{deJ} (he calls what are ``polarized (multi)chains'' here by the name ``systems of $\cO_S$-modules of type II''). 

Now we define the local model diagram for $Sh$.  We assume $\cO_E = \Z_p$ for simplicity.  Let us define $\widetilde{Sh}$ to be the $\Z_p$-scheme representing the functor whose points in a $\Z_p$-scheme $S$ is the set of pairs 
$$
(\{A_\Lambda \}, \bar{\lambda}, \bar{\eta}) \in Sh(S);\,\,\,\,\,\, \gamma_\Lambda: M(A_\Lambda) ~ \widetilde{\rightarrow} ~ \Lambda \otimes_{\Z_p} \cO_S,
$$
where $\gamma_\Lambda$ is an isomorphism of polarized multichains of $\cO_B \otimes \cO_S$-modules.  The morphism 
$$
\varphi: \widetilde{Sh} \rightarrow Sh
$$ 
is the obvious morphism which forgets $\gamma_\Lambda$.  By Theorem \ref{RZ_thm_3.16}, $\varphi$ is smooth (being a torsor for a smooth group scheme) and surjective.  Now we want to define
$$
\psi: \widetilde{Sh}(S) \rightarrow \Mloc(S).
$$
We define it to send an $S$-point $(\{A_\Lambda \}, \bar{\lambda}, \bar{\eta}, \gamma_\Lambda)$ to the morphism of functors
$$
\Lambda \otimes_{\Z_p} \cO_S \rightarrow {\rm Lie}(A_{\Lambda})
$$
induced by the composition $\gamma^{-1}_\Lambda : \Lambda \otimes_{\Z_p} \cO_S \cong M(A_\Lambda)$ with the canonical surjective morphism
$$
M(A_\Lambda) \rightarrow {\rm Lie}(A_\Lambda).
$$
It is not completely obvious that the morphisms $\Lambda \otimes \cO_S \rightarrow {\rm Lie}(A_{\Lambda})$ satisfy the condition (iii) of Definition \ref{general_local_model_def}.  We will explain it in the Siegel case below, as a consequence of Proposition 5.1.10 of \cite{BBM} (our Prop. \ref{BBM_5.1.10}).  We omit discussion of this point in other cases.

The theory of Grothendieck-Messing (\cite{Me}) shows that the morphism $\psi$ is formally smooth.  Since both schemes are of finite type over $\Z_p$, it is smooth.  In summary:

\begin{theorem}[\cite{RZ}, $\S3$] \label{local_model_diagram_exists}
The diagram 
$$
\xymatrix{
Sh & \ar[l]_{\varphi} \widetilde{Sh} \ar[r]^{\psi} & \Mloc}
$$
is a local model diagram.  The morphism $\varphi$ is a torsor for the smooth affine group scheme ${\mathcal G}$. 
\end{theorem}

\begin{proof} We have indicated why condition (1) of Definition \ref{def_local_model_diagrams} is satisfied.  Condition (2) is proved in \cite{RZ}, 3.30-3.35; see also \cite{deJ}, Cor. 4.6.
\end{proof}

We will describe the local model diagrams more explicitly for each of our two main examples next.  Our goal is to show that their local models are none other than the ones defined in section \ref{local_models}.

\subsubsection{Symplectic case} \label{local_model_Siegel_case}

Following \cite{deJ} and \cite{GN} we change conventions slightly and replace the de Rham homology functor with 
the cohomology functor
$$
A/S \mapsto H^1_{DR}(A/S).
$$

What kind of data do we get by applying the de Rham cohomology functor to a point in our moduli problem $Sh$ from section \ref{Siegel}?  For notational convenience, let us now number the chains of abelian varieties in the opposite order:
$$
\{A_{\Lambda_\bullet} \} = A_n \rightarrow A_{n-1} \rightarrow \cdots \rightarrow A_0.
$$

\begin{lemma} \label{applying_dR_cohomology} The result of applying $H^1_{DR}$ to a point $(\{A_{\Lambda_\bullet} \}, \lambda_0, \lambda_n)$ in $Sh(S)$ is a datum of form
$(M_0 \overset{\alpha}{\rightarrow} M_{1} \overset{\alpha}{\rightarrow} \cdots \overset{\alpha}{\rightarrow} M_n, q_0,q_n)$ satisfying
\begin{itemize}
\item $M_i$ is a locally free $\cO_S$-module of rank $2n$;
\item ${\rm Coker}(M_{i-1} \rightarrow M_{i})$ is a locally free $\cO_S/p\cO_S$-module of rank 1;
\item for $i = 0,n$,  $\,\, q_i : M_i \otimes  M_i \rightarrow \cO_S$ is a non-degenerate symplectic pairing;
\item for any $i$, the composition of 
$$
\xymatrix{
M_0  \ar[r]^{\alpha} & \cdots \ar[r]^{\alpha} &  M_n  \ar[d]^{q_n} \\
M_0^\vee \ar[u]^{q_0} & \ar[l]_{\alpha^\vee} \cdots & \ar[l]_{\alpha^\vee} M_n^\vee 
}
$$
starting and ending at $M_i$ or $M^\vee_i := {\rm Hom}(M_i, \cO_S)$ is multiplication by $p$.
\end{itemize}
\end{lemma}

\begin{proof}
The pairings $q_0, q_n$ come from the polarizations $\lambda_0, \lambda_n$.  The various properties are easy to 
check, using the canonical natural isomorphism $H^1_{DR}(\widehat{A}/S) = (H^1_{DR}(A/S))^\vee$; cf. Prop. \ref{BBM_5.1.10}.
\end{proof}

\smallskip

Our next goal is to rephrase Definition \ref{general_local_model_def} in terms of data similar to that in Lemma \ref{applying_dR_cohomology}, which will take us closer to the definition of $\Mloc$ in $\S \ref{local_models}$.

Let ${\mathcal L} = \Lambda_\bullet$ be the standard self-dual lattice chain in 
$V = \Q_p^{2n}$, with respect to the usual pairing $(x,y) = x^t \tilde{I}y$.  
Clearly we may rephrase Definition \ref{general_local_model_def} using the sub-objects 
$\omega'_\Lambda := {\rm ker}(\psi_\Lambda)$ of $\Lambda \otimes \cO_S$ rather than the 
quotients $t_\Lambda$.  Then condition (iii) becomes
\begin{enumerate}
\item [(iii')] $(\omega'_{\Lambda})^{perp} \subset \omega'_{\Lambda^{\perp}}, \,\,\,\, \mbox{which is equivalent to 
$(\omega'_{\Lambda})^{perp} = \omega'_{\Lambda^{\perp}}$,}$ 
\end{enumerate}
in other words, under the canonical pairing $(\cdot,\cdot):\Lambda \otimes \cO_S \times \Lambda^\perp \otimes \cO_S \rightarrow \cO_S$, the submodules $\omega'_{\Lambda}$ and $\omega'_{\Lambda^\perp}$ are perpendicular.  If $\Lambda = \Lambda^\perp$, this means 
$$
(\cdot, \cdot)|_{\omega'_{\Lambda} \times \omega'_{\Lambda}} \equiv 0,
$$
and if $\Lambda^\perp = p \Lambda$, this means
$$
p(\cdot, \cdot)|_{\omega'_{\Lambda} \times \omega'_{\Lambda}} \equiv 0,$$
 since the pairing on $\Lambda \otimes \cO_S \times \Lambda \otimes \cO_S$ is defined by composing the standard pairing on $\Lambda \otimes \cO_S \times \Lambda^\perp \otimes \cO_S$ with the periodicity isomorphism
$$
p: \Lambda ~ \widetilde{\rightarrow} ~ \Lambda^\perp
$$
in the second variable.  For the ``standard system'' $(\Lambda_0 \rightarrow \Lambda_1 \rightarrow \cdots \rightarrow \Lambda_n, q_0, q_n)$ as in Lemma \ref{applying_dR_cohomology}, the (perfect) pairings are given by
\begin{align*}
q_0 &= (\cdot, \cdot) : \Lambda_0 \times \Lambda_0 \rightarrow \Z_p \\
q_n &= p(\cdot, \cdot) : \Lambda_n \times \Lambda_n \rightarrow \Z_p.
\end{align*}

Note that if $\omega'_i := \omega'_{\Lambda_i}$, the identity 
$(\omega'_{i})^{perp} = \omega'_{-i}$ ((iii) of Def. \ref{general_local_model_def}) means that $\omega'_\bullet$ is uniquely determined by the elements $\omega'_0, \dots, \omega'_n$.  Conversely, suppose we are given $\omega'_0, \dots, \omega'_n$ such that $(\omega'_0)^{perp} = \omega'_0$ and $(\omega'_n)^{perp} = p \, \omega'_n =: \omega'_{-n}$.  Then we can {\em define} $\omega'_{-i} = (\omega'_i)^{perp}$ for $i = 0, \dots, n$, and then extend by periodicity to get an infinite chain $\omega'_\bullet$ as in Definition \ref{general_local_model_def} (condition (iii) being satisfied by fiat).  

We thus have the following reformulation of Definition \ref{general_local_model_def}, which shows that that definition agrees with the one in section \ref{local_models} for ${\rm GSp}_{2n}$.

\begin{lemma}
In the Siegel case, an $S$-point of  $\Mloc$ (in the sense of Definition \ref{general_local_model_def}) is a commutative diagram
$$
\xymatrix{
\Lambda_0 \otimes \cO_S \ar[r] & \Lambda_1 \otimes \cO_S \ar[r] & \cdots \ar[r] & \Lambda_n \otimes \cO_S \\
\omega'_0 \ar[u] \ar[r] & \omega'_1 \ar[u] \ar[r] & \cdots \ar[r] & \omega'_n, \ar[u]}
$$
such that 
\begin{itemize}
\item for each $i$, $\omega'_i$ is a locally free $\cO_S$-submodule of $\Lambda_i \otimes \cO_S$ of rank $n$;
\item $\omega'_0$ is totally isotropic for $(\cdot, \cdot)$ and $\omega'_n$ is totally isotropic for 
$p(\cdot, \cdot)$.
\end{itemize}
\end{lemma}

Finally, we promised to explain why the morphism $\psi: \widetilde{Sh} \rightarrow \Mloc$ really takes values in $\Mloc$.  We must also redefine it in terms of cohomology.  Recall we now have the Hodge filtration
$$
\omega_{A_{\Lambda}/S} \subset H^1_{DR}(A_{\Lambda}/S).
$$ 
We define $\psi$ to send $(\{A_0 \leftarrow \cdots \leftarrow A_n \},\lambda_0,\lambda_n,\bar{\eta},\gamma_\Lambda)$ to the locally free, rank $n$,  $\cO_S$-submodules
$$
\gamma_\Lambda(\omega_{A_\Lambda}) \subset \Lambda \otimes \cO_S,
$$
where now  $\gamma_\Lambda$ is an isomorphism of polarized multichains of $\cO_S$-modules
$$
\gamma_\Lambda: H^1_{DR}(A_\Lambda/S) ~ \widetilde{\rightarrow} ~ \Lambda \otimes \cO_S.
$$
The following result ensures that this map really takes values in $\Mloc$.

\begin{lemma} The morphism $\psi$ takes values in $\Mloc$, i.e., condition (iii') holds.
\end{lemma}
\begin{proof} 
Setting $\omega_{A_{\Lambda_i}} = \omega_i$, we need to see that the Hodge filtration $\omega_0$ resp. $\omega_n$ is totally isotropic 
with respect to the pairing $q_0$ resp. $q_n$ induced by the polarization $\lambda_0$ resp. $\lambda_n$.  But this is Proposition 5.1.10 of \cite{BBM} (our Prop. \ref{BBM_5.1.10}).  See also \cite{deJ}, Cor. 2.2.
\end{proof}

\noindent {\em Comparison of homology and cohomology local models.}  One further remark is in order.  Let us consider a point 
$$A = (A_0 \rightarrow \cdots \rightarrow A_n, \lambda_0, \lambda_n, \bar{\eta})
$$ 
in our moduli problem $Sh$.  Note that this data gives us another point in $Sh$, namely
$$
\widehat{A} = (\widehat{A}_n \rightarrow \cdots \rightarrow \widehat{A}_0, \lambda^{-1}_n, \lambda^{-1}_0, \bar{\eta}).
$$
(We need to use the assumption that the polarizations $\lambda_i$ are required to be {\em symmetric} isogenies $A_i \rightarrow \widehat{A}_i$, in the sense that $\widehat{\lambda}_i = \lambda_i$.)

The moduli problem $Sh$ is thus equipped with an automorphism of order 2, given by $A \mapsto \widehat{A}$.  

This comes in handy in comparing the ``homology'' and ``cohomology'' constructions of the local model diagram.  Namely, since $M(A_i) = H^1_{DR}(\widehat{A}_i)$ (Prop. \ref{BBM_5.1.10}), an isomorphism 
$
\gamma_\bullet: M(A_\bullet) \,\, \widetilde{\rightarrow} \,\, \Lambda_\bullet \otimes \cO_S
$
is simultaneously an isomorphism $\gamma_\bullet: H^1_{DR}(\widehat{A}_\bullet) \,\, \widetilde{\rightarrow} \,\, \Lambda_\bullet \otimes \cO_S$.  In the ``homology'' version, $\psi$ sends $(A,\gamma_\bullet)$ to the quotient chain
$$
\Lambda_\bullet \otimes \cO_S \rightarrow {\rm Lie}(A_\bullet),
$$
defined using $\gamma^{-1}_\bullet$.  On the other hand, in the ``cohomology'' version, $\psi$ sends $(\widehat{A}, \gamma_\bullet)$ to the sub-object chain
$$
\omega_{\widehat{A}_\bullet}  \subset \Lambda_\bullet \otimes \cO_S
$$
(identifying $\omega_\bullet$ with $\gamma_\bullet(\omega_\bullet)$).
But the exact sequence 
$$
0 \rightarrow \omega_{\widehat{A}} \rightarrow M(A) \rightarrow {\rm Lie}(A) \rightarrow 0
$$
(Prop. \ref{BBM_5.1.10}) means that the two chains correspond: they give exactly the same element of $\Mloc$.  In summary, we have the following result relating the ``homology'' and ``cohomology'' constructions of the local model diagram.

\begin{prop} \label{homology_vs_cohomology_local_models}
There is a commutative diagram
$$
\xymatrix{
\widetilde{Sh}^{hom} \ar[r]^{\psi^{hom}} \ar[d] & \Mloc \ar[d]^{=} \\
\widetilde{Sh}^{coh} \ar[r]^{\psi^{coh}} & \Mloc,
}
$$
where the left vertical arrow is the automorphism $(A, \gamma_\bullet) \mapsto (\widehat{A}, \gamma_\bullet)$.
\end{prop}

\subsubsection{``Fake'' unitary case} \label{fake_unitary_local_models}

Here the ``standard'' polarized multichain of $\cO_B \otimes \Z_p$-lattices is given by 
$\{ \Lambda_i \oplus \Lambda^*_i \}$, in the notation of section \ref{integral_data}.  Recall that
$$
\cO_B \otimes \Z_p \cong M^{opp}_n(\Z_p) \times M^{opp}_n(\Z_p),
$$
according to the decomposition of $B_{\Q_p}^{opp} = D_{\Q_p}$ 
$$
D_{\Q_p} = D_{\mathfrak p} \times D_{{\mathfrak p}^*} \cong M_n(\Q_p) \times M_n(\Q_p).
$$
Let $W$ (resp. $W^*$) be $\Z_p^n$ viewed as a left $\cO_B \otimes \Z_p$-module, via right multiplications by elements of the first (resp. second) factor of $M_n(\Z_p) \times M_n(\Z_p)$.  
The ring $B_{\Q_p}$ has two simple left modules: $W_{\Q_p}$ and $W^*_{\Q_p}$.  
We may write
$$
V_{E'} = V_1 \oplus V_2
$$
as before.  The determinant condition now implies (at least over $E'$) that 
$$
V_1 = W_{E'}^{d} \oplus (W^*_{E'})^{n-d};
$$
(comp. section \ref{mu}).  Using the ``sub-object'' variant of Definition \ref{general_local_model_def}, it follows that an $S$-point of $\Mloc$ is a commutative diagram (here $\Lambda_i$ being understood as $\Lambda_i \otimes 
\cO_S$)
$$
\xymatrix{
\Lambda_0 \oplus \Lambda^*_0 \ar[r] & \Lambda_1 \oplus \Lambda^*_1 \ar[r] & \cdots \ar[r] & \Lambda_n \oplus \Lambda^*_n \\
\cF_0 \oplus \cF^*_0 \ar[r] \ar[u] & \cF_1 \oplus \cF^*_1 \ar[r] \ar[u] & \cdots \ar[r]  & \cF_n \oplus \cF^*_n \ar[u],}
$$
where $\,\, \cF_i \oplus \cF^*_i$ is an $\cO_B \otimes \cO_S$-submodule of $\Lambda_i \oplus \Lambda^*_i$ which, locally on $S$, is a direct factor isomorphic to $W_{\cO_S}^{n-d} \oplus (W_{\cO_S}^*)^{d}$.

\smallskip

The analogue of condition  (iii'), which is imposed in Definition \ref{general_local_model_def}, is 
$$
(\cF_i \oplus \cF^*_i)^{perp} = \cF_{-i} \oplus \cF^*_{-i}.
$$
On the other hand, from the definition of $\langle \cdot,\cdot \rangle$ in section \ref{integral_data} it is immediate that
$$
(\cF_i \oplus \cF^*_i)^{perp} = \cF^{*,perp}_{i} \oplus \cF^{perp}_{i}.
$$
We see thus that the first factor $\cF_\bullet$ uniquely determines the second factor $\cF^*_\bullet$ (and vice-versa).  Thus $\Mloc$ is given by chains of right $M_n(\cO_S) = M_n(\Z_p) \otimes \cO_S$-modules
$$
\cF_0 \rightarrow \cF_1 \rightarrow \cdots \rightarrow \cF_n$$
which are locally direct factors in 
$$
M_n(\cO_S) \rightarrow {\rm diag}(p^{-1},1^{n-1}) M_n(\cO_S) \rightarrow \cdots \rightarrow p^{-1}M_n(\cO_S),$$
each term locally isomorphic to $(\cO_S^n)^{n-d}$.  By Morita equivalence, $\Mloc$ is just given by the definition in section \ref{local_models} (for the integer $d$).

\section{Flatness} \label{flatness}

Because of the local model diagram, the flatness of the moduli problem $Sh$ can be investigated by considering its local model.  The following fundamental result is due to U. G\"{o}rtz.  It applies to all parahoric subgroups.

\begin{theorem} [\cite{Go1}, \cite{Go2}]
Suppose ${\bf M}^{\rm loc}$ is a local model attached to a group ${\rm Res}_{F/\Q_p}({\rm GL}_n)$ or ${\rm Res}_{F/\Q_p}({\rm GSp}_{2n})$, where $F/\Q_p$ is an unramified extension.  Then ${\bf M}^{\rm loc}$ is flat over $\cO_E$.  Moreover, its special fiber is reduced, and has rational singularities.
\end{theorem}

We give the idea for the proof.  One reduces to the case where $F = \Q_p$.  In order to prove flatness over $\cO_E = \Z_p$ it is enough to prove the following facts (comp. \cite{Ha}, III.9.8):

\noindent 1)  The special fiber is reduced, as a scheme over $\F_p$;

\noindent 2)  The model is {\em topologically flat}: every closed point in the special fiber is contained in the the scheme-theoretic closure of the generic fiber.

The innovation behind the proof of 1) is to embed the special fiber into the affine flag variety and then to make systematic use of the theory of Frobenius-splitting to prove affine Schubert varieties 
are compatibly Frobenius-split.  See \cite{Go1}.

To prove 2), suppose $\mu$ is such that $\Mloc = M_{-w_0\mu}$.  It is enough by a result of Kottwitz-Rapoport (Theorem \ref{Perm=Adm}) to prove that the generic element in a stratum of the special fiber indexed by a translation element in ${\rm Adm}(-w_0\mu)$ can be lifted to characteristic zero.  This statement is checked by hand in \cite{Go1}.  

We will provide an alternative, calculation-free, proof by making use of nearby cycles \footnote{We emphasize that this proof is much less elementary than the original proof of G\"{o}rtz \cite{Go1}, relying as it does on the full strength of \cite{HN1}.}.  We freely make use of material on nearby cycles from $\S \ref{nearby_cycles}, \ref{ss_for_simple}$.

We fix an element $\lambda \in W(-w_0\mu)$ and consider the stratum of $\Mloc$ indexed by $t_\lambda$.  We want to show this stratum is in the closure of the generic fiber.

The nearby cycles sheaf $R\Psi^{\Mloc}(\Ql)$ is supported only on this closure (Theorem \ref{RPsi_properties}), and so it is enough to show that 
$$
{\rm Tr}^{ss}(\Phi^r_p, R\Psi_{t_\lambda}^{\Mloc}(\Ql)) \neq 0.
$$
But it is clear that
$$
z_{-w_0\mu,r}(t_\lambda) \neq 0
$$
from the definition of Bernstein functions (see \cite{Lu} or \cite{HKP}), and we are done by Theorem \ref{RPsi=z_mu} (which also holds for the group ${\rm GSp}_{2n}$, see \cite{HN1}).

As G. Pappas has observed \cite{P}, the local models attached (by \cite{RZ}) to the groups above can fail to be flat if $F/\Q_p$ is {\em ramified}.  In their joint works \cite{PR1}, \cite{PR2}, Pappas and Rapoport provide alternative definitions of local models in that case (in fact they treat nearly all the groups considered in \cite{RZ}), and these new models {\em are} flat.  However, these new models cannot always be described as the scheme representing a ``concrete'' moduli problem.

\section{The Kottwitz-Rapoport stratification} \label{KR_strata}

Let us assume $Sh$ is the model over $\cO_E$ for one of the Shimura varieties $Sh({\bf G},h)_{\bf K}$ discussed in 
$\S 5$, i.e. a ``fake'' unitary or a Siegel modular variety.  We assume (for simplicity of statements) that $K_p$ is an Iwahori subgroup of $G(\Q_p)$.  Let us summarize what we know so far.

The group $G_{\Q_p}$ is either isomorphic to ${\rm GL}_{n, \Q_p} \times \G_{m, \Q_p}$ or ${\rm GSp}_{2n, \Q_p}$.  These groups being split over $\Q_p$, we have $E = \Q_p$ and $\cO_E = \Z_p$.

The Shimura datum $h$ gives rise to a dominant minuscule cocharacter $\mu$ of ${\rm GL}_{n,\Q_p}$
or ${\rm GSp}_{2n,\Q_p}$, respectively. The functorial description of the local
model ${\bf M}^{\rm loc}$ shows that it can be embedded into the deformation $M$ from the affine Grassmannian $\Grass_{\Q_p}$ 
to the affine flag variety $\Fl_{\F_p}$ associated to $G$, and has generic fibre $\mathcal
Q_{-w_0\mu}$. Since the local model is flat with reduced special fiber \cite{Go1}, \cite{Go2} (see $\S \ref{flatness}$) and is closed in $M$, it coincides with the scheme-theoretic
closure $M_{-w_0\mu}$ of $\mathcal Q_{-w_0\mu}$ in this deformation.  We thus identify $\Mloc = M_{-w_0\mu}$.  (For all this, keep in mind we use the ``homology'' definition of the local model diagram.)

The relation between the model of the Shimura variety over $\Z_p$ and its
local model is given by a diagram 
$$\xymatrix{Sh & \widetilde{Sh} \ar[l]_\varphi \ar[r]^\psi & \Mloc
}$$
of $\Z_p$-schemes, where $\varphi$ is a torsor under the smooth affine group
scheme ${\mathcal G}$ of $\S \ref{RZ_thm_3.16}$ (also termed ${\rm Aut}$ in $\S \ref{local_models}$), and $\psi$ is smooth. The fibres of $\varphi$ are geometrically
connected (more precisely, this holds for the restriction of $\varphi$ to any 
geometric connected component of $\widetilde{Sh}$).  
One can show that \'etale-locally around each
point of the special fiber of $Sh$, the schemes $Sh$ and $\Mloc$ are isomorphic.

The stratification of the special fibre of $\Mloc$ (by Iwahori-orbits) induces stratifications of
the special fibers of $\widetilde{Sh}$ and $Sh$ (see below).  The resulting stratification of $Sh_{\F_p}$ is called the {\em Kottwitz-Rapoport (or KR-) stratification}.

\subsection{Construction of the KR-stratification} \label{construction_of_KR_stratification}

Essentially following \cite{GN}, we will recall the construction and basic properties of the KR-stratification.  The difference between their treatment and ours is that they construct local models in terms of 
de Rham {\em cohomology}, whereas here they are constructed in terms of de Rham {\em homology}.  This is done for compatibility with the computations in $\S \ref{ss_for_simple}$.

For later use in $\S \ref{ss_for_simple}$, we give a detailed treatment  here for the case of ``fake'' unitary Shimura varieties. 

Let $k$ denote the algebraic closure of the residue field of $\Z_p$, and let 
$\widetilde{\Lambda}_\bullet = \Lambda_\bullet \oplus \Lambda^*_\bullet$ denote the self-dual multichain of $\cO_B \otimes \Z_p$-lattices from $\S \ref{integral_data}$.   Recall that a point in ${\bf M}^{\rm loc}(k)$ is a ``quotient chain'' of $k$-vector spaces
$$
\widetilde{\Lambda}_{\bullet} \otimes k \rightarrow t_{\widetilde{\Lambda}_\bullet},
$$
self-dual in the sense of Definition \ref{general_local_model_def} (iii), and such that each $t_{\widetilde{\Lambda}_i}$ satisfies the determinant condition, that is,
$$
t_{\widetilde{\Lambda}_i} = W_k^d \oplus (W^*)_k^{n-d}
$$
as $\cO_B \otimes k$-modules.  We can identify this object with a lattice chain in the affine flag variety for ${\rm GL}_n(k((t)))$ as follows.  Let $\cV_{\bullet,k}$ denote the ``standard'' complete lattice chain from $\S 4.1$.  Using duality and Morita equivalence (see $\S \ref{fake_unitary_local_models}$), the quotient $t_{\widetilde{\Lambda}_\bullet}$ can be identified with a quotient $t_{\Lambda_\bullet}$ of the standard lattice chain $\cV_{\bullet,k} \subset k\xT^n$.  Then we may write
$$
t_{\Lambda_i} = \cV_{i,k}/\cL_i
$$
for a unique lattice chain $\cL_{\bullet} = (\cL_0 \subset \cdots \subset \cL_n = t^{-1}\cL_0)$ consisting of $k[[t]]$-submodules of $k((t))^n$ which satisfy for each $i= 0, \dots, n$,
\begin{itemize}
\item $t\cV_{i,k} \subset \cL_i \subset \cV_{i,k}$;
\item the $k$-vector space $\cV_{i,k}/\cL_i$ has dimension $d$ (determinant condition).
\end{itemize}

The set of such lattice chains $\cL_\bullet$ is the special fiber of the model $M_{-w_0\mu}$ attached to the dominant coweight $-w_0\mu = (1^d,0^{n-d})$ of ${\rm GL}_n$.  Indeed, the two conditions above mean that for each $i$,
$$
{\rm inv}_{K}(\cL_i,\cV_{i,k}) = \mu
$$
and thus
$$
{\rm inv}_{K}(\cV_{i,k}, \cL_i) = -w_0\mu.
$$
Here ${\rm inv}_K$ is the standard notion of relative position of $k\xt$-lattices in $k\xT^n$, relative to the base point $\cV_{0,k} = k\xt^n$: we say ${\rm inv}_K(g\cV_{0,k},g'\cV_{0,k}) = \lambda \in X_+(T)$ if $g^{-1}g' \in K\lambda K$, where $K = {\rm GL}_n(k\xt)$.  Recall we have identified $\mu$ with $(0^{n-d}, (-1)^{d})$ and have embedded coweights into the loop group by the rule $\lambda \mapsto  \lambda(t)$.

If $\cL_\bullet = x(\cV_{\bullet,k})$ for $x \in \widetilde{W}({\rm GL}_n)$, this means that $x \in {\rm Perm}(-w_0\mu)$, which is also the set ${\rm Adm}(-w_0\mu)$, see $\S \ref{local_models}$.

Recall that the Iwahori subgroup $I = I_{k \xT}$ fixing $\cV_{\bullet,k}$ preserves the subset $M_{-w_0\mu,k} \subset {\mathcal Fl}_k$ and so via the identification ${\bf M}^{\rm loc} = M_{-w_0\mu}$, it also acts on the local model.  The Iwahori-orbits give a cellular decomposition
$$
{\bf M}^{\rm loc}_{k} = \coprod_{w \in {\rm Adm}(\mu)} {\bf M}^{\rm loc}_w.
$$

Here we define ${\bf M}^{\rm loc}_w$ to be the set of $\cL_\bullet$ above such that 
$$
{\rm inv}_I(\cV_{\bullet,k}, \cL_\bullet) = w^{-1},
$$
or equivalently
$$
{\rm inv}_I(\cL_\bullet, \cV_{\bullet,k}) = w,
$$
for $w^{-1} \in {\rm Adm}(-w_0\mu)$ (which happens if and only if $w \in {\rm Adm}(\mu)$).  Here we define ${\rm inv}_I(g \cV_{\bullet,k}, g' \cV_{\bullet,k}) = w$ if $g^{-1}g' \in IwI$.  

Each stratum is smooth (in fact ${\bf M}^{\rm loc}_w = \A^{\ell(w)}$), and the closure relations are determined by the Bruhat order on $\widetilde{W}$; that is, ${\bf M}^{\rm loc}_w \subset \overline{{\bf M}^{\rm loc}_{w'}}$ if and only if $w \leq w'$.

There is a surjective homomorphism $I_{k \xT} \rightarrow {\rm Aut}_{k}$, where ${\rm Aut}$ is the group scheme $\mathcal G$ of Theorem \ref{RZ_thm_3.16} which acts on the whole local model diagram.  The action of $I_{k \xT}$ on ${\bf M}^{\rm loc}_{k}$ factors through ${\rm Aut}_{k}$, so that the strata above can also be thought of as ${\rm Aut}_{k}$-orbits.  The morphism $\psi$ is clearly equivariant for ${\rm Aut}_{k}$, hence we have a stratification of $\widetilde{Sh}_{k}$
$$
\widetilde{Sh}_{k} = \coprod_{w \in {\rm Adm}(\mu)} \psi^{-1}( {\bf M}^{\rm loc}_w),
$$
whose strata are non-empty (Lemma \ref{Genestier_surjectivity}), smooth, and stable under the action of ${\rm Aut}_{k}$.  Since $\varphi_{k}$ is a torsor for the smooth group scheme ${\rm Aut}_{k}$, the stratification descends to $Sh_{k}$:
$$
Sh_{k} = \coprod_{w \in {\rm Adm}(\mu)} Sh_w
$$
such that $\varphi^{-1}(Sh_w) = \psi^{-1}({\bf M}^{\rm loc}_w)$.  These strata are still smooth, non-empty, and satisfy closure relations determined by the Bruhat order.  

All statements above remain true over the base field $\F_p$ instead of its algebraic closure $k$.

\medskip

\subsection{Relating nearby cycles on local models and Shimura varieties}

We will need in $\S \ref{ss_for_simple}$ the following result relating the nearby cycles on $Sh$, $\widetilde{Sh}$, and ${\bf M}^{\rm loc}$, which follows immediately from the above remarks and Theorem \ref{RPsi_properties} below (cf. \cite{GN}):

\begin{lemma} \label{RPsi_Sh_vs_local_model}
There are canonical isomorphisms 
$$
\varphi^*R\Psi^{Sh}(\Ql) = R\Psi^{\widetilde{Sh}}(\Ql) = \psi^*R\Psi^{{\bf M}^{\rm loc}}(\Ql).
$$
Moreover, $R\Psi^{Sh}(\Ql)$ is constant on each stratum $Sh_w$, and if $\Phi_{\mathfrak p} \in {\rm Gal}(\overline{\Q}_p/E)$ is a geometric Frobenius element, then for any elements $x \in Sh_w(k_r)$ and $x_0 \in {\bf M}^{\rm loc}_w(k_r)$, we have
$$
{\rm Tr}^{ss}(\Phi^r_{\mathfrak p}, R\Psi^{Sh}_x(\Ql)) = {\rm Tr}^{ss}(\Phi^r_{\mathfrak p}, R\Psi_{x_0}^{{\bf M}^{\rm loc}}(\Ql)).
$$
\end{lemma}

Here we use the notion of semi-simple trace, which is explained below in $\S \ref{bad_problems}$.  The above lemma plays a key role in the determination of the semi-simple local zeta function in $\S$\ref{ss_for_simple}.

\subsection{The Genestier-Ng\^{o} comparison with $p$-rank}  
We assume in this section that $Sh$ is the Siegel modular variety with Iwahori-level structure from $\S \ref{Siegel}$.

Recall that any $n$-dimensional abelian variety $A$ over an algebraically closed field $k$ of characteristic $p$ has 
$$
\# A(k) = p^j,
$$
for some integer $0 \leq j \leq n$.  The integer $j$ is called the $p$-{\em rank} of $A$.  Ordinary abelian varieties are those whose $p$-rank is $n$, the largest possible.  The $p$-rank is constant on isogeny classes, and therefore it determines a well-defined function on the set of geometric points $Sh_{\bar{\F}_p}$.  The level sets determine the {\em stratification by $p$-rank}.

It is natural to ask how this stratification relates to the KR-stratification.
In \cite{GN}, Genestier and Ng\^{o} have very elegantly derived the relationship using local models and work of de Jong 
\cite{deJ}.  As they point out, their theorem yields interesting results even in the case of Siegel modular varieties having {\em  good reduction} at $p$: they derive a short and beautiful proof that the ordinary locus in such Shimura varieties is open and dense in the special fiber (comp. \cite{W1}).

To state their result, we define for $w \in \widetilde{W}({\rm GSp}_{2n})$ an integer $r(w)$ as follows.  
Its image $\overline{w}$ in the finite Weyl group $W({\rm GSp}_{2n})$ is a permutation of the set $\{1, \cdots, 2n \}$ commuting with the involution $i \mapsto 2n+1-i$.  The set of fixed points of $\overline{w}$ is stable under the involution, and therefore has even cardinality (the involution is without fixed-points).  Define 
$$
2r(w) = \#\{ \mbox{fixed points of $\overline{w}$} \}.
$$
 
\begin{theorem}[Genestier-Ng\^{o} \cite{GN}]
The $p$-rank is constant on each $KR$-stratum $Sh_w$.  More precisely, the $p$-rank of a point in $Sh_w$ is the integer $r(w)$.
\end{theorem}

\begin{cor} [\cite{GN}] \label{GN_ordinary_locus}
The ordinary locus in $Sh_{\F_p}$ is precisely the union of the KR-strata indexed by the 
translation elements 
in ${\rm Adm}(\mu)$, that is, the elements $t_\lambda$, for $\lambda \in W\mu$.  
Moreover, the ordinary locus is dense and open in $Sh_{\F_p}$.
\end{cor}

\begin{proof}
By the theorem, the $p$-rank is $n$ on $Sh_w$ if and only $r(w) = n$; writing $w = t_\lambda \overline{w}$, this is equivalent to $\overline{w} = 1$.  We conclude the first statement by noting that $w \in {\rm Adm}(\mu) \Rightarrow \lambda \in W\mu$.  

Finally, the union of the strata ${\bf M}_{t_\lambda}^{\rm loc}$ for $\lambda \in W\mu$ is clearly dense and open in ${\bf M}^{\rm loc}_{\F_p}$. 
\end{proof} 

\begin{rem}
It should be noted that in \cite{GN} the local model, and thus the KR-stratification, is defined in terms of the ``cohomology'' local model diagram, whereas here everything is stated using the ``homology'' version.  Furthermore, in \cite{GN} the ``standard'' lattice chain is ``opposite'' from ours, so that an element $w \in \widetilde{W}({\rm GSp}_{2n})$ is used to index a double coset $\bar{I} w \bar{I}/\bar{I}$, where $\bar{I}$ is an ``opposite'' Iwahori subgroup.  
Nevertheless, our conventions and those of \cite{GN} yield {\em the same answer}, that is, the $p$-rank on $Sh_w$ is given by $r(w)$ is both cases.  This may be seen by using the comparison between ``homology'' and ``cohomology'' local model diagrams in Prop. \ref{homology_vs_cohomology_local_models}, and by imitating the proof of \cite{GN} with our conventions in force.
\end{rem}

\subsection{The smooth locus of $Sh_{\F_p}$} \label{smooth_locus}

Also, in \cite{GN} one finds the proof of the following related fact.

\begin{prop}[Genestier-Ng\^o \cite{GN}] \label{GN_comparison}
The smooth locus of $Sh_{\F_p}$ is the union of the KR-strata indexed by $t_\lambda$, for $\lambda \in W\mu$ (in particular the smooth locus agrees with the ordinary locus).
\end{prop}

\subsubsection{The geometric proof of \cite{GN}}
The crucial observation is that any stratum $Sh_w$, where $w$ is {\em not} a translation element, is contained in the singular locus.  Genestier and Ng\^{o} deduce this by showing that for such $w$,
\begin{equation} \label{eq:RPsi_nontrivial}
{\rm Tr}^{ss}(\Phi^r_{\mathfrak p}, R\Psi^{Sh}_w(\Ql)) \neq 1,
\end{equation}
which by the general geometric principle explained below, shows that $w$ is singular.
Now (\ref{eq:RPsi_nontrivial}) itself is proved by combining the main theorems of \cite{HN1} and \cite{H2}, and by taking into account Lemma \ref{RPsi_Sh_vs_local_model}.

Here is the geometric principle implicit in \cite{GN} and a sketch of the proof from \cite{GN}.  We will use freely the material from sections \ref{bad_problems} and \ref{nearby_cycles} below.
\begin{lemma} Let $S = (S,s,\eta)$ be a trait, with $k(s) = {\mathbb F}_q$ a finite field.
Suppose $M \rightarrow S$ is a finite type flat model with $M_\eta$ smooth.  Then $x \in M({\mathbb F}_{q^r})$ is a smooth point of $M_{\bar s}$ only if ${\rm Tr}^{ss}(\Phi^r_{q}, R\Psi_x^M(\Ql)) = 1$.
\end{lemma}
\begin{proof}  
Let $M' \subset M$ be the open subscheme obtained by removing the singular locus of the special fiber of $M$.  We see that $M' \rightarrow S$ is smooth (since $M' \rightarrow S$ is flat of finite-type, it suffices to check the smoothness fiber by fiber, and by construction $M'_\eta = M_\eta$ and $M'_s$ are both smooth).  Now we invoke the general fact (Theorem \ref{RPsi_properties}) that for smooth models $M'$, $R\Psi^{M'}(\Ql) \cong \Ql$, the constant sheaf on the special fiber.  This implies that the semi-simple trace of nearby cycles at $x \in M'_s$ is 1.
\end{proof}

\noindent {\em Proof of Proposition \ref{GN_comparison}.}  We consider the stratum $Sh_w$, or equivalently, $\Mloc_w$, for $w \in {\rm Adm}(\mu)$.  We recall that $\Mloc = M_{-w_0\mu}$ and the stratum $\Mloc_w$ is the Iwahori-orbit indexed by $x := w^{-1}$, contained in $M_{-w_0\mu}$.  For such an $x$, 
we have from \cite{HN1}, \cite{H2}, \cite{HP} an explicit formula for the semi-simple trace of Frobenius on nearby cycles at $x$
$$
{\rm Tr}^{ss}(\Phi_q, R\Psi_x^{{\bf M}^{\rm loc}}(\Ql)) = (-1)^{\ell(t_\mu) + \ell(x)} \, R_{x,t_{\lambda(x)}}(q),
$$
where $\lambda(x)$ is the translation part of $x$ ($x = t_{\lambda(x)}w$, for $\lambda(x) \in X_*$, and $w \in W_0$), and $R_{x,y}(q)$ is the Kazhdan-Lusztig $R$-polynomial.  This polynomial can be computed explicitly, but we need only the fact that it is always a polynomial in $q$ of degree $\ell(y) - \ell(x)$.  It follows from this and the above lemma that whenever $x$ corresponds to a stratum of codimension $\geq 1$, every point of that stratum is singular. 

\subsubsection{A combinatorial proof}
There is however a more elementary way to proceed: we prove below that every codim $\geq 1$ stratum in ${\bf M}^{\rm loc}_{\F_p}$ is contained at least two irreducible components.  The same goes for $Sh_{\F_p}$, proving the proposition.  (There is even a third proof of the proposition, given in \cite{GH}.)

\begin{prop} Let $\mu$ be minuscule.  For any $x \in {\rm Adm}(\mu)$ of codimension 1, there exist exactly two distinct translation elements $\lambda_1, \lambda_2$ in $W\mu$ such that $x \leq t_{\lambda_i}$ (for $i = 1,2$).  Thus, any codimension 1 KR-stratum in the special fiber of a Shimura variety $Sh$ with Iwahori-level structure at $p$ is contained in exactly two irreducible components.   
\end{prop}
\begin{proof}
We give a purely combinatorial proof.  Suppose $x \in {\rm Adm}(\mu)$ has codimension 1.  We have $x < t_\nu$, for some $\nu \in W\mu$.  By properties of the Bruhat order, there exists an affine reflection $s_{\beta + k}$, where $\beta$ is a $B$-positive root, such that $x = t_\nu s_{\beta + k}$.  Since $s_{\beta + k} = t_{-k\beta^\vee}s_\beta$, this means 
$x = t_{\nu -k\beta^\vee}s_\beta$.  The translation part must lie in ${\rm Perm}(\mu) \cap X_* = W\mu$, hence we have $x = t_\lambda s_\beta$, for $\lambda \in W\mu$.  By comparing lengths, we have $t_\lambda s_\beta < t_\lambda$ and $t_\lambda s_\beta < s_\beta t_\lambda s_\beta = t_{s_\beta\lambda}$.  We claim that 
$\langle \beta, \lambda \rangle < 0$.  Indeed, writing $\epsilon_\beta \in [-1,0)$ for the infimum of the set $\beta({\bf a})$ (recalling that our base alcove ${\bf a}$ is contained in the $\bar{B}$-positive chamber) , we have
\begin{align*}
t_{\lambda}s_\beta < t_\lambda &\Leftrightarrow \mbox{${\bf a}$ and $t_{-\lambda}{\bf a}$ are on opposite sides of the hyperplane $\beta = 0$} \\
         &\Leftrightarrow \beta(-\lambda + {\bf a}) \subset [0, \infty) \\
         &\Leftrightarrow -\langle \beta, \lambda \rangle + (\epsilon_\beta,0) \subset [0, \infty) \\ 
         &\Leftrightarrow \langle \beta, \lambda \rangle < 0. \\
\end{align*}
We see that $s_\beta\lambda \neq \lambda$, and so $x$ precedes at least the two distinct translation elements $t_\lambda$ and $t_{s_\beta\lambda}$ in ${\rm Adm}(\mu)$.  It remains to prove that these are the only such translation elements.  So suppose now that $t_\lambda s_\beta < t_{\lambda'}$, where $\lambda' \in W\mu$; we will show that $\lambda' \in \{ \lambda, s_\beta\lambda \}$.  As above, there is an affine reflection $s_{\alpha + n}$ such that $t_\lambda s_\beta = t_{\lambda'}s_{\alpha + n} = t_{\lambda' - n\alpha^\vee}s_\alpha$, where $\alpha$ is $B$-positive.  We see that $\alpha = \beta$, and 
$\lambda' - n\beta^\vee = \lambda$.  Thus, $\lambda'$, $\lambda$, and $s_\beta\lambda$ all lie on the line $\lambda + {\mathbb R}\beta^\vee$.  Since all elements in $W\mu$ are vectors with the same Euclidean length, this can only occur if $\lambda' \in \{ \lambda, s_\beta\lambda \}$. 
\end{proof}

\section{Langlands' strategy for computing local $L$-factors}

The well-known general strategy for computing the local $L$-factor at $p$ of a Shimura variety in terms of automorphic 
$L$-functions 
is due to the efforts of many people, beginning with Eichler, Shimura, Kuga, Sato, and Ihara, and reaching its final conjectural form with Langlands, Rapoport, and Kottwitz.

Let us fix a rational prime $p$, and a compact open subgroup $K_p \subset G(\Q_p)$ at $p$; we consider the Shimura variety $Sh({\bf G},h)_{\bf K}$ as in $\S \ref{PEL_type_data}$.

Roughly, the method of Langlands is to start with a cohomological definition of the local factor of the Hasse-Weil zeta function for $Sh({\bf G},h)_{\bf K}$, and express its logarithm via the Grothendieck-Lefschetz trace formula as a certain sum of orbital integrals for the group ${\bf G}({\mathbb A})$.  
This involves both a process of counting points (with ``multiplicity'' -- the trace of the correspondence on the stalk of an appropriate sheaf), and then a ``pseudo-stabilization'' like that done to stabilize the geometric side of the Arthur-Selberg trace formula (we are ignoring the appearance of endoscopic groups other than ${\bf G}$ itself in this stage).  At this point, we can apply the Arthur-Selberg stable trace formula and express the sum as a trace of a function on automorphic representations appearing in the discrete part of 
$L^2({\bf G}({\mathbb Q})\backslash {\bf G}({\mathbb A}))$.  This equality of traces implies a relation like that in Theorem \ref{HN3_main_thm} below.  

More details on the general strategy as well as the precise conjectural description of $IH^i(Sh \times_E \overline{\Q}_p, \, \Ql)$ in terms of automorphic representations of ${\bf G}$ can be found in \cite{Ko90}, 
\cite{Ko92b}, and \cite{BR}.  These sources discuss the case of good reduction at ${\mathfrak p}$.

Some details of the analogous strategy in case of bad reduction will be given below in $\S \ref{ss_for_simple}$.  Let us explain more carefully the relevant definitions.  

\subsection{Definition of local factors of the Hasse-Weil Zeta function}

Let ${\mathfrak p}$ denote a prime of a number field ${\bf E}$, lying over $p$.  Let $X$ be a smooth $d$-dimensional variety over ${\bf E}$.  We define the Hasse-Weil zeta function (of a complex variable $s$) as an Euler product
$$
Z(s,X) = \prod_{\mathfrak p} Z_{\mathfrak p}(s,X),
$$
where the local factors are defined as follows.

\begin{definition} \label{zeta_fcn} The factor $Z_{\mathfrak p}(s,X)$ is defined to be
$$
\prod_{i=0}^{2d} {\rm det}(1 - N{\mathfrak p}^{-s} \Phi_{\mathfrak p} \, ; \, H^i_c(X \times_E \overline{\Q}_p, \, \overline{\Q}_\ell)^{\Gamma^0_{\mathfrak p}})^{(-1)^{i+1}},
$$
where 
\begin{itemize}
\item $\ell$ is an auxiliary prime, $\ell \neq p$;
\item $\Phi_{\mathfrak p}$ is the inverse of an arithemetic Frobenius element for the extension 
${\bf E}_{\mathfrak p}/{\Q}_p$;
\item $N{\mathfrak p} = {\rm Norm}_{{\bf E}_{\mathfrak p}/\Q_p}{\mathfrak p}$;
\item $\Gamma^0_{\mathfrak p} \subset \Gamma_{\mathfrak p} := {\rm Gal}(\overline{\Q}_p/{\bf E}_{\mathfrak p})$ is the inertia subgroup.
\end{itemize}
\end{definition}

It is believed that $Z(s,X)$ has good analytical properties (it should satisfy a functional equation and have a meromorphic analytic continuation on $\C$) and that analytical invariants (e.g. residues, special values, orders of zeros and poles) carry important arithmetic information about $X$.  Moreover, it is believed that the local factor is indeed independent of the choice of the auxiliary prime $\ell$.  At present these remain only guiding principles, as very little has been actually proved.  
As Taniyama originally proposed, a promising strategy for establishing the functional equation is to express the zeta function as a product of automorphic $L$-functions, whose analytic properties are easier to approach.  The hope that this can be done is at the heart of the Langlands program.

\subsection{Definition of local factors of automorphic $L$-functions} \label{local_L_fcn}

Let $\pi_p$ be an irreducible admissible representation of $G(\Q_p)$.  Let us assume that the local Langlands conjecture holds for the group $G(\Q_p)$.  Then associated to $\pi_p$ is a local Langlands parameter, that is, a homomorphism
$$
\varphi_{\pi_p}: W_{\Q_p} \times {\rm SL}_2(\C) \rightarrow \, ^LG,
$$
where $W_{\Q_p}$ is the Weil group for $\Q_p$ and, letting $\widehat{G}$ denote the Langlands dual group over $\C$ associated to $G_{\Q_p}$, the $L$-group is defined to be $^LG = W_{\Q_p} \ltimes \widehat{G}$.  Let $r = (r,V)$ be a rational representation of $^LG$.  The local $L$-function attached to $\pi_p$ is defined using $\varphi_{\pi_p}$ and $r$ as follows.

\begin{definition} \label{def_L_fcn}  We define $L(s, \pi_p, r)$ to be
$$
{\rm det}(1 - p^{-s}\, r \, \varphi_{\pi_p}(\Phi_p \times \tiny{\begin{bmatrix} p^{-1/2} & 0 \\ 0 & p^{1/2} \end{bmatrix}}); \, ({\rm ker} N)^{\Gamma^0_p})^{-1},
$$
where
\begin{itemize}
\item $\Phi_p \in W_{\Q_p}$ is the inverse of the arithmetic Frobenius for $\Q_p$;
\item $N$ is a nilpotent endomorphism on $V$ coming from the action of ${\mathfrak sl}_2$ on the representation $r \varphi_{\pi_p}$, namely $N := d(r \varphi_{\pi_p})(1 \times \tiny{\begin{bmatrix} 0 & 1 \\ 0 & 0 \end{bmatrix}})$; equivalently, $N$ is determined by $r\varphi_{\pi_p}(1 \times \tiny{\begin{bmatrix} 1 & 1 \\ 0 & 1 \end{bmatrix}}) = {\rm exp}(N)$;
\item the action of inertia $\Gamma^0_p \subset W_{\Q_p}$ on ${\rm ker}(N)$ is the restriction of $r  \varphi_{\pi_p}$ to $\Gamma^0_p \times {\rm id} \subset W_{\Q_p} \times {\rm SL}_2(\C)$.
\end{itemize}
\end{definition}

\begin{rem}
If $\pi_p$ is a spherical representation, then $(\rm ker N)^{\Gamma^0_p} = V$ and therefore in that case the local factor takes the more familiar form
$$
{\rm det}(1 - p^{-s}r({\rm Sat}(\pi_p)); \, V)^{-1},
$$
where ${\rm Sat}(\pi_p)$ is the {\rm Satake parameter} of $\pi_p$; see \cite{Bo}, \cite{Ca}.
\end{rem}

Any irreducible admissible representation $\pi$ of ${\bf G}(\A)$ has a tensor factorization $\pi = \otimes_v \pi_v$ ($v$ ranges over all places of $\Q$) where $\pi_v$ is an admissible representation of the local group ${\bf G}(\Q_v)$.  If $v= \infty$, there is a suitable definition of $L(s, \pi_\infty, r)$  (see \cite{Ta}).  We then define the automorphic $L$-function
$$
L(s, \pi, r) = \prod_v L(s, \pi_v, r).
$$

\subsection{Problems in case of bad reduction, and definition of semi-simple local factors} \label{bad_problems}

\subsubsection{Semi-simple zeta function}

 Let us fix ${\mathfrak p}$ and set $E = {\bf E}_{\mathfrak p}$.
In the case where $X$ possesses an integral model over $\cO_E$, one can study the local factor by reduction modulo ${\mathfrak p}$.  In the case of good reduction (meaning this model is smooth over $\cO_E$), the inertia group acts trivially and the cohomology of $X \times_E \overline{\Q}_p$ can be identified with that of the special fiber $X \times_{k_E} \overline{\F}_p$  
($k_E$ being the residue field of $E$).  
By the Grothendieck-Lefschetz trace formula, the local zeta function then satisfies the familiar identity
$$
{\rm log}(Z_{\mathfrak p}(s,X)) = \sum_{r=1}^{\infty} \#X(k_{E,r})\, \dfrac{N{\mathfrak p}^{-rs}}{r},
$$
where $k_{E,r}$ is the degree $r$ extension of $k_E$ in its algebraic closure $k$.
In the case of bad reduction, inertia can act non-trivially and the smooth base-change theorems of $\ell$-adic cohomology no longer apply in such a simple way.  Following Rapoport \cite{R1}, we bypass the first difficulty by working with the {\em semi-simple} local zeta function, defined below.  The second difficulty forces us to work with nearby cycles (see $\S 10$): if $X$ is a proper scheme over $\cO_E$, then there is a $\Gamma_{\mathfrak p}$-equivariant isomorphism
$$
H^i(X \times_E \overline{\Q}_p, \, \overline{\Q}_\ell) = H^i( X \times_{k_E} \overline{k}_E, \, R\Psi(\overline{\Q}_\ell)),
$$
so that the Grothendieck-Lefschetz trace formula gives rise to the problem of ``counting counts 
$x \in X(k_{E,r})$ with multiplicity'', i.e., to computing the semi-simple trace on the stalks of the complex of nearby cycles:
$$
{\rm Tr}^{ss}(\Phi^r_{\mathfrak p}, R\Psi(\overline{\Q}_\ell)_x),
$$
in order to understand the semi-simple zeta function.

How do we define semi-simple trace and semi-simple zeta functions?  We recall the discussion from \cite{HN1}.  Suppose $V$ is a finite-dimensional continuous $\ell$-adic representation of the Galois group $\Gamma_{\mathfrak p}$.  Grothendieck's local monodromy theorem states that this representation is necessarily {\em quasi-unipotent}: there exists a finite index subgroup $\Gamma_1 \subset \Gamma^0_{\mathfrak p}$ such that $\Gamma_1$ acts unipotently on $V$.  Hence there is a finite increasing $\Gamma_{\mathfrak p}$-invariant filtration $V_\bullet = (0 \subset V_1 \subset \cdots \subset V_m = V)$ such that $\Gamma^0_{\mathfrak p}$ acts on $\oplus_k gr_k V_\bullet$ through a finite quotient.  Such a filtration is called {\em admissible}.  Then we define
$$
{\rm Tr}^{ss}(\Phi_{\mathfrak p}; V) = \sum_k {\rm Tr}(\Phi_{\mathfrak p}; gr_k(V_\bullet)^{\Gamma^0_{\mathfrak p}}).
$$
This is independent of the choice of admissible filtration.   Moreover, the function $V \mapsto {\rm Tr}^{ss}(\Phi_{\mathfrak p}, V)$ factors through the Grothendieck group of the category of $\ell$-adic representations $V$, and using this one proves that it extends naturally to give a ``sheaf-function dictionary'' \`{a} la Grothendieck: a complex $\cF$ in the ``derived'' category 
$D^b_c(X \times_\eta s, \, \overline{\Q}_\ell)$, 
\footnote{ This is the ``derived category of $\overline{\Q}_{\ell}$-sheaves'' -- although this is somewhat misleading terminology (see \cite{KW} for a detailed account)  -- on $X \times_{k_E} \overline{k_E}$ equipped with a continuous action of $\Gamma_{\mathfrak p}$ compatible with the action of its quotient ${\rm Gal}(\overline{k_E}/k_E)$ on $X \times \overline{k_E}$.  See $\S \ref{nearby_cycles}$.}
gives rise to the $\overline{\Q}_\ell$-valued function 
$$
x \mapsto {\rm Tr}^{ss}(\Phi^r_{\mathfrak p}, \, \cF_x) 
$$
on $X(k_{E,r})$.  Furthermore, the formation of this function is compatible with the pull-back and proper-push-forward operations on the derived catgegory, and a Grothendieck-Leftschetz trace formula holds.  (For details, see \cite{HN1}.)  We can then {\em define} the semi-simple local zeta function 
$Z^{ss}_{\mathfrak p}(s,X)$ by the identity

\begin{definition} \label{ss_zeta}
$$
{\rm log}(Z^{ss}_{\mathfrak p}(s, X)) = \sum_{r = 1}^{\infty} \big( \sum_i (-1)^i \, 
{\rm Tr}^{ss}(\Phi^r_{\mathfrak p}; \, H^i_c(X \times_E \overline{\Q}_p, \, \overline{\Q}_\ell)) \big) ~ \dfrac{N{\mathfrak p}^{-rs}}{r}.
$$
\end{definition}

\begin{rem} \label{unipotence_remark} Note that in the case where $\Gamma^0_{\mathfrak p}$ acts unipotently (not just quasi-unipotently) on the cohomology of the generic fiber, then we have
$$
{\rm Tr}^{ss}(\Phi^r_{\mathfrak p}; \, H^i_c(X \times_E \overline{\Q}_p, \, \Ql)) =  
 {\rm Tr}(\Phi^r_{\mathfrak p}; \, H^i_c(X \times_E \overline{\Q}_p, \, \Ql)).
$$
As we shall see below in $\S \ref{unipotence}$, $\Gamma^0_{\mathfrak p}$ does indeed act unipotently on the cohomology of a proper Shimura variety with Iwahori-level structure at $p$.   The definition of $Z^{ss}$ thus simplifies in that case.
\end{rem}  
 
As before, the global semi-simple zeta function is defined to be an Euler product over all finite places of the local functions: $Z^{ss}(s,X) = \prod_{\mathfrak p} Z^{ss}_{\mathfrak p}(s,X)$.

\subsubsection{Semi-simple local $L$-functions}

Retain the notation of $\S \ref{local_L_fcn}$.  The Langlands parameters $\varphi = \varphi_{\pi_p}$ being used here have the property that the representation $r\varphi$ on $V$ arises from a representation $(\rho, N)$ of the Weil-Deligne group $W'_{\Q_p}$ on $V$, see \cite{Ta}.   In that case we have for $w \in W_{\Q_p}$
\begin{align*}
\rho(w) &= r\varphi(w \times \tiny{\begin{bmatrix} \|w\|^{1/2} & 0 \\ 0 & \|w\|^{-1/2} \end{bmatrix}}) \\
{\rm exp}(N) &= r\varphi(1 \times \tiny{\begin{bmatrix} 1 & 1 \\ 0 & 1 \end{bmatrix}}),
\end{align*}
where $\|w\|$ is the power to which $w$ raises elements in the residue field of $\Q_p$.

Now suppose that via some choice of isomorphism $\overline{\Q}_\ell \cong \C$, the pair $(\rho,N)$ comes from an $\ell$-adic representation $(\rho_\lambda,V_\lambda)$ of $W_{\Q_p}$, by the rule 
$$
\rho_\lambda(\Phi^n\sigma) = \rho(\Phi^n\sigma){\rm exp}(N \, t_\ell(\sigma)),
$$
where $n \in \Z$, $\sigma \in \Gamma^0_p$, and $t_\ell: \Gamma^0_p \rightarrow {\Q}_\ell$ is a nonzero homomorphism (cf. \cite{Ta}, Thm 4.2.1).  

For each $\sigma \in \Gamma^0_p$, $\rho_\lambda(\sigma) = \rho(\sigma) {\rm exp}(N t_\ell(\sigma))$ is the multiplicative Jordan decomposition of $\rho_\lambda(\sigma)$ into its semi-simple and unipotent parts, respectively.  Therefore a vector $v \in V$ is fixed by $\rho_\lambda(\sigma)$ if and only if it is fixed by both $\rho(\sigma)$ and ${\rm exp}(Nt_\ell(\sigma))$.  Thus we have the following result.

\begin{lemma} \label{rho_vs_rho_lambda}  If $(\rho, N, V)$ and $(\rho_\lambda,V_\lambda)$ are the two avatars above of a representation of the Weil-Deligne group $W'_{\Q_p}$, then
$$
({\rm ker} N)^{\Gamma^0_p} = V_\lambda^{\Gamma^0_p}.
$$
Furthermore $\rho$ is trivial on $\Gamma^0_p$ if and only if $\rho_\lambda(\Gamma^0_p)$ acts unipotently on $V_\lambda$, and in that case 
$$({\rm ker} N)^{\Gamma^0_p} = {\rm ker} N = V_\lambda^{\Gamma^0_p}.
$$
\end{lemma}

\begin{cor}  The local $L$-function can also be expressed as
$$
L(s, \pi_p, r) = {\rm det}(1 - p^{-s}\rho_\lambda(\Phi_p); \, V_\lambda^{\Gamma^0_p})^{-1}.
$$
\end{cor}

Note the similarity with Definition \ref{zeta_fcn}.  The representation $\rho_\lambda(\Gamma^0_p)$ being quasi-unipotent, we can define the semi-simple $L$-function as in Definition \ref{ss_zeta}.

\begin{definition} \label{ss_L_fcn}
$$
{\rm log}(L^{ss}(s, \pi_p, r)) = \sum_{r=1}^{\infty}{\rm Tr}^{ss}(\rho_\lambda(\Phi^r_p) ; \, V_\lambda) ~ 
\dfrac{p^{-rs}}{r}.
$$
\end{definition}
(The symbol $r$ occurs with two different meanings here, but this should not cause confusion.)

\begin{rem} \label{unipotence_remark_for_L_fcn}  In analogy with Remark \ref{unipotence_remark}, in the case that $(\rho,N,V)$ has $\rho(\Gamma^0_p) = 1$, in view of Lemma \ref{rho_vs_rho_lambda}, we have
$$
{\rm Tr}^{ss}(\rho_{\lambda}(\Phi^r_p); \, V) = {\rm Tr}(\rho_\lambda(\Phi^r_p); \, V),
$$
and the definition of $L^{ss}$ simplifies.  The dictionary set-up by the local Langlands correspondence asserts that if $\pi_p$ has an Iwahori-fixed vector, then $\rho(\Gamma^0_p)$ is trivial (cf. \cite{W2}).  Moreover, in the situation of parahoric level structure at $p$, the only representations $\pi_p$ which will arise necessarily have Iwahori-fixed vectors.  Hence the simplified definition of $L^{ss}$ will apply in our situation.
\end{rem}

Let $\pi = \otimes_v \pi_v$ be an irreducible admissible representation of ${\bf G}(\A)$.  It is to be hoped that a reasonable definition of $L^{ss}(s, \pi_v, r)$ exists for Archimedean places $v$, and if so we can then define
$$
L^{ss}(s, \pi, r) = \prod_v L^{ss}(s, \pi_v, r).
$$

\section{Nearby cycles} \label{nearby_cycles}

\subsection{Definitions and general facts}

Let $X$ be a scheme of finite type over a finite (or algebraically
closed) field $k$. (The following also works if we assume that $k$ is the fraction field of a discrete valuation ring $R$ with finite residue field, and that $X$ is finite-type over $R$, cf. \cite{Ma}.)  Denote by $\overline{k}$ an algebraic closure of $k$, and by $X_{\overline{k}}$ the base change $X \times_k \overline{k}$.

We denote by $D^b_c(X, \Ql)$ the 'derived' category of $\Ql$-sheaves on $X$.
Note that this is not actually the derived category of the category of
$\Ql$-sheaves, but is defined via a limit process. See \cite{BBD} 2.2.14
or \cite{Weil2} 1.1.2, or \cite{KW} for more details. 
Nevertheless, $D^b_c(X, \Ql)$ is a triangulated category which admits the
usual functorial formalism, and which can be equipped with a 'natural'
 $t$-structure having as its core the category of $\Ql$-sheaves.
If $f: X \lto Y$ is a morphism of schemes of finite type over $k$, we have the derived functors 
$f_\ast, \,\, f_!: D^b_c(X, \Ql) \to D^b_c(Y,\Ql)$ and 
 $f^\ast, \,\, f^! : D^b_c(Y, \Ql) \lto D^b_c(X, \Ql)$.  Ocassionally we denote these same functors using the symbols $Rf_\ast$, etc.

Let $(S,s,\eta)$ denote an Henselian trait: $S$ is the spectrum of a complete discrete valuation ring, with special point $s$ and generic point $\eta$.  The key examples for us are $S = {\rm Spec}(\Z_p)$  (the $p$-adic setting) and $S = {\rm Spec}(\F_p[[t]])$ (the function-field setting).  Let $k(s)$ resp. $k(\eta)$ denote the residue fields of $s$ resp. $\eta$.  

We choose a separable closure $\bar{\eta}$ of $\eta$ and define the Galois group $\Gamma = {\rm Gal}(\bar{\eta}/\eta)$ and the inertia subgroup $\Gamma_0 = 
{\rm ker}[{\rm Gal}(\bar{\eta}/\eta) \rightarrow {\rm Gal}(\bar{s}/s)]$, where $\bar{s}$ is the residue field of the normalization $\bar{S}$ of $S$ in 
$\bar{\eta}$.

Now let $X$ denote a finite-type scheme over $S$.
The category $D^b_c(X \times_{s} \eta, \Ql)$ is the category of sheaves
$\cF \in D^b_c(X_{\bar{s}}, \Ql)$ together with a continuous action of 
 $\Gal(\overline{\eta}/\eta)$ which is compatible with the action on $X_{\bar{s}}$.
(Continuity is tested on cohomology sheaves.)

For $\cF \in D^b_c(X_\eta,\Ql)$, we define the {\em nearby cycles sheaf} to be the object in $D^b_c(X \times_s \eta, \Ql)$ given by
$$
R\Psi^X({\cF}) = \bar{i}^*R\bar{j}_* (\cF_{\bar{\eta}}),
$$
where $\bar{i}: X_{\bar{s}} \hookrightarrow X_{\bar{S}}$ and $\bar{j} : X_{\bar{\eta}} \hookrightarrow X_{\bar{S}}$ are the closed and open immersions of the geometric special and generic fibers of $X/S$, and $\cF_{\bar{\eta}}$ is the pull-back of $\cF$ to $X_{\bar{\eta}}$.  

Here we list the basic properties of $R\Psi^X$, extracted from the standard references: \cite{BBD}, \cite{Il}, \cite{SGA7 I}, \cite{SGA7 XIII}.  The final listed property has been proved (in this generality) only recently, and is due to the author and U. G\"{o}rtz \cite{GH}.

\begin{theorem} \label{RPsi_properties}
The following properties hold for the functors
$$
R\Psi : D^b_c(X_\eta, \overline{\Q}_\ell) \rightarrow D^b_c(X \times_s \eta, \, \overline{\Q}_\ell):
$$
\begin{enumerate}
\item[(a)] $R\Psi$ commutes with proper-push-forward: if $f: X \rightarrow Y$ is a proper $S$-morphism, then the canonical base-change morphism of functors to $D^b_c(Y \times_s \eta, \overline{\Q}_\ell)$ is an isomorphism:
$$
R\Psi \, f_*  \,\, \widetilde{\rightarrow} \,\, f_* \, R\Psi;
$$
In particular, if $X \rightarrow S$ is proper there is a ${\rm Gal}(\bar{\eta}/\eta)$-equivariant isomorphism
$$
H^i( X_{\bar{\eta}}, \, \Ql) = H^i(X_{\bar{s}}, \, R\Psi(\Ql)).
$$
\item[(b)]  Suppose $f: X \rightarrow S$ is finite-type but not proper.  Suppose that there is a compactification $j: X \hookrightarrow \overline{X}$ over $S$ such that the boundary $\overline{X}\backslash X$ is a relative normal crossings divisor over $S$.  Then there is a ${\rm Gal}(\bar{\eta}/\eta)$-equivariant isomorphism 
$$
H^i_c( X_{\bar{\eta}}, \, \Ql) = H^i_c(X_{\bar{s}}, \, R\Psi(\Ql)).
$$
\item[(c)] $R\Psi$ commutes with smooth pull-back: if $p : X \rightarrow Y$ is a smooth $S$-morphism, then the base-change morphism is an isomorphism:
$$
p^*  R\Psi \,\, \widetilde{\rightarrow} \,\, R\Psi  p^*.
$$
\item[(d)]  If $\cF \in D^b_c(X, \Ql)$, we define $R\Phi(\cF)$ (the {\em vanishing cycles}) to be the cone of the canonical morphism
$$
\cF_{\bar{s}} \rightarrow R\Psi(\cF_{\eta});
$$
there is a distinguished triangle
$$
\cF_{\bar{s}} \rightarrow R\Psi(\cF_{\eta}) \rightarrow R\Phi(\cF) \rightarrow \cF_{\bar{s}}[1].
$$
If $X \rightarrow S$ is smooth, then $R\Phi(\Ql) = 0$; in particular, $\Ql \cong R\Psi(\Ql)$ in this case;
\item[(e)] $R\Psi$ commutes with Verdier duality, and preserves perversity of sheaves (for the middle perversity);
\item[(f)] For $x \in X_{\bar{s}}$, $R^i\Psi(\Ql)_x = H^i(X_{(\bar{x}) \bar{\eta}}, \Ql)$, where $X_{(\bar{x})\bar{\eta}}$ is the fiber over $\bar{\eta}$ of the strict henselization of $X$ in a geometric point $\bar{x}$ with center $x$.  In particular, the support of $R\Psi(\Ql)$ is contained in the scheme-theoretic closure of $X_\eta$ in $X_{\bar{s}}$; 
\item[(g)] If the generic fiber $X_\eta$ is non-singular, then the complex $R\Psi^X(\Ql)$ is mixed, in the sense of \cite{Weil2}.
\end{enumerate}
\end{theorem}

\begin{rem}
The ``fake'' unitary Shimura varieties $Sh({\bf G},h)_{\bf K}$ discussed in $\S \ref{simple_setup}$ are proper over $\cO_E$ and hence by (a) we can use nearby cycles to study their semi-simple local zeta functions.  
The Siegel modular schemes of $\S \ref{Siegel}$ are not proper over $\Z_p$, so (a) does not apply directly; in fact it is not a priori clear that there is a Galois equivariant isomorphism on cohomology with compact supports as in (b).  In order to apply the method of nearby cycles to study the semi-simple local zeta function, we need such an isomorphism.

\begin{conj} \label{H_c_conjecture} Let $Sh_{K_p}$ denote a model over $\cO_E$ for a PEL Shimura variety with parahoric level structure $K_p$, as in $\S \ref{PEL_parahoric_section}$.  Then the natural morphism
$$
H^i_c(Sh_{K_p} \times_{\cO_E} \overline{k_E}, \, R\Psi(\Ql)) \rightarrow 
H^i_c(Sh_{K_p} \times_{\cO_E} \overline{\Q}_p, \, \Ql)
$$
is an isomorphism.
\end{conj}

In the Siegel case one should be able to prove this by finding a suitably nice compactification, perhaps by adapting the methods of \cite{CF}.

Such an isomorphism would allow us to study the semi-simple local zeta function in the Siegel case by the same approach applied to the ``fake'' unitary case in $\S \ref{ss_for_simple}$.  The local geometric problems involving nearby cycles have already been resolved in \cite{HN1}, and Conjecture \ref{H_c_conjecture} 
encapsulates the remaining geometric difficulty (which is global in nature).  There will be additional group-theoretic problems in applying the Arthur-Selberg trace formula, however, due to endoscopy.
\end{rem}

\subsection{Concerning the inertia action on certain nearby cycles} \label{unipotence}

Let $Sh({\bf G},h)_{\bf K}$ denote a ``fake'' unitary Shimura variety as in $\S \ref{simple_setup}$, where $K_p$ is an Iwahori subgroup.  Suppose $Sh$ is its integral model over $\cO_E = \Z_p$ defined by the moduli problem in $\S \ref{simple_moduli_problem}$.  Our goal in $\S \ref{ss_for_simple}$ is to explain how to identify its $Z^{ss}_{\mathfrak p}$ with a product of semi-simple local $L$-functions (see Theorem 
\ref{HN3_main_thm}).  We first want to justify our earlier claim that the simplified definitions of $Z^{ss}$ and $L^{ss}$ apply to this case.  The key ingredient is a theorem of D. Gaitsgory showing that the inertia action on certain nearby cycles is unipotent.  

We first recall the theorem of Gaitsgory, which is contained in \cite{Ga}.  Let $\lambda$ denote a dominant coweight of $G_{\overline{\Q}_p}$ and consider the corresponding $G(\Q_p[[t]])$-orbit ${\mathcal Q}_\lambda$ in the affine Grassmannian ${\rm Grass}_{\overline{\Q}_p} = G(\overline{\Q}_p\xT)/G(\overline{\Q}_p\xt)$.  Let $M_\lambda$ denote the scheme-theoretic closure of ${\mathcal Q}_\lambda$ in the deformation $M$ of ${\rm Grass}_{\overline{\Q}_p}$ to ${\mathcal Fl}_{\overline{\F}_p}$, from Remark \ref{deformation_remark}.  Let $IC_\lambda$ denote the intersection complex of the closure $\overline{\mathcal Q}_\lambda$.  

\begin{theorem}  [Gaitsgory \cite{Ga}]
The inertia group $\Gamma^0_p$ acts unipotently on $R\Psi^{M_\lambda}(IC_\lambda)$.
\end{theorem}

See also \cite{GH}, $\S 5$, for a detailed proof of this theorem.  We remark that Gaitsgory proves this statement for nearby cycles taken with respect to Beilinson's deformation of the affine Grassmannian of 
$G(\F_p\xt)$ to its affine flag variety, but the same proof applies in the present $p$-adic setting; see loc. cit.

We can apply this to the local model ${\bf M}^{\rm loc} = M_{-w_0\mu}$.  Because the morphisms in the local model diagram are smooth and surjective (Lemma \ref{Genestier_surjectivity}), by taking \cite{GH}, Lemma 5.6 into account, we see that unipotence of nearby cycles on ${\bf M}^{\rm loc}$ implies unipotence of nearby cycles on $Sh$:
\begin{cor}
The inertia group $\Gamma^0_p$ acts unipotently on $R\Psi^{Sh}(\Ql)$.  Consequently, by Theorem \ref{RPsi_properties} (a), it also acts unipotently on $H^i(Sh \times_{\Q_p} \overline{\Q}_p,\, \Ql)$.
\end{cor}

This justifies the assertion made in Remark \ref{unipotence_remark}.  Together with Remark \ref{unipotence_remark_for_L_fcn}, we thus see that the simplified definitions of $Z^{ss}$ and $L^{ss}$ apply in this case, which will be helpful in $\S \ref{ss_for_simple}$ below.

\section{The semi-simple local zeta function for ``fake'' unitary \\ Shimura varieties} \label{ss_for_simple}
 
We assume in this section that $Sh$ is the model from $\S \ref{simple_setup}$: a ``fake'' unitary Shimura variety.  We also assume that $K_p$ is the ``standard'' Iwahori subgroup of ${\bf G}(\Q_p)$, i.e., the subgroup stabilizing the ``standard'' self-dual multichain of $\cO_B \otimes \Z_p$-lattices
$$
\widetilde{\Lambda}_\bullet = \Lambda_\bullet \oplus \Lambda^*_\bullet
$$
from $\S \ref{integral_data}$.

Following the strategy of Kottwitz \cite{Ko92}, \cite{Ko92b}, we will explain how to express $Z^{ss}_{\mathfrak p}(s, Sh)$ in terms of the functions 
$L^{ss}(s, \pi_p, r)$.

There are two equations to be proved.  The first equation is an expression for the semi-simple Lefschetz number
\begin{equation} \label{eq:Lef=sum}
\sum_{x \in Sh(k_r)} {\rm Tr}^{ss}(\Phi^r_{\mathfrak p}, R\Psi_x(\Ql)) = \sum_{\gamma_0} \sum_{(\gamma, \delta)} c(\gamma_0; \gamma,\delta) ~ O_{\gamma}(f^p) ~ TO_{\delta \sigma}(\phi_r).
\end{equation}
The left hand side is termed the semi-simple Lefschetz number ${\rm Lef}^{ss}(\Phi^r_{\mathfrak p})$.  The right hand side has exactly the same form (with essentially the same notation, see below) as \cite{Ko92}, 
p. 442 \footnote{Note that the factor $|{\rm ker}^1(\Q,{\bf G})| = |{\rm ker}^1(\Q, Z({\bf G}))|$ appearing in loc. cit. does not appear here since our assumptions on ${\bf G}$ guarantee that this number is $1$; see loc. cit. $\S 7$.}.  Recall that the twisted orbital integral is defined as
$$
TO_{\delta \sigma}(\phi_r) = \int_{{\bf G}_{\delta \sigma} \backslash {\bf G}(L_r)} \phi_r(x^{-1} \delta \sigma(x)),
$$
with an appropriate choice of measures (and $O_\gamma(f^p)$ is similarly defined).  However in contrast to loc. cit., here the Haar measure on ${\bf G}(L_r)$ is the one giving the standard Iwahori subgroup 
$I_r \subset {\bf G}(L_r)$ volume $1$.

The second equation relates the sum of (twisted) orbital integrals on the right to the spectral side of the Arthur-Selberg trace formula for ${\bf G}$:
\begin{equation} \label{eq:sum=spectral_side}
\sum_{\gamma_0} \sum_{(\gamma, \delta)} c(\gamma_0; \gamma,\delta) ~ O_{\gamma}(f^p) ~ TO_{\delta \sigma}(\phi_r) = 
\sum_{\pi} m(\pi) ~ {\rm Tr}\,\pi \,(f^p  f^{(r)}_p  f_\infty).
\end{equation}
All notation on the right hand side is as in \cite{Ko92b} (cf. $\S 4$) where the ``fake'' unitary Shimura varieties were analyzed in the case that $K_p$ is a {\em hyperspecial maximal compact} subgroup rather than an Iwahori.  In particular, $\pi$ ranges over irreducible admissible representations of ${\bf G}(\A)$ which occur in the discrete part of $L^2({\bf G}(\Q) A_{\bf G}(\R)^0 \backslash {\bf G}(\A))$, with multiplicity $m(\pi)$.  The function $f^p \in C^\infty_c({K}^p\backslash {\bf G}(\A^p_f)/{K}^p)$ is just the characteristic function ${\mathbb I}_{{K}^p}$ of ${K}^p$, whereas $f_\infty$ is the much more mysterious function from \cite{Ko92b}, $\S 1$ \footnote{Up to the sign $(-1)^{{\rm dim}(Sh_E)}$, a pseudo-coefficient of a representation $\pi_\infty^0$ in the packet of discrete series ${\bf G}(\R)$-representations having trivial central and infinitesimal characters.}.  The function $f^{(r)}_p$ is a kind of ``base-change'' of $\phi_r$ and will be further explained below.

The equality (\ref{eq:sum=spectral_side}) comes from the ``pseudo-stabilization'' of its left hand side, similar to that done in \cite{Ko92b}, $\S 4$.  One important ingredient in that is the ``base-change fundamental lemma'' between the test function $\phi_r$ and its ``base-change'' $f^{(r)}_p$; the novel feature here is that the function $\phi_r$ is more complicated than when $K_p$ is maximal compact: it is not a spherical function, but rather an element in the center of an Iwahori-Hecke algebra (see below).  The ``base-change fundamental lemma'' for such functions  is proved in \cite{HN3}, and further discussion will be omitted here.  Finally, after pseudo-stablilization and the fundamental lemma, we apply a simple form of the Arthur-Selberg trace formula, which produces the right hand side of 
(\ref{eq:sum=spectral_side}).  See \cite{HN3} for details. 

Our object here is to explain (\ref{eq:Lef=sum}), following the strategy of \cite{Ko92} which handles the case where 
$K_p$ is maximal compact (the case of good reduction).  The main difficulty is to identify the test function $\phi_r$ 
that appears in the right hand side.

\subsection{Finding the test function $\phi_r$ via the Kottwitz conjecture}

To understand the function $\phi_r$, we will use the full strength of our description of local models in 
$\S \ref{constructing_local_model_diagrams}$, in particular $\S \ref{fake_unitary_local_models}$, and also the material in the appendix $\S \ref{appendix}$.

In (\ref{eq:Lef=sum}), the index $\gamma_0$ roughly parametrizes polarized $n^2$-dimensional $B$-abelian varieties over $k_r$, up to $\overline{\Q}$-isogeny.  The index $(\gamma,\delta)$ roughly parametrizes those polarized $n^2$-dimensional $B$-abelian varieties, up to $\Q$-isogeny, which belong to the $\overline{\Q}$-isogeny class indexed by $\gamma_0$.  (For precise statements, see \cite{Ko92}.)  Therefore, the summand roughly counts (with ``multiplicity'') the elements in $Sh(k_r)$ which belong to a fixed $\Q$-isogeny class.  We will make this last statement precise, and also explain the crucial meaning of ``multiplicity'' here.

Let us fix a polarized $n^2$-dimensional $B$-abelian variety over $k_r$, up to isomorphism: $(A',\lambda',i')$.  We assume it possesses some $K^p$-level structure $\bar{\eta}'$.  We fix once and for all an isomorphism of skew-Hermitian $\cO_B \otimes \A^p_f$-modules
\begin{equation} \label{eq:adelic_Tate_module}
V \otimes \A^p_f = H_1(A',\A^p_f),
\end{equation}
(in the terminology of \cite{Ko92}, $\S 4$).  Since it comes from a level-structure $\bar{\eta}'$, this isomorphism is Galois-equivariant (with the trivial action on $V \otimes \A^p_f$).

Associated to $(A',\lambda', i')$ is also an $L_r$-isocrystal $(H'_{L_r},\Phi)$ as in $\S \ref{appendix}$.  In brief, $H' = H(A')$ is the the $W(k_r)$-dual of $H^1_{\rm crys}(A'/W(k_r))$, and $\Phi$ is the $\sigma$-linear bijection on $H'_{L_r}$ such that 
$p^{-1}H' \supset \Phi H' \supset H'$ (i.e., $\Phi$ is $V^{-1}$, where $V$ is the Verschiebung from Cor. \ref{sigma(Lie)}).  Because of (\ref{eq:adelic_Tate_module}) and the determinant condition, there is also an isomorphism of skew-Hermitian $\cO_B \otimes L_r$-modules
\begin{equation} \label{eq:isocrystal}
V \otimes_{\Q_p} L_r = H'_{L_r}
\end{equation}
which we also fix once and for all.  (See \cite{Ko92}, p. 430.)

From these isomorphisms we construct the elements $(\gamma_0;\gamma,\delta)$ that appear in (\ref{eq:Lef=sum}).  Namely, the absolute Frobenius $\pi_{A'}$ for $A'/k_r$ acting on $H_1(A', \A^p_f)$ induces the automorphism $\gamma^{-1} \in {\bf G}(\A^p_f)$, and the $\sigma$-linear bijection $\Phi$ acting on $H'_{L_r}$ induces the element $\delta \sigma$, for $\delta \in {\bf G}(L_r)$.  The element $\gamma_0 \in {\bf G}(\Q)$ is constructed from $(\gamma, \delta)$ as in \cite{Ko92}, $\S 14$.  The existence of $\gamma_0$ is proved roughly as follows.  By Cor. \ref{sigma(Lie)}, we have $\Phi^{r} = \pi^{-1}_{A'}$ acting on $H'_{L_r}$.  Hence the elements $\gamma_l \in {\bf G}(\Q_l)$ (for $l \neq p$) and $N\delta \in {\bf G}(L_r)$ 
come from $l$-adic (resp. $p$-adic) realizations of the endomorphism $\pi^{-1}_{A'}$ of $A'$.  
Using Honda-Tate theory (and a bit more), one can view $\pi_{A'}^{-1}$ as a semi-simple element $\gamma_0 \in {\bf G}(\Q)$, which is well-defined up to stable conjugacy.  By its very construction, $\gamma_0$ is stably conjugate to $\gamma$, resp. $N\delta$.

Let $(A_\bullet, \lambda, i, \bar{\eta}) \in Sh(k_r)$ be a point in the moduli problem ($\S \ref{simple_moduli_problem}$).  We want to classify those such that $(A_0, \lambda, i)$ is $\Q$-isogenous to $(A', \lambda', i')$.  Let us consider the category
$$
\{ (A_\bullet, \lambda, i, \xi) \}
$$
consisting of chains of polarized $\cO_B$-abelian varieties over $k_r$ (up to $\Z_{(p)}$-isogeny), equipped with a $\Q$-isogeny of polarized $\cO_B$-abelian varieties $\xi: A_0 \rightarrow A'$, defined over $k_r$.  Of course the integral ``isocrystal'' functor $A \mapsto H(A)$ of $\S\ref{appendix}$ is a covariant functor from this category to the category of  $W(k_r)$-free $\cO_B \otimes W(k_r)$-modules in $H'_{L_r}$ equipped with Frobenius and Verschiebung endomorphisms.  In fact, $(A_\bullet, \lambda, i, \xi) \mapsto \xi(H(A_\bullet))$ gives an isomorphism 
$$
\{ (A_\bullet, \lambda, i, \xi) \} \,\, \widetilde{\rightarrow} \,\, Y_p,
$$
where $Y_p$ is the category consisting of all type $(\widetilde{\Lambda}_\bullet)$ multichains of $\cO_B \otimes W(k_r)$-lattices $H_\bullet$ in $H'_{L_r} = V \otimes L_r$, self-dual up to a scalar in $\Q^{\times}$ \footnote{In particular, these multichains are {\em polarized} in the sense of \cite{RZ}, Def. 3.14.}, 
such that for each $i$, $p^{-1}H_i \supset \Phi H_i \supset H_i$, and $\sigma^{-1}(\Phi H_i/ H_i)$ satisfies the determinant condition.  

\subsubsection{Interpreting the determinant condition.}
The final condition on $H_i$ comes from the determinant condition on ${\rm Lie}(A_i)$, and Cor. \ref{sigma(Lie)}.  Let us see what this means more concretely.   By Morita equivalence (see $\S \ref{fake_unitary_local_models}$), a 
type $(\widetilde{\Lambda}_\bullet)$ multichain of $\cO_B \otimes W(k_r)$-lattices $H_\bullet$, self-dual up to a scalar in $\Q^\times$, inside $H'_{L_r} = V \otimes_{\Q_p} L_r$ can be regarded as a complete $W(k_r)$-lattice chain $H^0_\bullet$ in $L_r^n$.  
By working with $H^0_i$ instead of $H_i$, we can work in $L_r^n$ instead of $V \otimes_{\Q_p} L_r$ (which has ${\rm dim}_{L_r} = 2n^2$).  
Recall our minuscule coweight $\mu = (0^{n-d}, (-1)^{d})$ of ${\rm GL}_n$ ($\S \ref{mu}$),and write $\Phi$ for the Morita equivalent $\sigma$-linear bijection of $L_r^n$.   The determinant condition now reads
$$ \Phi H^0_i/ H_i^0 \cong \sigma(k_r)^d,
$$
that is, the relative position of the $W(k_r)$-lattices $H^0_i$ and $\Phi H^0_i$ in $L_r^n$ is given by 
$$
{\rm inv}_K(H^0_i, \Phi H^0_i) = \sigma(\mu(p)) = \mu(p).
$$
(We write $\mu(p)$ in place of $\mu$ to emphasize our convention that coweights $\lambda$ are embedded in ${\rm GL}_n(\Q_p)$ by the rule $\lambda \mapsto \lambda(p)$.)  The same identity holds for $H_i$ replacing $H^0_i$, when we interpret $\mu$ as a coweight for 
the group ${\bf G}(L_r) \subset {\rm Aut}_B(V \otimes_{\Q_p} L_r)$.   
 
By Theorem \ref{RZ_thm_3.16} (and the proof of \cite{Ko92}, Lemma 7.2), we may find $x \in {\bf G}(L_r)$ such that
$$H_\bullet = x\widetilde{\Lambda}_{\bullet,W(k_r)},
$$
where $\widetilde{\Lambda}_{\bullet, W(k_r)} = \widetilde{\Lambda}_\bullet \otimes_{\Z_p} W(k_r)$ is the ``standard'' self-dual multichain of $\cO_B \otimes W(k_r)$-lattices in $V \otimes_{\Q_p} L_r$.

The determinant condition now reads: for every index $i$ in the chain $\widetilde{\Lambda}_\bullet$,
\begin{equation} \label{eq:permissibility_condition}
{\rm inv}_K(\widetilde{\Lambda}_{i,W(k_r)}, \, x^{-1}\delta \sigma(x) \widetilde{\Lambda}_{i, W(k_r)}) = \mu(p).
\end{equation}

Letting $I_r \subset {\bf G}(L_r)$ denote the stabilizer of $\widetilde{\Lambda}_{\bullet,W(k_r)}$, we have the Bruhat-Tits decomposition 
$$
\widetilde{W}({\rm GL}_n \times {\mathbb G}_m) \cong \widetilde{W}({\bf G}) = I_r \backslash {\bf G}(L_r) / I_r.
$$

Equation (\ref{eq:permissibility_condition}) recalls the definition of the $\mu$-permissible set ($\S \ref{generic_and_special_fibers}$). 
The determinant condition can now be interpreted as: 
$$x^{-1} \delta \sigma(x) \in I_r w I_r
$$ 
for some $w \in {\rm Perm}^{\bf G}(\mu)$.  The  equality ${\rm Adm}^{\bf G}(\mu) = {\rm Perm}^{\bf G}(\mu)$ holds (this translates under Morita equivalence to the analogous statement for ${\rm GL}_n$), and therefore we have proved that the determinant condition can now be interpreted as:
$$
x^{-1} \delta \sigma (x) \widetilde{\Lambda}_{\bullet, W(k_r)} \in {\bf M}_\mu(k_r),
$$
where by definition ${\bf M}_\mu(k_r)$ is the set of type $(\widetilde{\Lambda}_\bullet)$ multichains of $\cO_B \otimes W(k_r)$-lattices in $V \otimes L_r$, self-dual up to a scalar in $\Q^\times$, of form 
$$
g \widetilde{\Lambda}_{\bullet, W(k_r)}
$$  
for some $g \in {\bf G}(L_r)$ such that $I_r g I_r =  I_r w I_r$ for an element $w \in {\rm Perm}^{\bf G}(\mu) = {\rm Adm}^{\bf G}(\mu)$.

\smallskip

Now let $I$ denote the $\Q$-group of self-$\Q$-isogenies of $(A',\lambda', i')$.
Our above remarks and the discussion in \cite{Ko92}, $\S 16$ show that there is a bijection from the set of points $(A_\bullet, \lambda, i, \bar{\eta}) \in Sh(k_r)$ such that $(A_0, \lambda, i)$ is $\Q$-isogenous to $(A', \lambda', i')$, to the set $I(\Q) \backslash (Y^p \times Y_p)$, where
\begin{align*}
Y^p &= \{ y \in {\bf G}(\A^p_f)/{K}^p ~ | ~ y^{-1}\gamma y \in {K}^p \},  \\
Y_p &= \{ x \in {\bf G}(L_r)/I_r \, | \, I_r x^{-1}\delta \sigma(x)\widetilde{\Lambda}_{\bullet, W(k_r)} \subset {\bf M}_\mu(k_r) \}.
\end{align*}

\subsubsection{Compatibility of ``de Rham'' and ``crystalline'' maps.}
 
\begin{lemma}  \label{de Rham_crystalline_compatibility}
Fix $w \in {\rm Adm}(\mu)$.
Let $A = (A_\bullet, \lambda, i, \bar{\eta}) \in Sh(k_r)$.   Suppose that $A$ is $\Q$-isogenous to $(A', \lambda', i')$ and that via a choice of $\xi: A_0 \rightarrow A'$ we have $\xi H(A_\bullet) = 
x \widetilde{\Lambda}_{\bullet, W(k_r)}$ as above, for $x \in 
{\bf G}(L_r)/I_r$ \footnote{Note that $x \in {\bf G}(L_r)/I_r$ depends on $\xi$, but its image in $I(\Q) \backslash {\bf G}(L_r)/ I_r$ is independent of $\xi$, and hence the double coset $I_r x^{-1} \delta \sigma(x) I_r$ is well-defined.}.

Then $A$ belongs to the KR-stratum $Sh_w$ if and only if
$$
I_r x^{-1} \delta \sigma(x) I_r = I_r w I_r.
$$
\end{lemma}

\begin{proof}
Recall that the local model ${\bf M}^{\rm loc}$ is naturally identified with the model $M_{-w_0\mu}$ and as such 
its special fiber carries an action of the standard Iwahori-subgroup of 
${\rm GL_n}(k_r \xt)$; further the KR-stratum $Sh_w$ is the set of points which give rise to a point in the Iwahori-orbit indexed by $w^{-1}$ under the ``de Rham'' map $\psi: \widetilde{Sh} \rightarrow \Mloc$.  Let us recall the definition of $\psi$.  We choose any isomorphism
$$
\gamma_\bullet : M(A_\bullet) \,\, \widetilde{\rightarrow} \,\, \widetilde{\Lambda}_{\bullet,k_r}
$$
of polarized $\cO_B \otimes k_r$-multichains.  The quotient $M(A_\bullet) \rightarrow {\rm Lie}(A_\bullet)$ then determines via Morita equivalence a quotient
$$
\dfrac{\cV_{\bullet,k_r}}{t \cV_{\bullet,k_r}} \rightarrow \dfrac{\cV_{\bullet, k_r}}{\cL_\bullet}
$$
for a uniquely determined $k_r\xt$-lattice chain $\cL_\bullet$ satisfying $t\cV_{\bullet, k_r} \subset \cL_\bullet \subset \cV_{\bullet,k_r}$ (see $\S \ref{construction_of_KR_stratification}$).  The ``de Rham'' map $\psi$ sends the point $(A, \gamma_\bullet) \in \widetilde{Sh}(k_r)$ to $\cL_\bullet$.  By definition $A \in Sh_w$ if and only if
$$
{\rm inv}_I(\cL_\bullet, \cV_{\bullet, k_r}) = w,
$$
where the invariant measures the relative position of complete $k_r\xt$-lattice chains in $k_r \xT^n$.  (Note that $w$ is independent of the choice of $\gamma_\bullet$.)  

In the previous paragraph the ambient group is ${\rm GL}_n(k_r\xT)$, the function-field analogue of ${\rm GL}_n(L_r)$.  But it is more natural to work directly with ${\bf G}(L_r)$.  Given $\cL_\bullet$ as above, there exists a unique type $(\widetilde{\Lambda}_\bullet)$ polarized multichain of 
$\cO_B \otimes W(k_r)$-lattices $\widetilde{\cL}_\bullet$ in $V \otimes_{\Q_p} L_r$ such that 
$p \widetilde{\Lambda}_{\bullet, W(k_r)} \subset \widetilde{\cL}_\bullet \subset 
\widetilde{\Lambda}_{\bullet, W(k_r)}$ and the polarized $\cO_B \otimes k_r$-multichain 
$\widetilde{\Lambda}_{\bullet, W(k_r)}/ \widetilde{\cL}_\bullet$ is Morita equivalent to the $k_r$-lattice 
chain $\cV_{\bullet, k_r}/ \cL_\bullet$.  Thus we can think of $\psi$ as the map 
$$
(A, \gamma_\bullet) \mapsto \widetilde{\cL}_\bullet.
$$
Thus, $A \in Sh_w$ if and only if 
$$
{\rm inv}_{I_r}(\widetilde{\cL}_\bullet \,, \, \widetilde{\Lambda}_{\bullet, W(k_r)}) = w.
$$

We have an isomorphism of polarized multichains of $\cO_B \otimes k_r$-lattices
$$
\dfrac{\widetilde{\Lambda}_{\bullet, W(k_r)}}{\widetilde{\cL}_\bullet} = {\rm Lie}(A_\bullet) = \sigma^{-1} \left(\dfrac{x^{-1}\delta \sigma(x) \widetilde{\Lambda}_{\bullet, w(k_r)}}{\widetilde{\Lambda}_{\bullet,W(k_r)}}\right).
$$
The second equality comes from Cor. \ref{sigma(Lie)}; the map sending $(A_\bullet, \xi)$ to the multi-chain $x\widetilde{\Lambda}_{\bullet, W(k_r)}$ may be termed ``crystalline'' -- it is defined using crystalline homology.  This equality shows that the ``de Rham'' and ``crystalline'' maps are compatible (and ultimately rests on Theorem \ref{Oda} of Oda).

Putting these remarks together, we see that $A \in Sh_w$ if and only if
$${\rm inv}_{I_r}(\widetilde{\Lambda}_{\bullet, W(k_r)} \, ,\,\,  x^{-1}\delta \sigma(x) \widetilde{\Lambda}_{\bullet, W(k_r)}) = w,
$$ which completes the proof.
  
\end{proof}

\begin{rem}
The crux of the above proof is the aforementioned compatibility between the ``de Rham'' and ``crystalline'' maps.  This compatibility can be rephrased as the commutativity of the diagram at the end of $\S 7$ in \cite{R2}, when the morphisms there are suitably interpreted.  
\end{rem}

\subsubsection{Identifying $\phi_r$.}

To find $\phi_r$ we need to count the points in the set $I(\Q)\backslash (Y^p \times Y_p)$  with the correct ``multiplicity''. The test function $\phi_r$ in the twisted orbital integral must be such that 
\begin{equation} \label{eq:Lef_on_isogeny_class}
\sum_A {\rm Tr}^{ss}(\Phi^r_{\mathfrak p}; R\Psi^{Sh}_A(\Ql)) = \int_{I(\Q)\backslash ({\bf G}(\A^p_f) \times {\bf G}(L_r))}  f^p(y^{-1} \gamma y) \, \phi_r(x^{-1} \delta \sigma(x)),
\end{equation}
where $A$ ranges over points $(A_\bullet, \lambda, i, \bar{\eta}) \in Sh(k_r)$ such that $(A_0, \lambda, i)$ is $\Q$-isogenous to $(A', \lambda', i')$  (other notation and measures as in \cite{Ko92}, $\S16$, with the exception that here the Haar measure on $G(L_r)$ gives $I_r$ volume $1$).

Now by Lemma \ref{de Rham_crystalline_compatibility}, equation (\ref{eq:Lef_on_isogeny_class}) will hold if $\phi_r$ is a function in the Iwahori-Hecke algebra of $\Ql$-valued functions
$$
\cH_{I_r} = C_c(I_r \backslash {\bf G}(L_r)/I_r)
$$
such that 
$$
\phi_r(I_r w I_r) = {\rm Tr}^{ss}(\Phi^r_p, R\Psi^{M_{-w_0\mu}}_{w^{-1}}(\Ql)),
$$
for elements $w \in {\rm Adm}(\mu)$, and zero elsewhere.  Note that on the right hand side, $w^{-1}$ really represents an Iwahori-orbit in the affine flag variety ${\mathcal Fl}_{\F_p}$ over the {\em function field}.  The nearby cycles are equivariant for the Iwahori-action in a suitable sense, so that the semi-simple trace function is constant on these orbits.  Hence the right hand side is a well-defined element of $\Ql$.  

We can simply {\em define} the function $\phi_r$ by this equality.  But such a description of $\phi_r$ will not be useful unless we can identify it with an explicit function in the Iwahori-Hecke algebra (we need to know its traces on representations with Iwahori-fixed vectors, at least, if we want to make the spectral side of (\ref{eq:sum=spectral_side}) explicit).  This however is possible, due to the following theorem.  This result was conjectured by Kottwitz in a more general form, which inspired Beilinson to conjecture that nearby cycles can be used to give a geometric construction of the center of the affine Hecke algebra for {\em any} reductive group in the function-field setting.  This latter conjecture was proved by Gaitsgory \cite{Ga}, whose ideas were adapted to prove the p-adic analogue in \cite{HN1}.

\begin{theorem}[The Kottwitz Conjecture; \cite{Ga}, \cite{HN1}] \label{RPsi=z_mu}
Let $G = {\rm GL}_n$ or ${\rm GSp}_{2n}$.  Let $\lambda$ be a minuscule dominant coweight of $G$, with corresponding $\Z_p$-model $M_\lambda$ (cf. Remark \ref{deformation_remark}).  Let $\cH_{k_r \xT}$ denote the Iwahori-Hecke algebra $C_c(I_r \backslash G(k_r \xT)/I_r)$.  Let $z_{\lambda,r} \in Z(\cH_{k_r \xT})$ denote the Bernstein function in the center of the Iwahori-Hecke algebra which is associated to $\lambda$.  
Then
$$
{\rm Tr}(\Phi_p^r, R\Psi^{M_\lambda}(\Ql)) = p^{r \, {\rm dim}(M_\lambda)/2} \, z_{\lambda,r}.
$$
\end{theorem}

Let $K_r$ be the stabilizer in ${\bf G}(L_r)$ of $\widetilde{\Lambda}_{0, W(k_r)}$; this is a hyperspecial maximal compact subgroup of ${\bf G}(L_r)$ containing $I_r$.  The Bernstein function $z_{\lambda,r}$ is characterized as the unique element in the center of $\cH_{I_r}$ such that the image of $p^{r \langle \rho, \lambda \rangle} z_{\lambda,r}$ under the Bernstein isomorphism
$$
-*{\mathbb I}_{K_r} : Z(\cH_{I_r}) \,\, \widetilde{\rightarrow} \,\, \cH_{K_r} := 
C_c(K_r \backslash {\bf G}(L_r) / K_r)
$$
is the spherical function $f_{\lambda,r} := {\mathbb I}_{K_r \lambda K_r}$.  Here, we have used the fact that $\lambda$ is minuscule.  Recall that $\rho$ is the half-sum of the $B$-positive roots of $G$, and that ${\rm dim}(M_\lambda)/2 = \langle \rho, \lambda \rangle$.    

The above theorem for $\lambda = -w_0\mu$ implies that 
$$
\phi_r(w) = p^{r \, {\rm dim}(Sh)/2} \, z_{-w_0\mu}(w^{-1}).
$$

Now invoking the identity $z_{\mu}(w) = z_{-w_0\mu}(w^{-1})$  (see \cite{HKP}, $\S 3.2$), we have proved the following result.

\begin{prop} \label{test_function_answer}
Let $z_{\mu,r}$ denote the Bernstein function in the center of the Iwahori-Hecke algebra $\cH_{L_r} = C_c(I_r \backslash {\bf G}(L_r)/ I_r)$ corresponding to $\mu$.  Then the test function is given by
$$\phi_r = p^{r \, {\rm dim}(Sh)/2} \, z_{\mu,r}.
$$
\end{prop} 

\noindent{\em Remarks:} 1) Because $\phi_r$ is central we can define its ``base-change'' function $b(\phi_r) =: f^{(r)}_p$, an element in the center of the Iwahori-Hecke algebra for ${\bf G}(\Q_p)$.  We define the base-change 
homomorphism for centers of Iwahori-Hecke algebras as the unique 
homomorphism $b :Z(\cH_{I_r}) \rightarrow Z(\cH_I)$ which induces, via the Bernstein isomorphism, the usual base-change homomorphism for spherical Hecke algebras $b: \cH_{K_r} \rightarrow \cH_K$ for any special maximal compact $K_r$ containing $I_r$ (setting $I = I_r \cap {\bf G}(\Q_p)$ and $K = K_r \cap {\bf G}(\Q_p)$).  This gives a well-defined homomorphism, independent of the choice of $K_r$.  Moreover, the pair of functions $\phi_r$, $f^{(r)}_p$ have matching (twisted) orbital integrals; see \cite{HN3}.

\noindent 2) We know how $z_{\mu,r}$ (and hence how its base-change $b(z_{\mu,r})$) acts on unramified principal series.   This plays a key role in Theorem \ref{HN3_main_thm} below.  In fact, we have the following helpful lemma.

\begin{lemma} \label{how_z_mu_acts} Let $I$ denote our Iwahori subgroup $K^{\bf a}_p \subset G(\Q_p)$, whose reduction modulo $p$ is $B(\F_p)$.  Suppose $\pi_p$ is an irreducible admissible representation of $G(\Q_p)$ with $\pi^I_p \neq 0$.  Let $r_\mu$ be the irreducible representation of $^L{\bf G}_{\Q_p}$ having extreme weight $\mu$.  Let $d = {\rm dim}(Sh_E)$.  Then
$$
{\rm Tr}\, \pi_p \, (p^{rd/2} \, b(z_{\mu,r})) = {\rm dim}(\pi^I_p) \,\, p^{rd/2} \, {\rm Tr}\big(r_\mu \varphi_{\pi_p}
\big(\Phi^r \times \begin{tiny} \begin{bmatrix} p^{-r/2} & 0 \\ 0 & p^{r/2} \end{bmatrix} \end{tiny} \big)\big). 
$$
\end{lemma}
 
\begin{proof}
Write $G = G(\Q_p)$ and $B = B(\Q_p)$ and suppose that $\pi_p$ is an irreducible subquotient of the normalized unramified principal series representation $i^{G}_{B}(\chi)$, for an unramified quasi-character $\chi: T(\Q_p) \rightarrow \C^\times$.  Suppose $K_r \supset I_r$ (thus also $K \supset I$) is a hyperspecial maximal compact,  and suppose that $\pi_\chi$ is the unique $K$-spherical subquotient of $i^{G}_{B}(\chi)$.  
Suppose $\varphi : W_{\Q_p} \rightarrow \, ^LG$ is the unramified parameter associated to $\pi_\chi$.  Our normalization of the correspondence $\pi_p \mapsto \varphi_{\pi_p}$ is such that for all $r \geq 1$, the element $\varphi(\Phi^r) \in \widehat{G} \rtimes W_{\Q_p}$ can be described as
\begin{equation} \label{eq:scalar} 
\varphi(\Phi^r) = \varphi_{\pi_p}\big(\Phi^r \times \begin{tiny} \begin{bmatrix} p^{-r/2} & 0 \\ 0 & p^{r/2} \end{bmatrix} \end{tiny}\big) = (\chi \rtimes \Phi)^r 
\end{equation}
where $\Phi$ is a {\em geometric} Frobenius element in $W_{\Q_p}$, and where we identify $\chi$ with an element in the dual torus $\widehat{T}(\C) \subset \widehat{G}(\C)$ and take the product on the right in the group $^LG$.  Note that our normalization of the local Langlands correspondence is the one compatible with Deligne's normalization of the reciprocity map in local class field theory, where a uniformizer is sent to a geometric Frobenius element (see \cite{W2}).  
 
Let $f_{\mu, r} = {\mathbb I}_{K_r \mu K_r}$ in the spherical Hecke algebra $\cH_{K_r}$; it is the image of $p^{rd/2}z_{\mu,r}$ under the Bernstein isomorphism 
$$
-* {\mathbb I}_{K_r} \,\, : \,\, Z(\cH_{I_r}) ~ \widetilde{\rightarrow} ~ \cH_{K_r}.
$$
Further, $b(p^{rd/2}z_{\mu,r}) * {\mathbb I}_{K} = b(f_{\mu,r})$.  

Now by \cite{Ko84}, Thm (2.1.3), we know that under the usual left action of $\cH_{K}$, $b(f_{\mu,r})$ acts on $\pi^{K}_\chi$ by the scalar
$$
p^{rd/2} \, {\rm Tr}(r_\mu \varphi(\Phi^r)).
$$ 
Taking (\ref{eq:scalar}) into account along with the well-known fact that elements in $Z(\cH_I)$ act on the entire space $i^G_B(\chi)^I$ by a scalar (see e.g. \cite{HKP}), we are done.
\end{proof}

\begin{rem}
For application in Theorem \ref{HN3_main_thm}, we need the ``$\ell$-adic'' analogue of this lemma, i.e., we need to work with the dual group $\widehat{G}(\Ql)$ instead of $\widehat{G}(\C)$.  For a discussion of how to do this, see \cite{Ko92b}, $\S 1$.
\end{rem}

\subsection{The semi-simple local zeta function in terms of semi-simple $L$-functions} 

The foregoing discussion culminates in the following result from \cite{HN3}, to which we refer for details of the proof.

\begin{theorem} [\cite{HN3}]  \label{HN3_main_thm} Suppose $Sh$ is a simple (``fake'' unitary) Shimura variety.  Suppose $K_p$ is an Iwahori subgroup.  Suppose $r_\mu$ is the irreducible representation of $^L({\bf G}_{{\bf E}_{\mathfrak p}})$ with extreme weight $\mu$, where $\mu$ is the minuscule coweight determined by the Shimura data.  
Then we have
$$Z^{ss}_{\mathfrak p}(s,Sh)=\prod_{\pi_f}L^{ss}(s-{d\over 2},
\pi_p,r_\mu)^{a(\pi_f)\dim(\pi_f^{\bf K})}$$
where the product runs over all admissible representations $\pi_f$ 
of ${\bf G}(\mathbb A_f)$, and the integer number $a(\pi_f)$ is given by
$$a(\pi_f)=\sum_{\pi_\infty\in \Pi_\infty} m(\pi_f\otimes \pi_\infty) {\rm Tr}
\,\pi_\infty(f_\infty),$$
where $m(\pi_f\otimes \pi_\infty)$ is the multiplicity of 
$\pi_f\otimes \pi_\infty$ in $L^2({\bf G}(\Q)A_{\bf G}(\R)^0 \backslash {\bf G}(\mathbb A))$. 
Here $\Pi_\infty$ is the set of admissible representations 
of $G(\R)$ whose central and infinitesimal
characters are trivial.  Also, $d$ denotes the dimension of $Sh_E$.
\end{theorem}

Let us remark that it is our firm belief that this result continues to hold when $K_p$ is a general parahoric subgroup of ${\bf G}(\Q_p)$, but some details are more difficult than the Iwahori case treated in \cite{HN3}, and remain to be worked out. 

Finally, recall our assumptions on $\mathfrak p$ implied that ${\bf E}_{\mathfrak p} = \Q_p$ and $N\mathfrak p = (p)$.  In more general circumstances one or both of these statements will fail to hold, and the result will be a slightly more complicated expression for $Z^{ss}_{\mathfrak p}(s,Sh)$.

\section{The Newton stratification on Shimura varieties over finite fields}

\subsection{Review of the Kottwitz and Newton maps}
As usual $L$ denotes the fraction field of the Witt vectors $W(\overline{\F}_p)$, with Frobenius automorphism $\sigma$.  Let $G$ denote a connected reductive group over $\Q_p$.   Then we have the pointed set $B(G) = B(G)_{\Q_p}$ consisting of $\sigma$-conjugacy classes in $G(L)$.  Let us also assume (only for simplicity) that the connected reductive $\Q_p$-group $G$ is {\em unramified}.  We have the Kottwitz map
$$
\kappa_G: B(G) \rightarrow X^*(Z(\widehat{G})^{\Gamma_p}),
$$
where $\Gamma_p = {\rm Gal}(\overline{\Q}_p/\Q_p)$ (see \cite{Ko97}, $\S 7$) .  We also have the Newton map
$$
\overline{\nu} : B(G) \rightarrow {\mathfrak U}^+
$$
where the notation is as follows.  We choose a $\Q_p$-rational Borel subgroup $B \subset G$ and a maximal $\Q_p$-torus $T$ which is contained in $B$. Then ${\mathfrak U} = X_*(T)_{\R}^{\Gamma_p}$, and  ${\mathfrak U}^+$ denotes the intersection of ${\mathfrak U}$ with the cone of $B$-dominant elements in $X_*(T)_{\R}$.  We call $b \in B(G)$ {\em basic} if $\overline{\nu}_b$ is central, i.e., $\overline{\nu}_b \in X_*(Z)_{\R}$.

Suppose $\lbrace \lambda \rbrace$ is a conjugacy class of one-parameter subgroups of $G$, defined over $\Q_p$.  We may represent the class by a unique cocharacter $\overline{\lambda} \in X_*(T) = X^*(\widehat{T})$ lying in the $B$-positive Weyl chamber of $X_*(T)_{\R}$.  The Weyl-orbit of $\overline{\lambda}$ is stabilized by $\Gamma_p$.  The notion of $B$-dominant being preserved by $\Gamma_p$ (since $B$ is $\Q_p$-rational), we see that $\overline{\lambda}$ is fixed by $\Gamma_p$, hence it belongs to ${\mathfrak U}^+$. Also, restricting $\lbrace \lambda \rbrace$ to $Z(\widehat{G})$ determines a well-defined element $\lambda^\natural \in X^*(Z(\widehat{G})^{\Gamma_p})$.  

We can now define the subset $B(G, \lambda) \subset B(G)$ to be the set of classes $[b] \in B(G)$ such that 
\begin{align*} 
\kappa_G(b) &= \lambda^\natural \\
\overline{\nu}_b &\preceq \overline{\lambda}.
\end{align*}
Here $\preceq$ denotes the usual partial order on ${\mathfrak U}^+$ for which $\nu \preceq \nu'$ if  $\nu' - \nu$ is a nonnegative linear combination of simple relative coroots.

\subsection{Definition of the Newton stratification}

Suppose that the Shimura variety $Sh_{K_p} = Sh({\bf G},h)_{K^pK_p}$ is given by a moduli problem of abelian varieties, and that $K_p \subset {\bf G}(\Q_p)$ is a parahoric subgroup.   Also, assume for simplicity that $E := {\bf E}_{\mathfrak p} = \Q_p$.  Let $G = {\bf G}_{\Q_p}$, which again for simplicity we assume is unramified.  Let $k$ denote as usual an algebraic closure of the residue field $\cO_E/{\mathfrak p} = \F_p$, and let $k_r$ denote the unique subfield of $k$ having cardinality $p^r$.  Let $L_r$ denote the fraction field of the Witt vectors $W(k_r)$.

We denote by ${\bf a}$ the base alcove and suppose $0 \in \overline{\bf a}$ is a hyperspecial vertex.  Let $K_p^{\bf a}$ (resp. $K^{0}_p$) denote the corresponding Iwahori (resp. hyperspecial maximal compact) subgroup of $G(\Q_p)$.  We will let $K_p$ denote a ``standard'' parahoric subgroup, i.e., one such that $K^{\bf a}_p \subset K_p \subset K^0_p$.

We can define a map 
$$
\delta_{K_p}: Sh_{K_p}(k) \rightarrow B(G)
$$
as follows.  A point $A_\bullet = (A_\bullet, \lambda, i, \bar{\eta}) \in Sh_{K_p}(k_r)$ gives rise to a 
$c$-polarized virtual $B$-abelian variety over $k_r$ up to prime-to-$p$ isogeny (cf. \cite{Ko92} $\S 10$), which we denote by $(A, \lambda, i)$.  That in turn determines an $L$-isocystal $(H(A)_{L}, \Phi)$ as in $\S \ref{appendix}$, cf. loc. cit.  It is not hard to see that $\delta_{K_p}$ takes values in the subset 
$B(G, \mu) \subset B(G)$.  In fact, if we choose an isomorphism $(H(A)_{L}, \Phi) = (V \otimes L, \delta \sigma)$ of isocrystals for the group $G(L)$, then $\delta \in G(L)$ satisfies
$$
\kappa_G(\delta) = \mu^\natural,
$$
as follows from the determinant condition $\sigma({\rm Lie}(A_0)) = \Phi H(A_0)/H(A_0)$, for any $A_0$ in the chain $A_\bullet$ (use the argument of \cite{Ko92}, p. 431).  Furthermore, the Mazur inequality
$$
\overline{\nu}_{\delta} \preceq \overline{\mu}
$$
can be proved by reducing to the case where $K_p = K_p^0$ (which has been treated by Rapoport-Richartz \cite{RR} -- see also \cite{Ko03}) \footnote{Note that our assumption that $K_p$ be {\em standard}, i.e., $K_p \subset K_p^0$, was used in the last step, because we want to invoke \cite{RR}.  The results of loc. cit. probably hold for nonspecial maximal compact subgroups, so this assumption is probably unnecessary.}. 

\begin{definition} \label{def_Newton_strata}
We call the fibers of $\delta_{K_p}$ the {\em Newton strata} of $Sh_{K_p}$.  The inverse image of the basic set
$$
\delta^{-1}_{K_p}\big( B(G,\mu) \cap B(G)_{basic}\big).
$$
is called the {\em basic locus} of $Sh_{K_p}$.  
\end{definition}
We denote the Newton stratum $\delta^{-1}_{K_p}([b])$ by  ${\mathcal S}_{K_p,[b]}$, or if $K_p$ is understood, simply by ${\mathcal S}_{[b]}$.

The following conjecture is fundamental to the subject.  It asserts that all the Newton strata are nonempty.

\begin{conj}[Rapoport, Conj. 7.1 \cite{R2}] \label{Newton_strata_are_nonempty}
The map
$$
\delta_{K_p}: Sh_{K_p}(k) \rightarrow B(G,\mu)
$$
is surjective.  In particular, the basic locus is nonempty.
\end{conj}

\begin{rem} \label{delta_in_spherical_case}
Note that ${\rm Im}(\delta_{K^0_p})$ can be interpreted purely in terms of group theory: $[b] \in B(G,\mu)$ lies in the image of $\delta_{K^0_p}$ if and only if for one (equivalently, for all sufficiently divisible) $r \geq 1$, $[b]$ contains an element $\delta \in G(L_r)$ belonging to a triple $(\gamma_0;\gamma, \delta)$, which satisfies the conditions in \cite{Ko90} $\S 2$ (except for the following correction, noted at the end of \cite{Ko92}: under the canonical map $B(G)_{\Q_p} \rightarrow X^*(Z(\widehat{G})^{\Gamma_p})$, the $\sigma$-conjugacy class of $\delta$ goes to $\mu|_{X^*(Z(\widehat{G})^{\Gamma_p})}$ and not its negative), and for which the following four conditions also hold:
\begin{enumerate}
\item[(a)] $\gamma_0 \gamma^*_0 = c^{-1}_0p^{-r}$, where $c_0 \in \Q^\times$ is a $p$-adic unit;
\item[(b)] the Kottwitz invariant $\alpha(\gamma_0;\gamma,\delta)$ is trivial;
\item[(c)] there exists a lattice $\Lambda$ in $V_{L_r}$ such that $\delta \sigma \Lambda \supset \Lambda$;
\item[(d)] $O_\gamma({\mathbb I}_{K^p}) \, TO_{\delta \sigma}({\mathbb I}_{K_r \mu K_r}) \neq 0$.
\end{enumerate} 
Here $K_r \subset {\bf G}(L_r)$ is the hyperspecial maximal compact subgroup such that $K_r \cap {\bf G}(\Q_p) = K^0_p$.

To see this, use \cite{Ko92}, Lemmas 15.1, 18.1 to show that the first three conditions ensure the existence of a $c_0 p^r$-polarized virtual $B$-abelian variety $(A',\lambda', i')$ over $k_r$ up to prime-to-$p$ isogeny, giving rise to $(\gamma_0;\gamma, \delta)$.  In the presence of the first three, the last condition shows that there exists a $k_r$-rational point $(A,\lambda,i, \bar{\eta}) \in Sh_{K^0_p}$ such that $(A,\lambda, i)$ is $\Q$-isogenous to $(A',\lambda',i')$.  Hence $\delta_{K^0_p}((A, \lambda, i, \bar{\eta})) = [\delta]$.

Note that by counting fixed points of $\Phi_{\mathfrak p}^r \circ f$ for any Hecke operator $f$ away from $p$ as in \cite{Ko92}, $\S 16$, we may replace condition (d) with 
\begin{enumerate}
\item[(d')] For some $g \in {\bf G}(\A^p_f)$, we have $O_\gamma({\mathbb I}_{K^p g^{-1} K^p}) \, 
TO_{\delta \sigma}({\mathbb I}_{K_r \mu K_r}) \neq 0$.
\end{enumerate}

Since we may always choose $g = \gamma^{-1}$, we may also replace (d) or (d') with  
\begin{enumerate}
\item[(d'')] $TO_{\delta \sigma}({\mathbb I}_{K_r \mu K_r}) \neq 0$.  
\end{enumerate}
\end{rem}

\medskip

In the Siegel case with $K_p$ hyperspecial, we have the following result of Oort \cite{Oo} (comp. \cite{R2}, Thm. 7.4) which in particular proves Conjecture \ref{Newton_strata_are_nonempty} in that case.

\begin{theorem} [Oort, \cite{Oo}] \label{Oort}
Suppose $Sh_{K^0_p}$ is a Siegel modular variety with maximal hyperspecial level structure $K_p^0$ at $p$.  Then the Newton strata are all nonempty and equidimensional (and the dimension is given by a simple formula in terms of the partially ordered set $B({\rm GSp}_{2n}, \mu)$).
\end{theorem}

In Corollary \ref{Im_independent_of_K_p} below, we shall generalize the part of Theorem \ref{Oort} that asserts the nonemptiness of the Newton strata: in the Siegel case, all Newton strata of $Sh_{K_p}$ are nonempty, when $K_p$ is an arbitrary standard parahoric subgroup.  In the ``fake'' unitary case, we shall later prove in Cor. \ref{basic_locus_nonempty} only that the basic locus is nonempty (again for standard parahorics $K_p$).

\subsection{The relation between Newton strata, KR-strata, and affine Deligne-Lusztig varieties}

Fix a $\sigma$-conjugacy class $[b] \in B(G)$, and fix an element $w \in \widetilde{W}$.  Let $\widetilde{K}^{\bf a}_p \subset G(L)$ denote the Iwahori subgroup such that $\widetilde{K}^{\bf a}_p \cap G(\Q_p) = K_p^{\bf a}$. 

\begin{definition} \label{DL_variety}
We define the {\em affine Deligne-Lusztig variety} \footnote{Strictly speaking, this is only a set, not a variety.  The sets are the affine analogues of the usual Deligne-Lusztig varieties in the theory of finite groups of Lie type.}
$$
X_w(b)_{\widetilde{K}^{\bf a}_p} = \lbrace x \in G(L)/\widetilde{K}^{\bf a}_p ~ | ~ x^{-1} b\sigma(x) \in \widetilde{K}^{\bf a}_p \, w \, \widetilde{K}^{\bf a}_p \rbrace,
$$
and for any dominant coweight $\lambda$, 
$$
X(\lambda, b)_{\widetilde{K}^{\bf a}_p} = \bigcup_{w \in {\rm Adm}(\lambda)} X_w(b)_{\widetilde{K}^{\bf a}_p}.
$$
A similar definition can be made for a parahoric subgroup replacing $K^{\bf a}_p$ (cf. \cite{R2}).
\end{definition}

A fundamental problem is to determine the relations between the Kottwitz-Rapoport and Newton stratifications.  The following shows how this problem is related to the nonemptiness of certain affine Deligne-Lusztig varieties.

\begin{prop} \label{KR_vs_Newton_strata_via_DL}

Let $\mu$ be the minuscule coweight attached to the Shimura data for $Sh = Sh_{K^{\bf a}_p}$.  Suppose $w \in {\rm Adm}(\mu)$.  Then for every $[b] \in {\rm Im}(\delta_{K^0_p})$, we have
$$
X_w(b)_{\widetilde{K}^{\bf a}_p} \neq \emptyset \Leftrightarrow Sh_{w} \cap 
{\mathcal S}_{[b]} \neq \emptyset.
$$
\end{prop}

\begin{proof}
By Remark \ref{delta_in_spherical_case}, for all sufficiently divisible $r \geq 1$, $[b]$ contains an element $\delta \in G(L_r)$ which is part of a Kottwitz triple $(\gamma_0;\gamma,\delta)$ satisfying the conditions (a-d).  We consider such a triple, up to equivalence (we say $(\gamma_0';\gamma',\delta')$ is equivalent to $(\gamma_0;\gamma,\delta)$ if $\gamma_0'$ and $\gamma_0$ are stably-conjugate, $\gamma'$ and $\gamma$ are conjugate, and $\delta'$ and $\delta$ are $\sigma$-conjugate).  Then by the arguments in $\S \ref{ss_for_simple}$ together with \cite{Ko92}, $\S 18,19$, we have the equality
\begin{equation}
\#\lbrace A_\bullet \in Sh_w(k_r) ~ | ~ A_\bullet \rightsquigarrow (\gamma_0;\gamma, \delta) \rbrace = 
({\rm vol}) \, O_\gamma({\mathbb I}_{K^p})  ~ \# X_w(\delta)_{I_r},
\end{equation}
Let us explain the notation.   The notation $A_\bullet \rightsquigarrow (\gamma_0;\gamma,\delta)$ means that $A_\bullet$ gives rise to the equivalence class of $(\gamma_0;\gamma,\delta)$; cf. \cite{Ko92}, $\S 18,19$.  The term ${\rm vol}$ denotes the nonzero rational number 
$|{\rm ker}^1(\Q, {\bf G})| c(\gamma_0;\gamma,\delta)$, where the second term  is the number defined in loc. cit.  Also $I_r = \widetilde{K}^{\bf a}_p \cap G(L_r)$.  

This equality would imply the proposition, if we knew that $O_\gamma({\mathbb I}_{K^p}) \neq 0$.  But this follows from condition (d) in Remark \ref{delta_in_spherical_case}.
\end{proof}

The following result of Wintenberger \cite{Wi} proves a conjecture of Kottwitz and Rapoport in a suitably unramified case (cf. \cite{R2}, Conj. 5.2, and the notes at the end).

\begin{theorem} [Wintenberger] \label{Wintenberger}
Let $G$ be any connected reductive group, defined and quasi-split over $L$.  Suppose $\lbrace \lambda \rbrace $ is a conjugacy class of 1-parameter subgroups, defined over $L$.  Suppose $[b] \in B(G)$ and let $K$ be any standard parahoric subgroup (that is, one contained in a special maximal parahoric subgroup).  Then 
$$X(\lambda, b)_K \neq \emptyset \Leftrightarrow [b] \in B(G,\lambda).$$
\end{theorem}

\begin{cor}  \label{Im_independent_of_K_p}
Let $Sh$ be either a ``fake'' unitary or a Siegel modular variety, as in $\S \ref{PEL_parahoric_section}$.
\begin{enumerate}
\item[(a)] In the ``fake'' unitary case,
for any two standard parahoric sugroups $K'_p \subset K''_p$, we have ${\rm Im}(\delta_{K'_p}) = {\rm Im}(\delta_{K''_p})$.
\item[(b)]  In the Siegel case, we have ${\rm Im}(\delta_{K_p}) = B(G,\mu)$ for every standard parahoric subgroup $K_p$.
\end{enumerate}
\end{cor}

\begin{proof} 
Consider first (a).  We need to prove ${\rm Im}(\delta_{K'_p}) \supset {\rm Im}(\delta_{K''_p})$.    Clearly it is enough to consider the case $K'_p = K^{\bf a}_p$ and $K''_p = K^0_p$.  The natural morphism 
$Sh_{K^{\bf a}_p} \rightarrow Sh_{K^0_p}$ is proper, surjective on generic fibers (look at $\C$-points), and the target is flat (even smooth).  Therefore in the special fiber the morphism is surjective.  
This completes the proof of (a).

Consider now part (b).  In the Siegel case the morphism $Sh_{K^{\bf a}_p} \rightarrow Sh_{K^0_p}$ is still projective with flat image, so the same argument combined with Theorem \ref{Oort} yields the stronger result of (b).  Here is another argument using Theorem \ref{Oort},  Proposition \ref{KR_vs_Newton_strata_via_DL} and Theorem \ref{Wintenberger}.  Let $G = {\rm GSp}_{2n}$ and $\mu = (0^n,(-1)^n)$.  By Oort's theorem, ${\rm Im}(\delta_{K^0_p}) = B(G, \mu)$, so it is enough to prove 
$B(G,\mu) \subset {\rm Im}(\delta_{K^{\bf a}_p})$.   Let $[b] \in B(G,\mu)$.  By Wintenberger's theorem, there exists $w \in {\rm Adm}(\mu)$ such that $X_w(b)_{\widetilde{K}^{\bf a}_p} \neq \emptyset$, which implies the result by Proposition \ref{KR_vs_Newton_strata_via_DL}.  

Note that a similar argument provides an alternative proof for part (a).
\end{proof}

\begin{rem} 
Let $G = {\rm GSp}_{2n}$, and suppose $\mu$ is minuscule.  We can give a proof of Theorem \ref{Wintenberger} in this case using Oort's theorem, as follows.  Let $[b] \in B(G,\mu)$.  By Corollary \ref{Im_independent_of_K_p}, there is a point $A = (A_\bullet, \lambda, i, \bar{\eta})$ such that $\delta_{\widetilde{K}^{\bf a}_p}(A) = [b]$.  Now $A$ belongs to {\em some} KR-stratum $Sh_w$, so $Sh_w \cap {\mathcal S}_{[b]} \neq \emptyset$.  Now Proposition \ref{KR_vs_Newton_strata_via_DL} implies that 
$X_w([b])_{\widetilde{K}^{\bf a}} \neq \emptyset$.
\end{rem}

\subsection{The basic locus is nonempty in the ``fake'' unitary case} \label{basic_locus_in_simple_case} \footnote{This non-emptiness is implicit in both articles of \cite{FM}, and can be justified (indirectly) from their main theorems.  Our object here is only to give a simple direct proof.} 

First consider a ``fake'' unitary variety $Sh$ with $\mu = (0^{n-d}, (-1)^d)$.  We shall consider both hyperspecial and Iwahori level structures.

Recall that the subset ${\rm Adm}(\mu) \subset \widetilde{W}({\rm GL}_n)$ contains a unique minimal element.  To find this element, we consider first the element
$$
\tau_1 = t_{(-1, 0^{n-1})} (1 2 \cdots n) \in \widetilde{W}({\rm GL}_n),
$$
where the cycle $(1 2 \cdots n)$ acts on a vector $(x_1, x_2 \dots, x_n) \in X_*(T) \otimes \R$ by sending it to $(x_n, x_1, \dots, x_{n-1})$.  Note that this element preserves our base alcove
$$
{\bf a} = \lbrace (x_1, \dots, x_n) \in X_*(T) \otimes \R ~ | ~ x_n - 1 < x_1 < \cdots < x_{n-1} < x_n \rbrace.
$$
Hence its $d$-th power
$$
\tau = t_{((-1)^d,0^{n-d})} (1 2 \cdots n)^d \in \widetilde{W}({\rm GL}_n),
$$
is the unique element of $\Omega$ which is congruent modulo $W_{\rm aff}$ to $t_{\mu}$, so is the desired element.  This element is identified with an element in ${\rm GL}_n(L)$ using our usual convention (the vector part is sent to ${\rm diag}((p^{-1})^d, 1^{n-d})$).  Via Morita equivalence, we can view it as an element of ${\bf G}(L)$.  In fact, via (\ref{eq:decomposition}), $\tau$ becomes the element
$$
\tau = (X, p^{-1}\chi^{-1}(X^t)^{-1}\chi),
$$
where by definition $X = {\rm diag}((p^{-1})^d, 1^{n-d}) (1 2 \cdots n)^d$. 

Next consider the Siegel case, where $\mu = (0^n, (-1)^n)$.  The unique element $\tau \in \Omega$ which is congruent to $t_\mu$ modulo $W_{\rm aff}$ is given by 
$$
\tau = t_{((-1)^n,0^n)} (12\cdots 2n)^n \in \widetilde{W}({\rm GSp}_{2n}).
$$

We will show that either case, $\tau \in {\bf G}(L)$ is {\em basic}.  We will also show that $[\tau]$ belongs to the image of $\delta_{K^0_p}$.  By virtue of Corollary \ref{Im_independent_of_K_p}, this will show that the basic locus of $Sh_{K_p}$ is nonempty for every standard parahoric $K_p$.

Let us handle the second statement first.

\begin{lemma} \label{tau_is_in_Im}
 Let $\delta = \tau \in {\bf G}(L)$.  Then there exists a Kottwitz triple $(\gamma_0;\gamma, \delta)$ satisfying the conditions in Remark \ref{delta_in_spherical_case}.  Hence $[\tau] \in {\rm Im}(\delta_{K^0_p})$.
\end{lemma}

\begin{proof}
Consider the ``fake'' unitary case.  Note that $\delta \in {\bf G}(L_n)$ is clearly fixed by $\sigma$, so that 
\begin{equation} \label{eq:Ndelta}
N\delta = \delta^n = {\rm diag}((p^{-d})^n) \times {\rm diag}((p^{d-n})^n),
\end{equation}
where the first factor is the diagonal matrix with the entry $p^{-d}$ repeated $n$ times (similarly for the second factor).
This is the image of a unique element $\gamma_0 \in {\bf G}(\Q)$ under the composition of the inclusion ${\bf G}(\Q) \hookrightarrow {\bf G}(\Q_p)$ and the isomorphism (\ref{eq:decomposition}).  In fact $\gamma_0$ belongs to the center of ${\bf G}$, and thus is clearly an elliptic element in ${\bf G}(\R)$.  Moreover, $\gamma_0 \gamma_0^* = p^{-n}$.  For all primes $l \neq p$, we define $\gamma_l$ to the image of $\gamma_0$ under the inclusion ${\bf G}(\Q) \hookrightarrow {\bf G}(\Q_l)$.  We set $\gamma = (\gamma_l)_l \in {\bf G}(\A^p_f)$.  The resulting triple $(\gamma_0; \gamma, \delta)$ is clearly a Kottwitz triple satisfying the conditions (a-c,d'') of Remark \ref{delta_in_spherical_case} (with $r = n$).  One can check that $\alpha(\gamma_0;\gamma,\delta) = 1$ from the definitions, but in fact this is not necessary, as the group to which $\alpha(\gamma_0;\gamma,\delta)$ belongs is itself trivial in the ``fake'' unitary case (see \cite{Ko92b}, Lemma 2).  Hence by Remark \ref{delta_in_spherical_case}, $\delta$ arises from a point in $Sh_{K^0_p}(k_n)$, i.e., $[\tau]$ is in the image of $\delta_{K^0_p}$.

In the Siegel case, the same argument works, if we let $\delta = \tau \in {\rm GL}_n(L_2)$ and note 
\begin{equation} \label{eq:Ndelta_Siegel}
N\delta = \delta^2 = {\rm diag}((p^{-1})^{2n}).
\end{equation}
\end{proof}

\begin{lemma} \label{tau_is_basic}  The element $\tau \in G(L)$ is basic.
\end{lemma}
 \begin{proof}
We want to use the following special case of the characterization of $\bar{\nu}_b$, for certain $b \in G(L)$: suppose we are given an element $b \in G(L)$ such that for sufficiently divisible $s \in \N$, we may write in the semidirect product $G(L) \rtimes \langle \sigma \rangle$ the identity
$$
(b\sigma)^s = (s\nu)(p) \,\, \sigma^s
$$
for a rational $B$-dominant cocharacter $s\nu: {\bf G}_m \rightarrow Z(G)$ defined over $\Q_p$.  Then in that case, $b$ is basic and
$$
\overline{\nu}_b = \frac{1}{s}(s\nu) \in X_*(T)_{\Q}^{\Gamma_p}.
$$
This follows immediately from the general characterization of $\bar{\nu}_b$ given in \cite{Ko85}, $\S 4.3$.

This characterization applies to the element $\tau$ because of the identities (\ref{eq:Ndelta}) and (\ref{eq:Ndelta_Siegel}).  In fact we see that, in ${\rm GL}_n$, the Newton point of $\tau$ is 
$$
\overline{\nu}_\tau = ((\frac{-d}{n})^n) \in X_*(T({\rm GL}_n))_{\Q}^{\Gamma_p},
$$
which clearly factors through the center of $G$.  In the Siegel case,
$$
\overline{\nu}_\tau = ((\frac{-1}{2})^{2n}) \in X_*(T({\rm GSp}_{2n}))_{\Q}^{\Gamma_p}.
$$ 
Thus $\tau$ is basic in each case.
\end{proof}

We get the following, which of course we already knew in the Siegel case, by Theorem \ref{Oort}.

\begin{cor} \label{basic_locus_nonempty}
Let $Sh$ be a ``fake'' unitary or Siegel modular Shimura variety with level structure given by a standard parahoric subgroup $K_p$.  Then the basic locus of $Sh$ is nonempty.
\end{cor}

\medskip

For (more detailed) information concerning the Newton stratification on some other Shimura varieties, the reader might consult the work of Andreatta-Goren \cite{AG} and Goren-Oort \cite{GOo}.

\section{The number of irreducible components in $Sh_{\overline{\F}_p}$}

Let $Sh$ be a ``fake'' unitary or a Siegel modular variety with Iwahori-level structure as in $\S \ref{PEL_parahoric_section}$.  Recall that ${\bf M}^{\rm loc}_{\F_p}$ is a connected variety with an Iwahori-orbit stratification indexed by the finite subset ${\rm Adm}(\mu) \subset \widetilde{W}$, where $\mu = (0^{n-d}, (-1)^{d})$, resp. $(0^n,(-1)^n)$.  Its irreducible components are indexed by the maximal elements in this set, namely by the translation elements 
$t_\lambda$ for  $\lambda$ belonging to the Weyl-orbit $W\mu$ of the coweight $\mu$.  

It is natural to hope that similar statements apply in some sense to $Sh$.  The varieties $Sh_{\F_p}$ and $\widetilde{Sh}_{\F_p}$ are not geometrically connected (the number of connected components of $Sh_{\overline{\F}_p}$ depends on the choice of subgroup ${K}^p \subset {\bf G}(\A^p_f)$; see below).  Nevertheless the following two Lemmas (\ref{Genestier_surjectivity} and \ref{connectedness_of_fibers}) show that in every connected component of $Sh$, all the KR-strata are nonempty.  

\begin{lemma} [\cite{Ge}, Prop. 1.3.2, in the Siegel case] \label{Genestier_surjectivity}
If $Sh$ is either a ``fake'' unitary or a Siegel modular Shimura variety, and $K_p$ is any standard parahoric subgroup, then the morphism $\psi: \widetilde{Sh} \rightarrow {\bf M}^{\rm loc}_{K_p}$ is surjective.
\end{lemma}

\begin{proof}
First we claim that it suffices to prove the lemma for $K^{\bf a}_p$, our standard Iwahori subgroup.  
Indeed, it is enough to observe that the following diagram commutes
$$
\xymatrix{
\widetilde{Sh}_{K^{\bf a}_p} \ar[r]^{\psi} \ar[d] & {\bf M}^{\rm loc} \ar[d]^p \\
\widetilde{Sh}_{K_p} \ar[r]^{\psi} & {\bf M}^{\rm loc}_{K_p},
}
$$
and that the right vertical arrow $p$ is surjective on the level of KR-strata: every ${\rm Aut}_{K_p}$-orbit in ${\bf M}^{\rm loc}_{K_p}$ contains an element in the image of $p$; this follows from \cite{KR}, Prop. 9.3, 10.6.  See also \cite{Go3}.

Next consider the Siegel case with Iwahori level structure, where this result is due to Genestier, loc. cit.  We briefly recall his argument.  It is easy to see that $\psi$ is surjective on generic fibers, because it is ${\rm Aut}$-equivariant, and the generic fiber of ${\bf M}^{\rm loc}$ is a single orbit under ${\rm Aut}$.  Because $\psi$ is smooth, the complement of its image is an ${\rm Aut}$-invariant, Zariski-closed subset of the special fiber of ${\bf M}^{\rm loc}$.  On the other hand, there is a unique closed (zero-dimensional) ${\rm Aut}$-orbit (denoted here by $\tau^{-1}$) in that fiber which belongs to the Zariski-closure of every other ${\rm Aut}$-orbit, and one can show (by writing down an explicit chain of supersingular abelian varieties) that the point $\tau^{-1}$ belongs to the image of $\psi$.  It follows that $\psi$ is surjective.

In the ``fake'' unitary case with Iwahori-level structure, consider the element $\tau$ from $\S \ref{basic_locus_in_simple_case}$ above.   Then the element $\tau^{-1}$ indexes the unique zero-dimensional ${\rm Aut}$-orbit in ${\bf M}^{\rm loc} = M_{-w_0\mu}$.  By Genestier's argument, it is enough to prove that $\tau^{-1}$ belongs to the image of $\psi$, or equivalently, the stratum $Sh_\tau$ is nonempty.  But we have proved in $\S \ref{basic_locus_in_simple_case}$ that $[\delta] = [\tau]$ belongs to the image of $\delta_{K^0_p}$.  Furthermore, it is obvious that $X_\tau(\tau)_{\widetilde{K}^{\bf a}_p} \neq \emptyset$, and hence by Proposition \ref{KR_vs_Newton_strata_via_DL}, we conclude that $Sh_\tau \cap {\mathcal S}_{[\tau]} \neq \emptyset$.   ({\em Note:}  Using the element 
$\tau \in \widetilde{W}({\rm GSp}_{2n})$, this argument applies just as well to the Siegel case, thus providing an alternative to Genestier's step of finding an explicit chain of abelian varieties in the moduli problem which maps to $\tau^{-1}$.)
\end{proof}

Recall that $Sh$ is not a connected scheme.  In fact, for our two examples, the set of connected components of the geometric generic fiber $Sh_{\overline{\Q}_p}$ carries a simply transitive action by the finite abelian group 
$$
\pi_0 = \Z^+_{(p)} \backslash (\A^p_f)^\times / c(K^p) = \Q^+ \backslash \A^\times_f/ c(K^p)\Z^\times_p  = \Q^\times \backslash \A^\times / c({\bf K}K_\infty),
$$
where $c: {\bf G} \rightarrow {\mathbb G}_m$ is the similitude homomorphism, $K_\infty$ denotes the centralizer of $h_0$ in ${\bf G}(\R)$, and the superscript $+$ designates the positive elements of the set; cf. \cite{Del}.  To see the groups above are actually isomorphic, use the fact that $c(K_\infty) \supset \R^+$ and $c(K_p) = \Z^\times_p$, as one can easily check for each of our two examples.  Fixing isomorphisms $\A^p_f(1) = \A^p_f$ and $\overline{\Q}_p = \C$  once and for all, an element $(A_\bullet, \lambda, i, \bar{\eta}) \in Sh(\overline{\Q}_p)$ belongs to the connected component indexed by $a \in \pi_0$ if and only if the Weil-pairing $(\cdot, \cdot)_\lambda$ on $H_1(A, \A^p_f)$ pulls-back via $\bar{\eta}$ to the pairing $a(\cdot, \cdot)$ on $V \otimes \A^p_f$; see \cite{H3} $\S 2$.     

Let $\overline{\Z}_p \subset \overline{\Q}_p$ denote the subring of elements integral over $\Z_p$, and fix $a \in \pi_0$.  Let $Sh^0$ denote the moduli space (over e.g. $\overline{\Z}_p$) of points $(A_\bullet,\lambda, i, \bar{\eta}) \in Sh$ such that 
$$
\bar{\eta}^* \, (\cdot, \cdot)_\lambda = a (\cdot, \cdot).
$$

\begin{lemma}  \label{connectedness_of_fibers}
Suppose $K_p$ is any standard parahoric subgroup.  Then the fibers of $Sh^0 \rightarrow {\rm Spec}(\overline{\Z}_p)$ are connected.  Furthermore, the morphism $\psi: \widetilde{Sh^0} \rightarrow \Mloc_{K_p}$ is surjective.
\end{lemma}

\begin{proof} 
By \cite{Del}, the generic fiber $Sh^0_{\overline{\Q}_p}$ is connected.  In the ``fake'' unitary case, $Sh^0 \rightarrow {\rm Spec}(\overline{\Z}_p)$ is proper and flat, and hence by the Zariski connectedness principle (comp. \cite{Ha}, Ex. III.11.4), the special fiber is also connected.  In the Siegel case $Sh^0 \rightarrow {\rm Spec}(\overline{\Z}_p)$ is flat but not proper, so the same argument does not apply (but see \cite{CF} IV.5.10 when $K_p$ is maximal hyperspecial).  However, the connectedness of the special fiber still holds and can be proved in an indirect way from the $p$-adic monodromy theorem of \cite{CF}; for details see \cite{Yu}.

The statement regarding surjectivity follows from the proof of Lemma \ref{Genestier_surjectivity}: the local model diagram does not ``see''  $\bar{\eta}$, and so we can arrange matters so that all constructions occur within $Sh^0$ and $\widetilde{Sh^0}$.  
\end{proof}

From now on we assume $K_p = K^{\bf a}_p$.  The above lemma proves that in any connected component $Sh^0$, all KR-strata $Sh^0_w$ are nonempty.

In fact, because the KR-stratum $Sh_\tau$ is zero-dimensional, the nonemptiness of $Sh^0_\tau \cap {\mathcal S}_{[\tau]}$  proves slightly more.  The following statement is in some sense the opposite extreme of the result of Genestier-Ng\^{o} in Corollary \ref{GN_ordinary_locus}.

\begin{cor}
In the ``fake'' unitary or the Siegel case with $K_p = K^{\bf a}_p$, let $Sh^0_\tau$ denote the the zero-dimensional KR-stratum in a connected component $Sh^0$.  Then $Sh^0_\tau$ is nonempty and is contained in the basic locus of $Sh_{K^{\bf a}_p}$.
\end{cor}

How can we describe the irreducible components in $Sh^0_{\overline{\F}_p}$?  These are clearly just the closures of the irreducible components of the KR-strata $Sh^0_{t_\lambda}$, as $\lambda$ ranges over the Weyl-orbit $W\mu$.  A priori, each of these maximal KR-strata might be the (disjoint) union of several irreducible components, all having the same dimension. 

\begin{cor}  In the Siegel or ``fake'' unitary case, $Sh^0_{\overline{\F}_p}$ is equidimensional, and the number of irreducible components is at least $\# W\mu$.
\end{cor}

In the Siegel case, a much more precise statement has been established by C.-F. Yu \cite{Yu}, answering in the affirmative a question raised in \cite{deJ}.

\begin{theorem} [\cite{Yu}]  
In the Siegel case with $K_p = K^{\bf a}_p$, each maximal $KR$-stratum $Sh^0_{t_\lambda}$ is irreducible.  Hence $Sh^0_{\overline{\F}_p}$ has exactly $2^n$ irreducible components.  An analogous statement holds for any standard parahoric subgroup $K_p$.
\end{theorem}

It is reasonable to expect that similar methods will apply to the ``fake'' unitary case to prove that the number of irreducible components in $Sh^0_{\overline{\F}_p}$ is exactly $\# W\mu$.  In fact, it would be interesting to determine whether this last statement remains true for any PEL Shimura variety attached to a group whose $p$-adic completion is unramified.

\section{Appendix: Summary of Dieudonn\'{e} theory and de Rham and crystalline cohomology for abelian varieties} \label{appendix}

This summary is extracted from some standard references -- \cite{BBM}, \cite{Dem}, \cite{Fon}, \cite{Il}, \cite{MaMe}, \cite{Me}, and \cite{O} -- as well as from \cite{deJ}.

\subsection{de Rham cohomology and the Hodge filtration}

To an abelian scheme $a: A \rightarrow S$ of relative dimension $g$ is associated the de Rham complex $\Omega^\bullet_{A/S}$ of $\cO_A$-modules.  We define the de Rham cohomology sheaves 
$$
R^ia_*(\Omega^\bullet_{A/S}).
$$
The first de Rham cohomology sheaf 
$$
R^1a_*(\Omega^\bullet_{A/S})
$$
is a locally free ${\mathcal O}_S$-module of rank $2g$.  If $S$ is the spectrum of a Noetherian ring $R$, then 
$$
H^1_{DR}(A/S) := H^1(A, \Omega^\bullet_{A/S}) = \Gamma(S,R^1a_*(\Omega^\bullet_{A/S}))
$$
is a locally free $R$-module of rank $2g$.

The Hodge-de Rham spectral sequence degenerates at $E_1$ (\cite{BBM}, $\S 2.5$), yielding the exact sequence
$$
0 \rightarrow R^0a_*(\Omega^1_{A/S}) \rightarrow  R^1a_*(\Omega^\bullet_{A/S}) \rightarrow R^1a_*(\cO_A) \rightarrow 0.
$$
We define $\omega_A := R^0a_*(\Omega^1_{A/S})$, a locally free sub-$\cO_S$-module of rank $g$.  The term $R^1a_*(\cO_A)$ may be identified with ${\rm Lie}(\widehat{A})$, the Lie algebra of the dual abelian scheme $\widehat{A}/S$.  It is also locally free of rank $g$.  Thus we have the {\em Hodge filtration} on de Rham cohomology
$$ 
0 \rightarrow \omega_A \rightarrow  R^1a_*(\Omega^\bullet_{A/S}) \rightarrow {\rm Lie}(\widehat{A}) \rightarrow 0.
$$

Recall that in our formulation of the moduli problem defining ${\rm Sh}({\bf G}, h)_{K_p}$, the important determinant condition refers to the Lie algebra ${\rm Lie}(A)$, and not to ${\rm Lie}(\widehat{A})$.  Because of this it is convenient (although not, strictly-speaking, necessary) to work with the {\em covariant} analogue $M(A)$ of 
$R^1a_*(\Omega^\bullet_{A/S})$. To define it, recall that for any $\cO_S$-module $N$, we define the dual $\cO_S$-module $N^\vee$ by
$$
N^\vee := {\mathcal Hom}_{\cO_S}(N, \cO_S).
$$
Let $M(A) := (R^1a_*(\Omega^\bullet_{A/S}))^\vee$ be the dual of de Rham cohomology.
This is a locally free $\cO_S$-module of rank $2g$.  By the proposition below, we can identify $\omega^\vee_{A} = {\rm Lie}(A)$ and so the {\em Hodge filtration} on $M(A)$ takes the form
$$
0 \rightarrow {\rm Lie}(\widehat{A})^\vee \rightarrow M(A) \rightarrow {\rm Lie}(A) \rightarrow 0.
$$

It is sometimes convenient to denote ${\mathbb D}(A)_S := R^1a_*(\Omega_{A/S})$  (this notation refers to crystalline cohomology, see \cite{BBM}, \cite{Il}).

\begin{prop} [\cite{BBM}, Prop. 5.1.10]  \label{BBM_5.1.10}  There is a commutative diagram whose vertical arrows are isomorphisms 
$$
\xymatrix{
0 \ar[r] & {\rm Lie}(\widehat{A})^\vee \ar[r] \ar[d]_{\cong} & {\mathbb D}(A)_S^\vee \ar[r] \ar[d]_{\cong} & 
\omega_A^\vee \ar[r] \ar[d]_{\cong} & 0 \\ 0 \ar[r] & \omega_{\widehat{A}} \ar[r] & {\mathbb D}(\widehat{A})_S 
\ar[r] & {\rm Lie}(A) \ar[r] & 0.}
$$
\end{prop}

\subsection{Crystalline cohomology}

Let $k_r = \mathbb F_{p^r}$ be a finite field with ring of Witt vectors $W(k_r)$.  The fraction field $L_r$ of $W(k_r)$ is an unramified extension of ${\mathbb Q}_p$ and its Galois group is the cyclic group of order $r$ generated by the Frobenius  element $\sigma: x \mapsto x^p$; note also that $\sigma$ acts on Witt vectors by the rule $\sigma(a_0,a_1, \dots) = (a_0^p, a_1^p, \dots)$.  

Let $A$ be an abelian variety over $k_r$ of dimension $g$.  We have the {\em integral isocrystal} associated to $A/k_r$, given by the data 
$$
{\mathbb D}(A) = (H^1_{\rm crys}(A/W(k_r)), F, V).
$$

Here the crystalline cohomology group $H^1_{\rm crys}(A/W(k_r))$ is a free $W(k_r)$-module of rank $2g$, equipped with a $\sigma$-linear endomorphism $F$ (``Frobenius'')  and the $\sigma^{-1}$-linear endomorphism $V$ (``Verschiebung'')  which induce bijections on $H^1_{\rm crys}(A/W(k_r)) \otimes_{W(k_r)} L_r$.  We have the identity $FV = VF = p$ (by {\em definition} of $V$), hence the inclusions of $W(k_r)$-lattices 
$$
H^1_{\rm crys}(A/W(k_r)) \supset FH^1_{\rm crys}(A/W(k_r)) \supset p H^1_{\rm crys}(A/W(k_r))
$$
(as well as the analogous inclusions for $V$ replacing $F$).

The endomorphism $F$ has the property that $F^r = \pi_A$ on $H^1_{\rm crys}$, where $\pi_A$ denotes the absolute Frobenius morphism of $A$ relative to the field of definition $k_r$ (on projective coordinates $x_i$ for $A$, $\pi_A$ induces the map $x_i \mapsto x_i^{p^r}$).

\subsection{Relation with Dieudonn\'{e} theory}

The crystalline cohomology of $A/k_r$ is intimately connected to the (contravariant) Dieudonn\'{e} module of the $p$-divisible group $A[p^\infty] := \varinjlim A[p^n]$, the union of the sub-groupschemes $A[p^n] = {\rm ker}(p^n: A \rightarrow A)$.   Recall that the classical contravariant Dieudonn\'{e} functor $G \mapsto D(G)$ establishes an exact anti-equivalence between the categories
$$
\left\{ \mbox{$p$-divisible groups $G = \varinjlim G_n$ over $k_r$} \right\}
$$
and
$$ 
\{ 
\mbox{free $W(k_r)$-modules $M = \varprojlim M/p^nM$, equipped with operators $F,V$} \};
$$
see \cite{Dem}, \cite{Fon}.   Here $F$ and $V$ are $\sigma$- resp. $\sigma^{-1}$-linear, inducing bijections on $M \otimes_{W(k_r)} L_r$.

The crystalline cohomology of $A/k_r$, together with the operators $F$ and $V$, is the same as the Dieudonn\'{e} module of the $p$-divisible group $A[p^\infty]$, in the sense that there is a canonical isomorphism  
\begin{equation} \label{eq:crys=Dieudonne}
H^1_{\rm crys}(A/W(k_r)) \cong D(A[p^\infty])
\end{equation}
which respects the endomorphisms $F$ and $V$ on both sides, cf. \cite{BBM}.  Moreover, we have the following identifications
\begin{equation} \label{eq:crys=DR=Dieudonne}
{\mathbb D}(A)_{k_r} := H^1_{\rm crys}(A/W(k_r)) \otimes_{W(k_r)} k_r \cong H^1_{DR}(A/k_r) \cong D(A[p]).
\end{equation}
The second isomorphism is due to Oda \cite{O}; see below.  The first isomorphism is a standard fact (\cite{BBM}), but can also be deduced via Oda's theorem by reducing equation (\ref{eq:crys=Dieudonne}) modulo $p$: the exactness of the functor $D$ implies that $D(A[p]) = D(A[p^\infty])/pD(A[p^\infty]) = D(A[p^\infty]) \otimes_{W(k_r)} k_r$.  In particular, the $k_r$-vector space $H^1_{DR}(A/k_r)$ inherits endomorphisms $F$ and $V$ ($\sigma$- resp. $\sigma^{-1}$-linear).  

The theorem of Oda \cite{O} includes as well the relation between the Hodge filtration on the de Rham cohomology of $A$ and a suitable filtration on the isocrystal ${\mathbb D}(A)$.

\begin{theorem} [\cite{O}, Cor. 5.11] \label{Oda}
There is a natural isomorphism $\psi: {\mathbb D}(A)_{k_r} \tilde{\rightarrow} H^1_{DR}(A/k_r)$, and under this isomorphism, $V{\mathbb D}(A)_{k_r}$ is taken to $\omega_A$.  In particular there is an exact sequence
$$
0 \rightarrow V {\mathbb D}(A) \rightarrow {\mathbb D}(A) \rightarrow {\rm Lie}(\widehat{A}) \rightarrow 0.
$$  
\end{theorem}

\subsection{Remarks on duality}

We actually make use of this in a dual formulation.  Define $H = H(A)$ to be the $W(k_r)$-linear dual of the isocrystal ${\mathbb D}(A)$
\begin{equation} \label{eq:dual_isocrys}
H(A) = {\rm Hom}_{W(k_r)}( {\mathbb D}(A), W(k_r)),
\end{equation}
and define $H_{L_r} = H(A)_{L_r} = H(A)  \otimes_{W(k_r)} L_r$.  Letting $\langle \,\, , \,\, \rangle :  H \times {\mathbb D}(A) \rightarrow W(k_r)$ denote the canonical pairing, we define $\sigma$- resp. $\sigma^{-1}$-linear injections $F$ resp. $V$ on $H$ (they are bijective on $H_{L_r}$)  by the formulae
\begin{align} \label{eq:dual_FV}
\langle Fu, a \rangle &= \sigma \langle u, Va \rangle \\
\langle Vu, a \rangle &= \sigma^{-1} \langle u, Fa \rangle 
\end{align}
for $u \in H$ and $a \in {\mathbb D}(A) = H^1_{\rm crys}(A/W(k_r))$.  

Of course $H(A)_{k_r} := H \otimes_{W(k_r)} k_r$ is the $k_r$-linear dual of ${\mathbb D}(A)_{k_r}$, hence
$$
H(A)_{k_r} = M(A/k_r).
$$
    
\begin{lemma}  Let $S = {\rm Spec}(k_r)$.  Equip $H = {\mathbb D}(A)^\vee$ with operators $F,V$ as in 
(\ref{eq:dual_FV}).  Then the isomorphism
$$
{\mathbb D}(A)_{k_r}^\vee \,\, \widetilde{\rightarrow} \,\, {\mathbb D}(\widehat{A})_{k_r}
$$
of Proposition \ref{BBM_5.1.10} is an isomorphism of $W(k_r)[F,V]$-modules.
\end{lemma}

\begin{proof}
From \cite{Dem}, Theorem 8 (p. 71), for a p-divisible group $G$ with Serre dual $G'$ (loc. cit., p. 46) there is a duality pairing in the category of $W(k_r)[F,V]$-modules
$$
\langle \,\, , \,\, \rangle: D(G') \times D(G) \rightarrow W(k_r)(-1)
$$
where $W(k_r)(-1)$ is the isocrystal with underlying space $W(k_r)$ and $\sigma$-linear endomorphism $p \sigma$.  
That is, we have the identity
$$
\langle F_{G'} x, F_{G} y \rangle = p \sigma \langle x, y \rangle,
$$
or 
$$
\langle F_{G'} x, y \rangle = \sigma \langle x, pF_{G}^{-1} y \rangle = \sigma \langle x, V_{G}y \rangle.
$$

There is a canonical identification 
$$
(A[p^\infty])' = \widehat{A}[p^\infty].
$$
Now from the above pairing with $G = A[p^\infty]$ and (\ref{eq:crys=Dieudonne}) we deduce a duality pairing of $W(k_r)[F,V]$-modules
$$
{\mathbb D}(\widehat{A}) \times {\mathbb D}(A) \rightarrow W(k_r)(-1),
$$
which induces the isomorphism of Proposition \ref{BBM_5.1.10}.  The lemma follows from these remarks.
\end{proof}    

Applying Oda's theorem (\ref{Oda}) to $\widehat{A}$ and invoking the above lemma gives ${\rm Lie}(\widehat{A})^\vee = VM(A/k_r)$, thus there is an exact sequence
$$
0 \rightarrow V \, H(A) \rightarrow H(A) \rightarrow {\rm Lie}(A) \rightarrow 0.
$$

Thus, we have 

\begin{cor}\label{sigma(Lie)} 
Let $F,V$ be the $\sigma$- resp. $\sigma^{-1}$-linear endomorphisms of $H = H(A)$ defined in (\ref{eq:dual_FV}).  There is a natural isomorphism
$$
\dfrac{V^{-1}\, H}{H} = \sigma({\rm Lie}(A)).
$$
Moreover, on $H_{L_r}$ we have the identity $V^{-r} = \pi^{-1}_A$.
(Comp. \cite{Ko92}, $\S 10, 16$.)
\end{cor}

\end{document}